\documentclass[11pt]{llncs}
\usepackage{graphicx}

\usepackage{graphicx}
\usepackage{amsmath}
\usepackage{amssymb}
\usepackage{mathrsfs}
\usepackage[T1]{fontenc}
\usepackage{varioref}
\usepackage{calc}
\usepackage{ifthen}
\usepackage{enumerate}
\usepackage{hyperref}
\usepackage{tikz}
\usetikzlibrary{calc,through,backgrounds,arrows}
\usepackage{color}
\newcommand{\F}{\mathbb{F}}
\newcommand{\G}{\mathbb{G}}
\newcommand{\DN}{\mathbb{S}}
\newcommand{\TDN}{\mathbb{T}}

\newcommand{\N}{\mathbb{N}}
\newcommand{\Z}{\mathbb{Z}}
\newcommand{\abs}[1]{\lvert#1\rvert}

\newcommand{\mset}[1]{\mathcal{#1}}

\newcommand{\teich}{\mathrm{T}}

\newcommand{\eye}{\ensuremath{\!\!\begin{array}{c}\frown\\[-12pt]\smile\end{array}\!\!}}

\newcommand{\points}{\mathcal{P}}
\newcommand{\lines}{\mathcal{L}}

\newcommand{\incidence}{I}

\newcommand{\pclasses}{\overline{\mathcal{P}}}
\newcommand{\lclasses}{\overline{\mathcal{L}}}
\newcommand{\iclasses}{\overline{\incidence}}
\newcommand{\psegments}{\mathfrak{P}}
\newcommand{\lsegments}{\mathfrak{L}}
\newcommand{\isegments}{\mathfrak{I}}

\newcommand{\maxsize}{\mathrm{m}}

\DeclareMathOperator{\PG}{PG}
\DeclareMathOperator{\AG}{AG}
\DeclareMathOperator{\PHG}{PHG}

\DeclareMathOperator{\Aut}{Aut}

\DeclareMathOperator{\GR}{GR}

\DeclareMathOperator{\rad}{Rad}
\DeclareMathOperator{\Char}{char}
\DeclareMathOperator{\supp}{supp}

\newcommand{\opp}{\circ}

\newtheorem{fact}{Fact}
\date{}

\begin{document}

\title{New Results on Arcs in Projective Hjelmslev Planes \\
	over Small Chain Rings}

\author{Thomas Honold\inst{1}, Michael Kiermaier\inst{2}, Ivan Landjev\inst{3,4}}

\institute{
	\emph{Zhejiang Provincial Key Laboratory of Information Network Technology}\\ 
	\emph{and Department of Information and Electronic Engineering,}\\ 
	\emph{Zhejiang University, 38 Zheda Road, 310027 Hangzhou, China} \\
	\email:{honold@zju.edu.cn}) \and
		\emph{Lehrstuhl Mathematik II, Universit\"at Bayreuth,}\\ 
	\emph{95447 Bayreuth, Germany}
	\email{michael.kiermaier@uni-bayreuth.de} \and	
	\emph{Bulgarian Academy of Sciences, Institute of Mathematics and Informatics,} \\
	\emph{8 Acad G. Bonchev str., 1113 Sofia, Bulgaria}\\
	\email{ivan@math.bas.bg} \and
	\emph{New Bulgarian University, 21 Montevideo str., 1618 Sofia,
		Bulgaria} \\
	\email{i.landjev@nbu.bg}
}
	
\maketitle

\begin{abstract}
  We present various new constructions and bounds for arcs in
  projective Hjelmslev planes over finite chain rings of nilpotency
  index $2$. For the chain rings of cardinality at most $25$ we give updated
  tables with the best known upper and lower bounds for the maximum size of
  such arcs.
  
\end{abstract}

%================================================================

\section{Introduction}\label{sec:intro}

Let $R$ be a finite chain ring,
and let $\PHG(2,R)$ denote the right projective Hjelmslev plane over $R$
\cite{klingenberg54a,klingenberg55,kleinfeld59,artmann69,drake70,alex91a,it:chain}.
The maximal arc problem for coordinate projective Hjelmslev planes can be formulated
as follows: Given $R$ and an integer $n\geq 0$, find the largest size
$k$ of a $(k,n)$-arc in $\PHG(2,R)$.

The problem of finding good arcs is related to the construction of good linear codes over $R$.
These are linear codes over $R$ which, by means of a suitable Gray map, correspond to good
$q$-ary, not necessarily linear, codes over an alphabet of size which is a prime power.
The distance related to this problem is the usual Hamming distance.
An illustrative example is provided by the hyperovals in $\PHG(2,\mathbb{Z}_4)$ which correspond to the optimal non-linear binary $(14,2^6,6)$-codes.
The problem of finding the size of largest $(k,n)$-arc for fixed $n$ 
was studied extensively in the classical projective planes $\PG(2,q)$
(cf. for instance \cite{ball-tables,mbraun19} and the references there).

In \cite{it:kent,it:deadfin} the maximal arc problem was studied for
finite chain rings of length (nilpotency index) two.
Tables of exact values and bounds
for arcs in the Hjelmslev planes over small chain rings $R$,
($|R|=4,9,16,25$) were given. In this paper, we present new
constructions and bounds for such arcs, which improve considerably upon
\cite{it:deadfin,it:kent}. Parts of the work presented in this
paper have been announced in
\cite{ivan-silvia04,mt:phgarcs,kiermaier-kohnert07,kiermaier-koch09},
respectively, are reflected in the online tables \cite{arctables}
maintained by the second author.

The paper is organized as follows. In Section~2 we introduce the
necessary terminology about chain rings and projective Hjelmslev planes.
In Section~3 we present various constructions for
arcs. In Section~4 we present improvements on the known
upper bounds including several nonexistence proofs for important sets
of parameters. Section~5 contains the improved tables as well as some open
problems that we consider to be of interest.

\section{Preliminaries and earlier results}
\label{sec:preliminaries}

Everywhere in this paper, $R$ will denote a chain ring of nilpotency
index two with residue field $R/\rad R\cong\F_q$, where $q=p^r$
is a prime power. All such rings have cardinality $q^2$. For fixed $q$
there are $r+1$ isomorphism classes of such rings
\cite{cronheim78,raghavendran69}, the Galois ring
$\G_q=\Z_{p^2}[X]/(h)$ of characteristic $p^2$, where
  $h\in\Z_{p^2}[X]$ is a monic polynomial of degree $r$ which is
  irreducible modulo $p$,  and $r$ truncated (skew) polynomial rings
$\TDN_q^{(i)}=\F_q[X;\sigma_i]/(X^2)$, $\sigma_i\in\Aut(\F_q)$, of
characteristic $p$.
The chain rings of length $2$ and cardinality
$q^2\leq 25$ are
\begin{tabbing}
$q^2=4$:\qquad\=$\G_2=\Z_4$,\quad $\DN_2=\F_2[X]/(X^2)$;\\
$q^2=9$:\>$\G_3=\Z_9$,\quad $\DN_3=\F_3[X]/(X^2)$;\\
$q^2=16$:\>$\G_4=\Z_4[X]/(X^2+X+1)$,\quad
$\DN_4=\F_4[X]/(X^2)$,\\\>$\TDN_4=\F_4[X;a\mapsto a^2]/(X^2)$;\\
$q^2=25$:\>$\G_5=\Z_{25}$,\quad $\DN_5=\F_5[X]/(X^2)$.
\end{tabbing}

Denote by $R^3_R$ the free right $R$-module of rank $3$ and by
$\points$ (resp.\ $\lines$) the set of free rank $1$ (resp.\ rank $2$)
submodules of $R^3_R$.  The right projective Hjelmslev plane over $R$,
denoted by $\PHG(2,R)$, is defined as the incidence structure
$(\points,\lines,\subseteq)$ with point set $\points$, line set
$\lines$ and incidence relation $\subseteq$ (set inclusion). The left
projective Hjelmslev plane over $R$, defined in an analogous manner
using submodules of ${}_RR^3$, is isomorphic to the right projective
Hjelmslev plane $\PHG(2,R^\opp)$ over the opposite ring
$R^\opp$ with  multiplication operation $a\circ b= ba$.  
In the particular case of chain rings $R$ with $\abs{R}\leq 25$ (the
case considered here) we always
have $R\cong R^\opp$ and hence also $\PHG(2,R)\cong\PHG(2,R^\opp)$,
i.e.\ we need not distinguish between left and right projective
Hjelmslev planes. 

For standard facts about the planes $\PHG(2,R)$, many of which will not be
repeated here, we refer to \cite{it:hamburg,it:chain,it:kent,it:deadfin}.

We will use the standard terminology for multisets in point-line incidence
structures $\Pi=(\points,\lines,\incidence)$ as developed, for
example, in \cite{it:hamburg,it:chain}. A \emph{multiset in $\Pi$} is
a mapping $\mset{K}\colon\points\to\N_0$, assigning to every point
$x\in\points$ a nonnegative integer $\mset{K}(x)$ called the \emph{multiplicity}
of $x$ in $\mset{K}$. The definition of multiplicity extends to
the power set of $\points$ by setting $\mset{K}(S)=\sum_{x\in S}\mset{K}(x)$
for $S\subseteq\points$. Identifying lines
with sets of points (as will always be possible in the cases under
consideration), this defines, for every line $L\in\lines$, a
corresponding multiplicity $\mset{K}(L)$. The multiplicity
$\abs{\mset{K}}=\mset{K}(\points)$ of the whole point set $\points$ is
called the \emph{cardinality} (or \emph{size}) of $\mset{K}$. 
An important invariant of $\mset{K}$ is its \emph{spectrum},
defined as the sequence of integers
$a_i=a_i(\mset{K})=\#\bigl\{L\in\lines;\mset{K}(L)=i\bigr\}$,
$i\in\N_0$.

A multiset in $\Pi$ is said to be \emph{projective} if
$\mset{K}(x)\in\{0,1\}$ for all points $x\in\points$; i.e., the
multiset $\mset{K}$ is equal to the characteristic function of its
\emph{support} and thus may be considered as a set.
A multiset in $\Pi$ is called a \emph{$(k,n)$-arc}
(resp. \emph{$(k,n)$-blocking multiset}) if
$\abs{\mset{K}}=k$ and $\mathfrak{K}(L)\le n$ 
for every line $L\in\mathcal{L}$ (resp. $\mathfrak{K}(L)\ge n$
for every $L\in\mathcal{L}$). 

A $(k,n)$-arc in $\Pi$ is said to be \emph{maximal} if
$k'\leq k$ for all $(k',n)$-arcs in $\Pi$, and
\emph{complete} if increasing the multiplicity in $\mset{K}$ of any
point $x\in\points$ by one results in a line of multiplicity $n+1$.
Blocking multisets with the analogous properties are referred to as
\emph{minimal}, resp.\ \emph{irreducible}.

If $\mathcal{K}$ is a 
projective $(k,n)$-arc in $\Pi$ then
$1-\mset{K}$  is a projective blocking multiset with
parameters $(|\mathcal{P}|-k,l-n)$, where $l$ is the size of a line
in $\Pi$.
%Thus projective arcs and blocking sets are essentially
%equivalent concepts, the arc point-of-view being more convenient if
%$n$ is small and less convenient if it close to the line size.
%than the blocking set point-of-view if $n$ is small compared with $l$.

The present work is devoted to the maximal arc problem in projective
Hjelmslev planes. Following \cite{it:kent,it:deadfin}
we restrict ourselves to projective arcs
and accordingly define, for a finite chain ring $R$ of length $2$ and an
integer $n\in\{0,1,\dots,q^2+q\}$ the number
$\maxsize_n(R)$ as the largest size $k$ for which
a projective $(k,n)$-arc in $\PHG(2,R)$ exists.
  
  It is reasonable to allow for non-projective arcs in the
  definition of $\maxsize_n(R)$ (and thereby extend the domain of
  $n\mapsto\maxsize_n(R)$ to all integers $n\geq 0$). In this paper, we have
  refrained from doing so. Changing the definition of $\maxsize_n(R)$
  might lead to inconsistencies with previously published
  work, since (at least at present) we do not know whether for
  all integers $n$ in the range $0\leq n\leq q^2+q$ non-projective
  $(k,n)$-arcs necessarily satisfy the inequality
  $k\leq\maxsize_n(R)$.

  In subsequent sections we will give upper and lower bounds for
$\maxsize_n(R)$, 
strengthening considerably those of
the earlier work \cite{it:kent,it:deadfin}, and we will determine the
exact values of $\maxsize_n(R)$ in several new cases. As in
\cite{it:kent,it:deadfin}, our analysis relies on
three important structure theorems for $\PHG(2,R)$, which relate
the maximal arc problem for $\PHG(2,R)$ to the corresponding well-studied
problem for the classical planes $\PG(2,q)$ and
$\AG(2,q)$. We denote the largest sizes of
$(k,n)$-arcs in $\PG(2,q)$ and $\AG(2,q)$ 
by $\maxsize_n(q)$ resp.\ $\maxsize'_n(q)$. Detailed information about these  
numbers, especially for small values of $q$, is available in
\cite[Ch.~12--14]{hirschfeld98}, as well as in \cite{ball-tables,mbraun19}.

We will now review the aforementioned
three structural properties of $\PHG(2,R)$, along with introducing
some further convenient notation. Proofs of the structure theorems may
be found in the original literature on Hjelmslev geometry cited at the
beginning of Section~\ref{sec:intro} and also in \cite{it:deadfin}.

In $\Pi=\PHG(2,R)=(\points,\lines,\subseteq)$ there is defined a
neighbour relation $\eye$ on points and lines by $x\eye x'$ for
$x,x'\in\points$ (resp.\ $L\eye L'$ for $L,L'\in\lines$) iff the
images in $\F_q^3$ (under the canonical projection)
of the corresponding $R$-submodules of $R^3_R$ coincide.
Neighbourhood $\eye$ defines an equivalence relation
on both $\points$ and $\lines$, and the quotient incidence structure
$\overline{\Pi}=(\pclasses,\lclasses,\iclasses)$ formed by the point
and line classes of $\PHG(2,R)$ is isomorphic to
$\PG(2,q)$.\footnote{This ``first structure theorem'' is not really a
  theorem, but rather part of the definition of a projective Hjelmslev
  plane.}
These classes are denoted as usual by $[x]$ resp.\ $[L]$.

Intersections $L\cap[x]\neq\emptyset$ of lines with point neighbour
classes are called \emph{line segments}. The points in a fixed point class
$[x]$ and the line segments contained in $[x]$ form an incidence
structure $\Pi_{[x]}\cong\AG(2,q)$.  

Given $[x]$ and a line segment
$S\subseteq[x]$, the lines $L$ for which $L\cap[x]$ is parallel to $S$
form a neighbour class of lines called the \emph{direction} of
$S$. The line segments with a fixed direction $[L]$ and the lines in
$[L]$ form an incidence structure
$\Pi_{[L]}=(\psegments,\lsegments,\isegments)$, which is
isomorphic to a punctured projective plane $\PG(2,q)\setminus
p_\infty$. The $q+1$ point neighbour classes $[x]$ incident with $[L]$
form the corresponding transversals (considered as sets of line
segments), so that $(\psegments,\lsegments,\isegments)$ may be
completed to a copy of $\PG(2,q)$ by adding the ``infinite'' point
$p_\infty$ to each such class $[x]$.

Now suppose $\mset{K}$ is a set of points in $\PHG(2,R)$ and
$[x]\in\pclasses$ is a point class. If $S$ is a line segment in $[x]$ joining
two distinct points of $\mset{K}$, we say that $S$ (or the direction
$[L]$ of $S$) is \emph{determined} by the points of $\mset{K}$ in $[x]$.

The point set $\mset{K}$ assigns to
each parallel class $\{S_1,\dots,S_q\}$ of lines of the affine plane
$\Pi_{[x]}$ a corresponding sequence
$\mset{K}(S_1),\dots,\mset{K}(S_q)$ of multiplicities called the \emph{type
of $[x]$ (with respect to $\mset{K}$)}.\footnote{To make types unique,
  we arrange the integers $\mset{K}(S_i)$ in non-increasing order.}
Since each parallel class $\{S_1,\dots,S_q\}$ is associated with a
unique direction $[L]$, types may also be considered as being defined
for flags of the quotient plane $(\pclasses,\lclasses,\iclasses)$.
Further it will be convenient to arrange the multiplicities of
the $q(q+1)$ line segments with a fixed direction $[L]$ (the
``finite'' points of $\Pi_{[L]}$) into a
$q\times(q+1)$ array $(k_{ij})$ whose columns are the types of the
various classes $[x_j]$ incident with $[L]$.\footnote{Here $[x_j]$ is again
  considered as a parallel class of line segments.} This array will
be referred to as the \emph{type of $[L]$ (with respect to
  $\mset{K}$)}.\footnote{We leave the order of the columns unspecified
  (although uniqueness of the type of $[L]$ could be easily enforced by
  ordering the columns appropriately).} Its entries can also be seen as
the values of the multiset in $\Pi_{[L]}$ induced by $\mset{K}$.

In \cite{it:kent} the following general upper bound for the size of
$(k,n)$-arcs was proved.\footnote{This bound is valid also for
  non-projective arcs.}

\begin{fact}[{\cite[Thm.~3.1]{it:kent}}]
\label{thm:bound}
Let $\mathfrak{K}$ be a $(k,n)$-arc in $\PHG(R_R^3)$ and $[x]$ a
neighbour class 
of points in $\PHG(R_R^3)$. Let further
$u=\mset{K}([x])$ and $u_i$, $i=1,2,\ldots,q+1$,
the maximum number of points on a line from the $i$-th parallel
class in the affine plane defined on $[x]$. Then
\[k \le q(q+1)n - q\sum_{i=1}^{q+1} u_i + u.\]
\end{fact}
Using the bounds $\sum_{i=1}^{q+1}u_i\geq q+u$ (the sum of the
multiplicities of all lines through a point of $\mset{K}$ in the
affine plane $\Pi_{[x]}$) and $u_i\geq\lceil u/q\rceil$, this
theorem has the following corollary:
% which in light of \
% Rem.~\ref{rmk:projective} is a slight extension of
% \cite[Cor.~3.3]{it:kent}.
\begin{fact}[{\cite[Cor.~1]{it:deadfin} or \cite[Cor.~3.3]{it:kent}}]
%\begin{corollary}
\label{cor:bound}
For $1\leq n\leq q^2+q$ we have
\begin{multline*}
  \label{eq:genineq}
\maxsize_n(R)\le\max_{1\le u\le q^2}\min\{u(q^2+q+1),\;
q(q+1)n-q^2-(q-1)u,\\
q(q+1)n-q(q+1)\lceil u/q\rceil+u\}. 
\end{multline*}
%\end{corollary}
\end{fact}
% \noindent {\it Proof.}
%   Suppose that $\mset{K}$ is a $(k,n)$-arc in $\PHG(R_R^3)$. We apply
%   Fact~\ref{thm:bound} to a point $x$ with $\mset{K}(x)>0$ and
%   $\mset{K}([x])=\max\bigl\{\mset{K}([y]);y\in\points\bigr\}$. The
%   sum $\sum_{i=1}^{q+1}u_i$ is lower bounded
%   both by $q\cdot\mset{K}(x)+u$ (the sum of the
%   multiplicities of all lines through $x$ in the affine plane defined
%   on $[x]$) and by $(q+1)\lceil u/q\rceil$ (using
%   $u_i\geq\lceil u/q\rceil$). In the first case we obtain
%   \begin{equation}
%     \label{eq:genineq1}
%     k\leq q(q+1)n-q^2\cdot\mset{K}(x)-(q-1)u,
%   \end{equation}
% which in view of $\mset{K}(x)\geq 1$
% leads to the stated upper bound for $k$, provided that
% $u\leq q^2$. In the case $u=q^2$ the minimum in \eqref{eq:genineq} is
% $q(q+1)n-q^3$ and is attained by the
% second and third expression in \eqref{eq:genineq} (here we need $n\leq
% q^2+q$). Since the second expression, $q(q+1)n-q^2-(q-1)u$, is a
% nonincreasing function of $u$, the stated upper bound for $k$ is also
%   valid in the case $u>q^2$.
% \qed
% \begin{remark}
%   \label{rmk:projective2}
%   A more thorough analysis along the lines of the preceding proof
%   shows that arcs attaining the bound of Cor.~\ref{cor:bound} are
%   necessarily projective.\footnote{Since the bound may not be
%     attained, this does not imply that maximal $(k,n)$-arcs with
%     $n\leq q^2+q$ are projective.}
% \end{remark}

This bound can often be strengthend by using more precise estimates
for the numbers $u_i$ occuring in the general bound
(Fact~\ref{thm:bound}). For example, if $q$ is odd and $u=q+2$, then
in view of the non-existence of $(q+2,2)$-arcs in $\AG(2,q)$ we cannot
have $\sum_{i=1}^{q+1}u_i=2q+2$.  In Section~\ref{ssec:newbound} we
will give tight lower bounds on $\sum_{i=1}^{q+1}u_i$ for $q\leq 5$ (so that
Fact~\ref{thm:bound} can be used in its full strength for these values
of $q$), and we will also make the bound in Fact~\ref{cor:bound} more
explicit.

\section{Constructions}\label{sec:const}

\subsection{General Constructions}\label{ssec:genconst}

In this section we present new general constructions of arcs in
$\PHG(2,R)$. These constructions yield larger
arcs than those given in previous
papers; they are defined for arbitrary finite chain rings of length
$2$ and do not distinguish between rings of the same cardinality.

We start with constructions of (projective) $(k,n)$-arcs in $\PHG(2,R)$
for large values of $n$. 
% In that case it is more convenient to consider the complementary
% $(q^4+q^3+q^2-k,q^2+q-n)$-blocking set.
For $q^2\leq n\leq q^2+q$ the exact
value of $\maxsize_n(R)$ is known and equal to 
$q^4+q^3+q^2-(q^2+q-n)(q^2+q)=q(q+1)n-q^3$; cf.\
\cite[Cor.~2]{it:deadfin}. The result in \cite{it:deadfin}
  was stated for projective arcs. It remains true for not
  necessarily projective arcs, since the bound of Cor.~\ref{cor:bound}
  is attained in this case (with $u=q^2$).
Corresponding maximal arcs can be obtained by removing from the point
set $\points$ a minimal $\bigl(s(q^2+q),s\bigr)$ blocking set with
$s=q^2+q-n\in\{0,1,\dots,q\}$.
These blocking sets have the same parameters as the
(in general non-projective) sum of $s$ lines.

\begin{theorem}
  \label{thml:s=q+1} For every finite chain ring $R$ there exist
  $(q^4-q^2-q,q^2-1)$-arcs in
  $\PHG(2,R)$; equivalently $\maxsize_n(R)\geq q(q+1)n-q^3$ for $n=q^2-1$.
\end{theorem}

\noindent
{\it Proof.}
  The existence of such arcs is equivalent to the existence
  of blocking sets with parameters $\bigl((q+1)(q^2+q),q+1\bigr)$.
  We shall assume $q\geq 4$, since $(18,3)$-blocking sets (equivalently,
  $(10,3)$-arcs) in the planes over $\Z_4$ and $\DN_2$ as well as
  $(48,4)$-blocking sets (equivalently, $(69,8)$-arcs) in the planes
  over $\Z_9$ and $\DN_3$ are already known; see \cite{it:deadfin}
  respectively \cite{imt:two_exact_values} or
  Remark~\ref{rmkl:q=3:n=8}.

  We choose $3$ points $x_0$, $x_1$, $x_2$ and
  $3$ lines $L_0$, $L_1$, $L_2$ determining a triangle in
  $\overline{\Pi}$ (i.e. the neighbour classes $[x_i]$ and $[L_i]$ are
  distinct) and arrange indices such that $L_i\cap L_{i+1}=x_{i+2}$
  for $i=0,1,2$ (indices taken mod $3$).

  Now consider the $\bigl((q+1)(q^2+q),q+1\bigr)$-blocking multiset
  $\mset{K}=(q-2)\delta_{L_0}+2\delta_{L_1}+\delta_{L_2}$ in $\Pi$
  (where $\delta_{L_i}$ denotes the characteristic function of $L_i$).
  We transform $\mset{K}$ into a blocking set $\mset{K}'$
  with the same parameters as follows.

  In $[x_2]$, which has multiplicity $\mset{K}\bigl([x_2])=q^2$, we
  take all points as members of $\mset{K}'$. In $[x_0]$ and $[x_1]$, 
  which have multiplicities $\mset{K}\bigl([x_0])=3q$ resp.\ 
  $\mset{K}\bigl([x_1])=(q-1)q$, we take as members of $\mset{K}'$ 
  all points contained in $3$ resp.\ $q-1$ parallel line segments with
  direction $[x_0x_2]=[L_1]$ resp.\ $[x_1x_2]=[L_0]$. In the remaining
  point classes incident with some line class $[L_i]$ we replace all
  multiple line segments of $\mset{K}$
  by the appropriate number of single, parallel line segments in the
  same direction. 

  By construction, the resulting point set $\mset{K}'$ (which is not
  yet completely specified!) has cardinality $(q+1)(q^2+q)$ and blocks
  each line outside $[L_1]\cup[L_2]\cup[L_3]$ exactly $q+1$
  times. Lines in $[L_2]$ are blocked at least $q+2$ times (the number
  of line segments in $[x_0]\cup[x_1]$). The proof is completed by
  arranging the $(q-2)q+1\geq 2q+1$ line segments in $[L_0]\setminus[x_2]$ and
  the $2q+1$ line segments in $[L_1]\setminus[x_2]$ in such a way
  that they form blocking sets in the corresponding affine planes. 
  This can clearly be done.\footnote{In
    $[L_0]\setminus[x_2]$, for example, we can take all line
    segments contained in two non-parallel lines first and
    then add appropriate further line segments.}

\begin{figure}[htp!]
	\begin{center}
\begin{tikzpicture}[line width=1pt, scale=0.4]
		\draw[gray] (0,0)--(3,0)--(3,3)--(0,3)--(0,0) [fill=yellow!30];
		\draw[gray] (16,0)--(19,0)--(19,3)--(16,3)--(16,0) [fill=yellow!30];
		\draw[gray] (8,10)--(11,10)--(11,13)--(8,13)--(8,10) [fill=yellow!30];

\draw[gray] (0,3) .. controls (0.5,6.5) and (4,11.5) .. (8,13)--(8,10) .. controls (5.5,8) and (3.5,5.5) .. (3,3)--(0,3) [fill=blue!20];
%\draw[gray] (3,3) .. controls (3.5,5.5) and (5.5,8) .. (8,10);
\draw[gray] (16,3) .. controls (15.5,5.5) and (13.5,8) .. (11,10)--(11,13) .. controls (15,11.5) and (18.5,5.5) .. (19,3) [fill=blue!20];
\draw[gray] (3,0)--(16,0)--(16,3)--(3,3)--(3,0) [fill=blue!20];	

\draw[gray] (6.25,0)--(6.25,3);	\draw[gray] (9.5,0)--(9.5,3);
\draw[gray] (12.75,0)--(12.75,3);	
\draw[gray] (3.6,4.7)--(1.0,6.0);
\draw[gray] (4.8,6.7)--(2.45,8.45);
\draw[gray] (6.6,8.7)--(4.7,10.9);

\draw[gray] (15.4,4.7)--(17.8,5.9);
\draw[gray] (14.2,6.7)--(16.4,8.3);
\draw[gray] (12.4,8.8)--(14.3,10.8);

\draw[black] (3.3,1)--(5.95,1);
\draw[black] (6.55,1)--(9.2,1);
\draw[black] (9.8,1)--(12.45,1);
\draw[black] (13.05,1)--(15.7,1);
\draw[black] (0.3,0.3)--(0.3,3); \draw[black] (0.6,0.3)--(0.6,2.7);
\draw[black] (0.9,0.3)--(0.9,2.7);
%\draw[black] (16.3,0.3)--(16.3,2.7); \draw[black] (16.6,0.3)--(16.6,2.7);
\draw[black] (16.9,0.3)--(16.9,2.7); \draw[black] (17.2,0.3)--(17.2,2.7); 
\draw[black] (17.5,0.3)--(17.5,2.7); \draw[black] (17.8,0.3)--(17.8,2.7);
\draw[black] (18.1,0.3)--(18.1,2.7); \draw[black] (18.4,0.3)--(18.4,2.7); 
\draw[black] (18.7,0.3)--(18.7,3); 

\draw[black] (0.3,3) .. controls (0.8,5.8) and (4.3,11.2) .. (8,12.7);
\draw[black] (0.3,3)--(0.6,3.3) .. controls (1.1,5.5) and (4.3,10.9) .. (8,12.4);

\draw[black] (8.3,10.3)--(10.7,10.3); \draw[black] (8.3,10.6)--(10.7,10.6);
\draw[black] (8.3,10.9)--(10.7,10.9); \draw[black] (8.3,11.2)--(10.7,11.2);
\draw[black] (8.3,11.5)--(10.7,11.5); \draw[black] (8.3,11.8)--(10.7,11.8);
\draw[black] (8.3,12.1)--(10.7,12.1); \draw[black] (8.3,12.4)--(10.7,12.4);
\draw[black] (8.3,10.3)--(10.7,10.3); \draw[black] (8.3,10.6)--(10.7,10.6);
\draw[black] (8.3,12.7)--(10.7,12.7); 
\draw[black] (8.3,10.3)--(8.3,12.7); \draw[black] (8.6,10.3)--(8.6,12.7);
\draw[black] (8.9,10.3)--(8.9,12.7); \draw[black] (9.2,10.3)--(9.2,12.7);
\draw[black] (9.5,10.3)--(9.5,12.7); \draw[black] (9.8,10.3)--(9.8,12.7);
\draw[black] (10.1,10.3)--(10.1,12.7); \draw[black] (10.4,10.3)--(10.4,12.7);
\draw[black] (10.7,10.3)--(10.7,12.7); 

\draw[black] (0.9,3.3) .. controls (1.2,4.35) and (1.2,4.35) .. (1.75,5.4);
\draw[black] (1.9,5.7) .. controls (2.45,6.7) and (2.45,6.7) .. (3.15,7.75);
\draw[black] (3.3,7.95) .. controls (4.1,9) and (4.1,9) .. (5.2,10.1);
\draw[black] (5.4,10.3) .. controls (7.0,11.7) and (7.0,11.55) .. (8,12.1);
\draw[red] (-0.2,-0.2)--(3.2,-0.2)--(3.2,3.2) .. controls (4.1,5.7) and (4.8,7.2) .. 
(8.2,9.8)--(8.2,13.2) .. controls (4.2,11.8) and (0.3,6.7) .. (-0.2,3.2)--(-0.2,-0.2);

\draw[black] (18.7,3) .. controls (18.2,5.8) and (14.7,11.2) .. (11,12.7);
\draw[black] (18.7,3)--(18.4,3.3) .. controls (17.9,5.5) and (14.7,10.9) .. (11,12.4);
\draw[black] (18.1,3.3) .. controls (17.8,4.35) and (17.8,4.35) .. (17.25,5.4);
\draw[black] (17.1,5.7) .. controls (16.55,6.7) and (16.55,6.7) .. (15.85,7.75);
\draw[black] (15.7,7.95) .. controls (14.9,9) and (14.9,9) .. (13.8,10.1);
\draw[black] (13.6,10.3) .. controls (12,11.7) and (12,11.55) .. (11,12.1);

\draw[black] (17.8,3.3) .. controls (17.5,4.25) and (17.5,4.25) .. (17.04,5.25);
\draw[black] (17.5,3.3) .. controls (17.3,4.05) and (17.3,4.05) .. (16.8,5.1);
\draw[black] (17.2,3.3) .. controls (17,3.95) and (17,3.95) .. (16.55,5);
\draw[black] (16.9,3.3) .. controls (16.8,3.75) and (16.8,3.75) .. (16.3,4.9);

\draw[black] (16.85,5.6) .. controls (16.3,6.5) and (16.3,6.5) .. (15.65,7.55);
\draw[black] (16.55,5.45) .. controls (16.1,6.3) and (16.1,6.3) .. (15.45,7.35);
\draw[black] (16.3,5.35) .. controls (15.9,6.1) and (15.9,6.1) .. (15.2,7.2);
\draw[black] (16.05,5.23) .. controls (15.7,5.9) and (15.7,5.9) .. (15,7);

\draw[black] (15.45,7.75) .. controls (14.55,8.9) and (14.55,8.9) .. (13.7,9.85);
\draw[black] (15.2,7.6) .. controls (14.35,8.7) and (14.35,8.7) .. (13.5,9.6);
\draw[black] (15,7.4) .. controls (14.15,8.5) and (14.15,8.5) .. (13.3,9.4);
\draw[black] (14.8,7.2) .. controls (13.95,8.3) and (13.95,8.3) .. (13.1,9.2);

\draw[black] (13.4,10.1) .. controls (12,11.4) and (12,11.2) .. (11,11.8);
\draw[black] (13.2,9.9) .. controls (12,11) and (12,10.9) .. (11,11.5);
\draw[black] (13,9.7) .. controls (12,10.7) and (12,10.6) .. (11,11.2);
\draw[black] (12.8,9.5) .. controls (12,10.35) and (12,10.25) .. (11,10.9);

\draw (-0.8,-0.8) node{\small{\color{blue} $[x_0]$}};
\draw (19.8,-0.8) node{\small{\color{blue} $[x_1]$}};
\draw (9.5,13.8) node{\small{\color{blue} $[x_2]$}};
\draw (1.5,-4.3) node{\small{\color{blue} $[L_1]$}};
\draw (17.5,-4.3) node{\small{\color{blue} $[L_0]$}};
\draw (10.5,-4.3) node{\small{\color{blue} $q-2$ segments}};
\draw (23.3,1.5) node{\small{\color{blue} $[L_2]$}};
\draw (-0.5,8) node{\small{\color{red} $\AG(2,q)$}};
\draw[gray,dotted,->] (1.5,-0.5)--(1.5,-3.5);
\draw[gray,dotted,->] (17.5,-0.5)--(17.5,-3.5);
\draw[gray,dotted,->] (19.5,1.5)--(22.5,1.5);
\draw[gray, dotted,->] (13.5,-3.5)--(16.9,-0.5);

\end{tikzpicture}		
\end{center}
\caption{The construction from Theorem~\ref{thml:s=q+1}}
\end{figure}
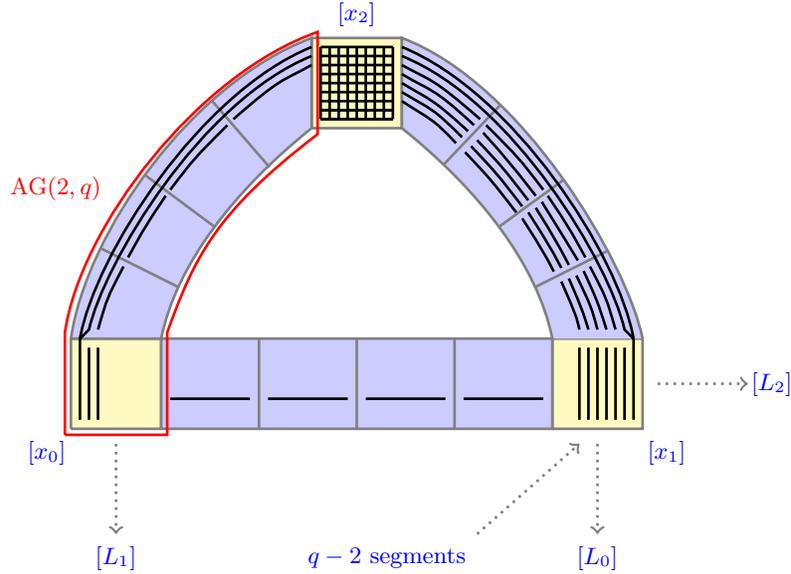

The proof in the case $s=q+2$ is virtually the same, except that
 $\mset{K}$ is defined as
 $\mset{K}=(q-2)\delta_{L_0}+2\delta_{L_1}+2\delta_{L_2}$. As for
 the necessary exceptions, the case $q=2$ is completely trivial and $q=3$ is
 covered by the (improved) construction in Th.~\ref{thml:q=3:n=7}
\qed

\begin{remark}
  \label{rmk:s=q+1}
  We know of no single instance where the bound of
  Th.~\ref{thml:s=q+1} 
  %for the case $n=q^2-1$ 
  can be improved. So
  one might conjecture that $\maxsize_{q^2-1}(R)=q^4-q^2-q$ is true
  for all $R$. 
\end{remark}

It is not difficult to construct $\bigl(q(q+1)n-q^3,n\bigr)$-arcs
also for $q^2-q\leq n\leq q^2-2$ by using the same
ideas. The following
theorem, however, improves on this for $n$ in the range $q^2-q\leq
n\leq q^2-3$.

\begin{theorem}
  \label{thml:s=q+2..2q}
  For every finite chain ring $R$ and every integer $n$ satisfying
  $q^2-q\leq n\leq q^2-2$
  there exist
  $(q^2n-2q,n)$-arcs
  in $\PHG(2,R)$. Consequently, 
  \begin{equation*}
  \maxsize_n(R)\geq q(q+1)n-q^3+q(q-2-t),
\end{equation*}
where $t\in\{0,1,\dots,q-2\}$ is defined by $n=q^2-q+t$.
\end{theorem}

\noindent {\it Proof.}
  Setting $n=q^2-q+t$, we have $q^2n-2q=q^4-q^3+tq^2-2q$ and $q(q+1)n-q^3+q(q-2-t)
  =q^4-q^2-q^3+t(q^2+q)+q^2-2q-tq=q^4-q^3+tq^2-2q$, 
  so both assertions in the theorem are indeed equivalent.
  We will construct the required arcs as unions of $qn-2$ line segments.
  The case $q=2$ is trivial and will be excluded from now
  on.
  
  We fix a point class $[z]$ in $\PHG(2,R)$ and denote the line classes
  through $[z]$ by $[L_0]$, $[L_1]$, \ldots, $[L_q]$.
  Writing $n=q^2-q+t$ with $0\leq t\leq q-2$, we choose $t+2$ further
  point classes $[x_1]$, \ldots, $[x_{t+2}]$ on the line class
  $[L_0]$.

  Now we define a point set $\mset{K}$ in $\PHG(2,R)$ as follows. In
  each point class $[x]$ contained in $[L_1]\cup[L_2]\cup\dots\cup[L_{t+2}]$
  and different from $[z]$ we choose $q-1$ parallel line segments. The
  direction of these line segments is determined as $[xx_i]$, where
  $[L_i]$ is the line class containing $[x]$. 

  In each class $[x]$
  contained in $[L_{t+3}]\cup[L_{t+4}]\cup\dots\cup[L_q]$ and
  different from $[z]$ we choose
  $q-2$ parallel line segments in the direction $[L_i]$, where again
  $i$ is determined by $[x]\in[L_i]$. The line segments incident with
  a fixed line class $[L_i]$ are arranged in such a way that they form
  a $(q^2-2q,q-1)$-arc in the affine plane defined on $[L_i]\setminus[z]$.

  Finally, in the classes $[x]$ on $[L_0]$ we take line segments in
  the direction of $[L_0]$ determined as follows: $q-2$
  line segments if $[x]$ is distinct from
  $[z]$ and also distinct from $[x_i]$ for $1\leq i\leq t+2$; $q-3$
  line segments if $[x]=[x_i]$ for some $1\leq i\leq t+2$; $t$ line
  segments if $[x]=[z]$. We can arrange these line segments in such a
  way that they form, together with the missing point $p_\infty$,
  a $(q^2-2q-1,q-1)$-arc in the
  projective plane defined on $[L_0]$. It is easier to obtain this arc by deleting a 
    $(3q+2,2)$-blocking set from the projective plane of order  $q$
    defined on $[L_0]$. The latter can be taken as a triangle having $p_\infty$
    as a $1$-point, from which $t+1$
    points $p_1=p_\infty,$ \ldots, $p_{t+1}$ on one side of the
    triangle are deleted and 
    $t+2$ further points $q_1,\dots,q_{t+2}$ 
    are added in such a way that they block all
    lines $p_ir$, $1\leq i\leq t+1$, connecting a point $p_i$ to the
    opposite vertex $r$ of the triangle, and such that the lines
    $p_\infty q_i$, $1\leq i\leq t+2$, are all different.

  The total number of line segments in $\mset{K}$ is 
  \begin{multline*}
    (t+2)(q^2-q)+(q-2-t)(q^2-2q)+(q^2-2q-2)\\
    =q^3-q^2-2q-2+(t+2)q=q(q^2-q+t)-2=qn-2
  \end{multline*}
  as required.

  Lines in $[L_i]$, $1\leq i\leq t+2$, have mutliplicity exactly $n$,
  since the restriction of $\mset{K}$ to $[L_i]$ consists of $n$ line
  segments having a direction different from $[L_i]$. Lines in
  $[L_i]$, $t+2\leq i\leq q$, have mutliplicity at most $q(q-1)+t=n$,
  since they meet $\mset{K}\cap[z]$ in $t$ points and
  $\mset{K}\setminus[z]$ in no more than $q-1$ line segments.
  Similarly, lines in $[L_0]$ have mutliplicity at most $q(q-1)\leq
  n$.

  Finally let $L$ by a line with $[z]\notin[L]$. If also
  $[x_i]\notin[L]$ for $1\leq i\leq t+2$, then $L$ contains no line
  segment of $\mset{K}$ and so 
  has multiplicity $(t+2)(q-1)+(q+1-t-2)(q-2)=q^2-q+t=n$.
  Otherwise $[x_i]\in[L]$ for a unique $i\in\{1,\dots,t+2\}$ and $L$
  contains at most one line segment of $\mset{K}$, whence $L$
  has multiplicity $q+(t+1)(q-1)+(q-t-2)(q-2)+(q-3)=n$ or $n-q$.  
\qed
% \begin{remark}
%   \label{rmk:s=q+2..2q}
%   For $n=q^2-q+t$ with $0\leq t\leq q-2$ we have
%   \begin{equation*}
%     q^2n-2q=q(q+1)n-q^3+q(q-t-2),
%   \end{equation*}
%   showing that Th.~\ref{thml:s=q+2..2q} improves for $q^2-q\leq n\leq
%   q^2-3$ on the construction of $\bigl(q(q+1)n-q^3,n\bigr)$-arcs
%   outlined in Remark~\ref{rmk:s=q+1..q+2}.
% \end{remark}

Next we exhibit a construction of $(k,n)$-arcs for a selected range of
small values of $n$. For $n<2q$ no satisfactory general constructions for
are known except for the case $n=2$, where we have the hyperovals
($(q^2+q+1,2)$-arcs) in the planes over Galois rings of characteristic
$4$ (see \cite{it:ovals}) and the ovals ($(q^2,2)$-arcs) in the planes
over the chain rings $\F_q[X;\sigma]/(X^2)$ of prime characteristic (see
\cite{mt:ovals}).

We now consider the range $2q\leq n\leq 3q-2$.

\begin{theorem}
  \label{thml:n=2q..3q-2}
  Let $R$ be a finite chain ring of length $2$.
  There exist
    $\bigl(q^3+(t+4)q^2+(t-4)q-2t,2q+t\bigr)$-arcs for every integer
    $t\in\{0,1,\dots,q-5\}$ in $\PHG(2,R)$, 
    as well as $(2q^3+q^2-12q+16,3q-4)$-arcs, $(2q^3+q^2-5q,3q-3)$-arcs, and
  $(2q^3+q^2-2q,3q-2)$-arcs.
  Consequently, for $n=2q+t$, $0\leq t\leq q-2$, we have
  \begin{equation*}
    \maxsize_n(R)\geq
    \begin{cases}
      q(q+1)n-q^3+2q^2-4q-2t&\text{if $t\leq q-5$},\\
      q(q+1)n-q^3+2q^2-8q+16&\text{if $t=q-4$},\\
      q(q+1)n-q^3+q^2-2q&\text{if $t=q-3$},\\
      q(q+1)n-q^3&\text{if $t=q-2$}.
    \end{cases}
  \end{equation*}
\end{theorem}

\noindent {\it Proof.}
  The construction is similar to that in
  Th.~\ref{thml:s=q+2..2q}. Using the same notation as in its proof and
  writing $n=2q+t$ with $0\leq t\leq q-2$, we define a multiset
  $\mset{K}$ as follows. In each point class $[x]\neq[z]$ contained in
  $[L_1]\cup\dots\cup[L_{t+2}]$ we choose two parallel line segments with
  direction of $[xx_i]$. In each class $[x]\neq[z]$ contained in
  $[L_{t+2}]\cup\dots[L_q]$ we choose one line segment with
  direction $[xz]$, say $[xz]=[L_i]$, and arrange the line segments in
  $[L_i]$ to form a $(q,2)$-arc in the affine plane defined on
  $[L_i]\setminus[z]$.
  
  In $[L_0]$ we choose $(q-2)(2q+t)$ points contained in three fixed
  classes $[y_1]$, $[y_2]$, $[y_3]$ distinct from $[z]$ and all $[x_i]$.
  The points are chosen in such a way that no line in
  $[L_0]$ contains more than $2q+t$ of them and no line outside $[L_0]$
  contains more than $q-2$ of them (i.e.\ a $(q-2)\times(2q+t)$ ``rectangle'').
  The construction is possible as long as $t\leq q-5$. 

  For $t=q-4$
  there are only two such classes $[y_1]$, $[y_2]$ and we modify the
  construction as follows: We take $q-2$ parallel line segments with
  direction $[L_0]$ in each of $[y_1]$, $[y_2]$. Further, we take
  $(q-4)^2$ points meeting each line in at most $q-4$ points
  in the affine plane defined on $[z]$ (i.e.\
  a $(q-4)\times(q-4)$ ``square'').
  
  For $t=q-3$ there is only one such class $[y_1]$, and we take $q-2$
  line segments in $[y_1]$ and $q-3$ line segments in $[z]$, all with
  direction $[L_0]$. 
  
  In the remaining case $t=q-2$ we take
  $q-2$ line segments in $[z]$ with direction $[L_0]$.
\qed

\begin{remark}
  \label{rmk:n=2q..3q-2}
  In the special case $n=2q$,
  Theorem~\ref{thml:n=2q..3q-2} gives $\maxsize_{2q}(R)\geq
  q^3+4q^2-4q$ and improves on the bound 
  $\maxsize_{2q}(R)\geq q^3+2q^2$ from
  \cite[Ex.~4.1]{it:kent}.

  For large $t$ (more precisely, for $t\geq q-3$ if $q$ is
  sufficiently large) there is a better
  construction placing $2q+t$ line segments resp.\ $q$ line
  segments in about half of the classes $[L_i]$. This construction
  covers also the case $n=3q-1$ which is not covered by the above theorem. Both
  Theorem~\ref{thml:n=2q..3q-2} and the alternative construction
  just mentioned can also be generalized to the case $n\geq
  3q$. 
\end{remark}

\subsection{Special Constructions for $q=3,4,5$}
  \label{ssec:speconst}
  
%\subsection{Several special constructions}

In this subsection we present several sporadic constructions for arcs,
which do not generalize to larger values of $q$.

We start with a construction of $(60,7)$-arcs in $\PHG(2,R)$,
$\abs{R}=9$, which improves upon the $(57,7)$-arcs of
Theorem~\ref{thml:s=q+2..2q}. The construction does not depend 
on the particular ring. We state the result in its
blocking set form.

\begin{theorem}
\label{thml:q=3:n=7}
There exist $(57,5)$-blocking sets in $\PHG(2,\Z_9)$ and
$\PHG(2,\DN_3)$.
\end{theorem}

\noindent {\it Proof.}
  Let $R$ be one of the rings $\Z_9$, $\DN_3$ and
  consider the factor geometry $\overline{\Pi}=\pi/\eye\cong\PG(2,3)$ of $\Pi=\PHG(2,R)$.
  Denote by $[x_1]$, $[x_2]$, $[x_3]$, $[x_4]$ the points of a fixed oval
  (quadrangle) in $\overline{\Pi}$, by $[L_i]$ the tangent to the oval in
  $[x_i]$, by $[S_{ij}]=[x_ix_j]$ the secants of the oval, by
  $[y_{12}]$, $[y_{13}]$, $[y_{14}]$, $[y_{23}]$, $[y_{24}]$, $[y_{34}]$
  the external points of the oval, indexed in such a way that
  $[y_{ij}]=[L_i]\cap[L_j]$, and by $[E_1]$, $[E_2]$, $[E_3]$ resp.\
  $[z_{12}]$, $[z_{13}]$, $[z_{23}]$ the external lines resp.\
  internal points of the oval, indexed in such a way that 
  $[z_{ij}]=[E_i]\cap[E_j]$.
  
  Define a set $\mset{K}$ of $57$ points contained in $19$ line
  segments of $\PHG(2,R)$ as follows.

  In each point class $[x_i]$ on the oval place $1$ line segment having the direction 
  of the tangent at this point.
  In each class $[y_{ij}]$ place $2$ parallel line segments in one of
  the two tangent directions according to the following rule: 
  In $[y_{i4}]$, $i=1,2,3$, choose
  $[L_i]$; in $[y_{12}]$ choose $[L_1]$; in $[y_{13}]$ choose $[L_3]$;
  in $[y_{23}]$ choose $[L_2]$.
  (This ``cyclic'' arrangement of the line segments in $[y_{12}]$,
  $[y_{13}]$, $[y_{23}]$ ensures
  that each tangent to the oval has the same
  direction as the line segments in exactly $0$ or $2$ of the three classes
  $[y_{ij}]$ incident with it.) Further, for $i=1,2,3$ 
  arrange the $5$ line segments having direction $[L_i]$ to form a
  blocking set in the affine plane defined on $[L_i]\setminus[y]$,
  where $[y]$ is the point class on $[L_i]$ whose line segments do not
  have direction $[L_i]$.
  Finally, in each class $[z_i]$ place $1$ line segment in one of the
  two external
  directions, again arranged in a ``cyclic fashion'' (choose
  directions $[E_1]$ for
  $[z_{12}]$, $[E_3]$ for $[z_{13}]$ and $[E_2]$ for $[z_{23}]$).

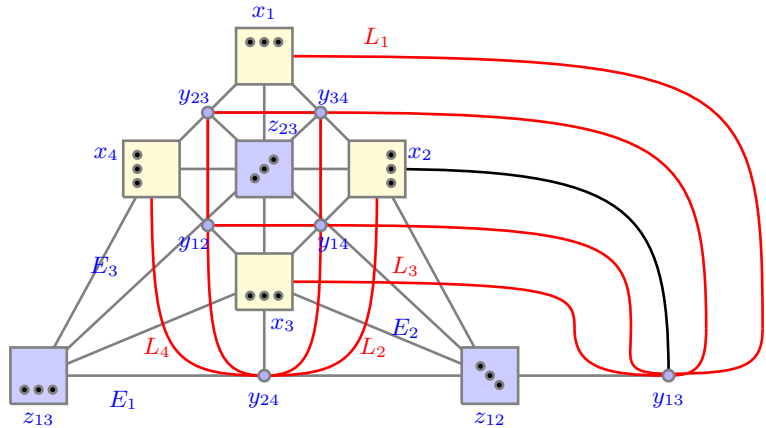
\begin{figure}[h!]
	\begin{center}
		\begin{tikzpicture}[line width=1pt, scale=0.25]

\draw[gray] (1.5,1.5)--(7.5,12.5)--(13.5,18.5);
\draw[gray] (1.5,1.5)--(13.5,6.5)--(19.5,12.5);
\draw[gray] (25.5,1.5)--(19.5,12.5)--(13.5,18.5);
\draw[gray] (25.5,1.5)--(13.5,6.5)--(7.5,12.5);
\draw[gray] (1.5,1.5)--(35,1.5);
\draw[gray] (13.5,1.5)--(13.5,18.5);
\draw[gray] (7.5,12.5)--(19.5,12.5);
\draw[gray] (1.5,1.5)--(16.5,15.5);
\draw[gray] (25.5,1.5)--(10.5,15.5);
\draw[gray] (10.5,9.5)--(16.5,9.5)--(16.5,15.5)--(10.5,15.5)--(10.5,9.5);
\draw[black] (19.5,12.5) .. controls (33,12.5) and (35,10.5) .. (35,1.5);
\draw[red] (10.5,15.5)--(10.5,9.5) .. controls (10.5,3) and (11,1.5) .. (13.5,1.5) ..
controls (19,1.5) and (19.5,3) .. (19.5,12.5);
\draw[red] (16.5,15.5)--(16.5,9.5) .. controls (16.5,3) and (16,1.5) .. (13.5,1.5) ..
controls (8,1.5) and (7.5,3) .. (7.5,12.5);
\draw[red] (10.5,15.5)--(16.5,15.5) .. controls (35,15.5) and (37,14) .. (37,4);
\draw[red] (37,4) .. controls (37,2) and (36.5,1.5) .. (35,1.5) .. controls 
(31.5,1.5) and (30,2) ..(30,4) .. controls (30,6) and (29.5,6.5) .. (13.5,6.5);
\draw[red] (10.5,9.5)--(16.5,9.5) .. controls (32,9.5) and (33,8.5) .. (33,4) ..
controls (33,2.3) and (32.5,1.6) .. (35,1.6) .. controls (39,1.6) and (40,2) .. (40,4);
\draw[red] (40,4) .. controls (40,17.5) and (39,18.5) .. (13.5,18.5);

\draw[gray] (0,0)--(3,0)--(3,3)--(0,3)--(0,0)  [fill=blue!20];
\draw[gray] (24,0)--(27,0)--(27,3)--(24,3)--(24,0) [fill=blue!20];
\draw[gray] (12,5)--(15,5)--(15,8)--(12,8)--(12,5) [fill=yellow!20];
\draw[gray] (6,11)--(9,11)--(9,14)--(6,14)--(6,11) [fill=yellow!20];
\draw[gray] (12,11)--(15,11)--(15,14)--(12,14)--(12,11) [fill=blue!20];
\draw[gray] (18,11)--(21,11)--(21,14)--(18,14)--(18,11) [fill=yellow!20];
\draw[gray] (12,17)--(15,17)--(15,20)--(12,20)--(12,17) [fill=yellow!20];

\draw[gray]  (13.5,1.5) circle (0.3cm) [fill=blue!30];
\draw[gray]  (35,1.5) circle (0.3cm) [fill=blue!30];
\draw[gray]  (10.5,15.5) circle (0.3cm) [fill=blue!30];
\draw[gray]  (16.5,15.5) circle (0.3cm) [fill=blue!30];
\draw[gray]  (10.5,9.5) circle (0.3cm) [fill=blue!30];
\draw[gray]  (16.5,9.5) circle (0.3cm) [fill=blue!30];

\draw(13.5,0.2) node{\small {\color{blue} $y_{24}$}};
\draw(35,0.2) node{\small {\color{blue} $y_{13}$}};
\draw(1.5,-0.8) node{\small {\color{blue} $z_{13}$}};
\draw(25.5,-0.8) node{\small {\color{blue} $z_{12}$}};
\draw(14.5,4.2) node{\small {\color{blue} $x_{3}$}};
\draw(13.5,20.8) node{\small {\color{blue} $x_{1}$}};
\draw(5.1,13.3) node{\small {\color{blue} $x_{4}$}};
\draw(21.8,13.3) node{\small {\color{blue} $x_{2}$}};
\draw(14.5,14.7) node{\small {\color{blue} $z_{23}$}};

\draw(9.8,16.3) node{\small {\color{blue} $y_{23}$}};
\draw(17.2,16.3) node{\small {\color{blue} $y_{34}$}};
\draw(9.8,8.5) node{\small {\color{blue} $y_{12}$}};
\draw(17.2,8.5) node{\small {\color{blue} $y_{14}$}};

\draw(19.5,19.5) node{\small {\color{red} $L_{1}$}};
\draw(21,7.3) node{\small {\color{red} $L_{3}$}};
\draw(7.8,3) node{\small {\color{red} $L_{4}$}};
\draw(19.3,3) node{\small {\color{red} $L_{2}$}};
\draw(6,0.2) node{\small {\color{blue} $E_{1}$}};
\draw(5,7.3) node{\small {\color{blue} $E_{3}$}};
\draw(21,4) node{\small {\color{blue} $E_{2}$}};

\draw[gray]  (0.75,0.75) circle (0.2cm) [fill=black];
\draw[gray]  (1.5,0.75) circle (0.2cm) [fill=black];
\draw[gray]  (2.25,0.75) circle (0.2cm) [fill=black];
\draw[gray]  (26,1) circle (0.2cm) [fill=black];
\draw[gray]  (25.5,1.5) circle (0.2cm) [fill=black];
\draw[gray]  (25,2) circle (0.2cm) [fill=black];
\draw[gray]  (14,13) circle (0.2cm) [fill=black];
\draw[gray]  (13.5,12.5) circle (0.2cm) [fill=black];
\draw[gray]  (13,12) circle (0.2cm) [fill=black];

\draw[gray]  (12.75,5.75) circle (0.2cm) [fill=black];
\draw[gray]  (13.5,5.75) circle (0.2cm) [fill=black];
\draw[gray]  (14.25,5.75) circle (0.2cm) [fill=black];
\draw[gray]  (12.75,19.25) circle (0.2cm) [fill=black];
\draw[gray]  (13.5,19.25) circle (0.2cm) [fill=black];
\draw[gray]  (14.25,19.25) circle (0.2cm) [fill=black];

\draw[gray]  (6.75,11.75) circle (0.2cm) [fill=black];
\draw[gray]  (6.75,12.5) circle (0.2cm) [fill=black];
\draw[gray]  (6.75,13.25) circle (0.2cm) [fill=black];
\draw[gray]  (20.25,11.75) circle (0.2cm) [fill=black];
\draw[gray]  (20.25,12.5) circle (0.2cm) [fill=black];
\draw[gray]  (20.25,13.25) circle (0.2cm) [fill=black];

\end{tikzpicture}		
\end{center}
\caption{The construction from Theorem~\ref{thml:q=3:n=7}}
\end{figure}

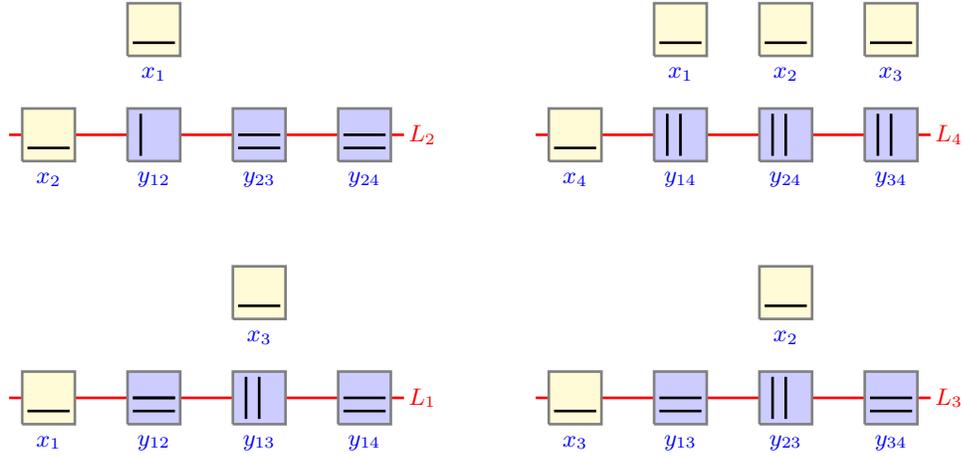
\begin{figure}[h!]
\begin{center}
\begin{tikzpicture}[line width=1pt, scale=0.35]
\draw[red] (-0.5,1)--(14.5,1);
\draw[red] (19.5,1)--(34.5,1);
\draw[red] (-0.5,11)--(14.5,11);
\draw[red] (19.5,11)--(34.5,11);
		
\draw[gray] (0,0)--(2,0)--(2,2)--(0,2)--(0,0) [fill=yellow!20];
\draw[gray] (4,0)--(6,0)--(6,2)--(4,2)--(4,0) [fill=blue!20];
\draw[gray] (8,0)--(10,0)--(10,2)--(8,2)--(8,0) [fill=blue!20];
\draw[gray] (12,0)--(14,0)--(14,2)--(12,2)--(12,0) [fill=blue!20];
\draw[gray] (8,4)--(10,4)--(10,6)--(8,6)--(8,4) [fill=yellow!20];

\draw[gray] (20,0)--(22,0)--(22,2)--(20,2)--(20,0) [fill=yellow!20];
\draw[gray] (24,0)--(26,0)--(26,2)--(24,2)--(24,0) [fill=blue!20];
\draw[gray] (28,0)--(30,0)--(30,2)--(28,2)--(28,0) [fill=blue!20];
\draw[gray] (32,0)--(34,0)--(34,2)--(32,2)--(32,0) [fill=blue!20];
\draw[gray] (28,4)--(30,4)--(30,6)--(28,6)--(28,4) [fill=yellow!20];

\draw[gray] (0,10)--(2,10)--(2,12)--(0,12)--(0,10) [fill=yellow!20];
\draw[gray] (4,10)--(6,10)--(6,12)--(4,12)--(4,10) [fill=blue!20];
\draw[gray] (8,10)--(10,10)--(10,12)--(8,12)--(8,10) [fill=blue!20];
\draw[gray] (12,10)--(14,10)--(14,12)--(12,12)--(12,10) [fill=blue!20];
\draw[gray] (4,14)--(6,14)--(6,16)--(4,16)--(4,14) [fill=yellow!20];

\draw[gray] (20,10)--(22,10)--(22,12)--(20,12)--(20,10) [fill=yellow!20];
\draw[gray] (24,10)--(26,10)--(26,12)--(24,12)--(24,10) [fill=blue!20];
\draw[gray] (28,10)--(30,10)--(30,12)--(28,12)--(28,10) [fill=blue!20];
\draw[gray] (32,10)--(34,10)--(34,12)--(32,12)--(32,10) [fill=blue!20];
\draw[gray] (24,14)--(26,14)--(26,16)--(24,16)--(24,14) [fill=yellow!20];
\draw[gray] (28,14)--(30,14)--(30,16)--(28,16)--(28,14) [fill=yellow!20];
\draw[gray] (32,14)--(34,14)--(34,16)--(32,16)--(32,14) [fill=yellow!20];

\draw(1,-0.7) node{\small {\color{blue} $x_{1}$}};
\draw(5,-0.7) node{\small {\color{blue} $y_{12}$}};
\draw(9,-0.7) node{\small {\color{blue} $y_{13}$}};
\draw(13,-0.7) node{\small {\color{blue} $y_{14}$}};
\draw(21,-0.7) node{\small {\color{blue} $x_{3}$}};
\draw(25,-0.7) node{\small {\color{blue} $y_{13}$}};
\draw(29,-0.7) node{\small {\color{blue} $y_{23}$}};
\draw(33,-0.7) node{\small {\color{blue} $y_{34}$}};

\draw(1,9.3) node{\small {\color{blue} $x_{2}$}};
\draw(5,9.3) node{\small {\color{blue} $y_{12}$}};
\draw(9,9.3) node{\small {\color{blue} $y_{23}$}};
\draw(13,9.3) node{\small {\color{blue} $y_{24}$}};
\draw(21,9.3) node{\small {\color{blue} $x_{4}$}};
\draw(25,9.3) node{\small {\color{blue} $y_{14}$}};
\draw(29,9.3) node{\small {\color{blue} $y_{24}$}};
\draw(33,9.3) node{\small {\color{blue} $y_{34}$}};

\draw(9,3.3) node{\small {\color{blue} $x_{3}$}};
\draw(29,3.3) node{\small {\color{blue} $x_{2}$}};
\draw(5,13.3) node{\small {\color{blue} $x_{1}$}};
\draw(25,13.3) node{\small {\color{blue} $x_{1}$}};
\draw(29,13.3) node{\small {\color{blue} $x_{2}$}};
\draw(33,13.3) node{\small {\color{blue} $x_{3}$}};

\draw(15.2,1) node{\small {\color{red} $L_{1}$}};
\draw(35.2,1) node{\small {\color{red} $L_{3}$}};
\draw(15.2,11) node{\small {\color{red} $L_{2}$}};
\draw(35.2,11) node{\small {\color{red} $L_{4}$}};

\draw[black] (0.2,0.5)--(1.8,0.5);
\draw[black] (4.2,0.5)--(5.8,0.5);
\draw[black] (4.2,1)--(5.8,1);
\draw[black] (12.2,0.5)--(13.8,0.5);
\draw[black] (12.2,1)--(13.8,1);
\draw[black] (20.2,0.5)--(21.8,0.5);
\draw[black] (24.2,0.5)--(25.8,0.5);
\draw[black] (24.2,1)--(25.8,1);
\draw[black] (32.2,0.5)--(33.8,0.5);
\draw[black] (32.2,1)--(33.8,1);
\draw[black] (8.5,0.2)--(8.5,1.8);
\draw[black] (9,0.2)--(9,1.8);
\draw[black] (28.5,0.2)--(28.5,1.8);
\draw[black] (29,0.2)--(29,1.8);

\draw[black] (0.2,10.5)--(1.8,10.5);
\draw[black] (4.5,10.2)--(4.5,11.8);
\draw[black] (4.2,1)--(5.8,1);
\draw[black] (8.2,10.5)--(9.8,10.5);
\draw[black] (8.2,11)--(9.8,11);
\draw[black] (12.2,10.5)--(13.8,10.5);
\draw[black] (12.2,11)--(13.8,11);
\draw[black] (20.2,10.5)--(21.8,10.5);
\draw[black] (24.5,10.2)--(24.5,11.8);
\draw[black] (25,10.2)--(25,11.8);
\draw[black] (28.5,10.2)--(28.5,11.8);
\draw[black] (29,10.2)--(29,11.8);
\draw[black] (32.5,10.2)--(32.5,11.8);
\draw[black] (33,10.2)--(33,11.8);

\draw[black] (8.2,4.5)--(9.8,4.5);
\draw[black] (28.2,4.5)--(29.8,4.5);
\draw[black] (4.2,14.5)--(5.8,14.5);
\draw[black] (24.2,14.5)--(25.8,14.5);
\draw[black] (28.2,14.5)--(29.8,14.5);
\draw[black] (32.2,14.5)--(33.8,14.5);

\end{tikzpicture}		
\end{center}
\caption{The directions of the line segments in the neighbour classes $[x_i]$
(Theorem~\ref{thml:q=3:n=7})}
\end{figure}
\vspace*{1cm}

  It is easily checked that the so-defined point set $\mset{K}$ is
  indeed a $(57,5)$-blocking set.\footnote{Lines in $[L_1]$ are blocked
  exactly $6$ times, lines in $[L_2]$, $[L_3]$, $[L_4]$ either $5$ or
  $8$ times, lines in $[S_{ij}]$ exactly $5$ times, and lines in
  $[E_i]$ either $5$, $8$ or $11$ times.}
\qed

\begin{remark}
\label{rmkl:q=3:n=8}
Omitting from $\mset{K}$ the three line segments contained in the
classes $[z_{ij}]$ we obtain $(48,4)$-blocking sets in $\PHG(2,\Z_9)$
resp.\ $\PHG(2,\DN_3)$, which are known to be optimal; cf.\
\cite{imt:two_exact_values}.
\end{remark}
 
Now we turn to the case $q=4$. First we improve the construction of
Th.~\ref{thml:n=2q..3q-2} for this particular case.

\begin{theorem}
  \label{thml:q=4:n=8..10}
  Suppose $R$ is one the chain rings $\G_4$, $\DN_4$, $\TDN_4$ of
  cardinality $16$. There exist arcs with parameters $(120,8)$,
  $(140,9)$ and $(152,10)$ in $\PHG(2,R)$.
\end{theorem}

\noindent {\it Proof.}
  Choose a non-incident point-line pair $[z]$, $[L]$ in the factor
  plane $\overline{\Pi}=\PG(2,4)$ and denote by $[x_i]$ and
  $[L_i]$, $0\leq i\leq 4$, the point classes on
  $[L]$ respectively the line classes through $[z]$.

  In each line class $[L_i]$ choose $6$ line segments with direction
  $[L_i]$ forming a
  hyperoval in the projective plane (of order $4$) defined on $[L_i]$
  and such that $[z]$ and $[L]$ remain empty.  (Since there exist
  hyperovals in $\PG(2,4)$ which are disjoint from any given pair of
  lines, this can be done.) Then every point class distinct from
  $[z]$ and all $[x_i]$ containes two line segments, 
  and it is immediately checked that the
  so-defined point set is a $(120,8)$-arc.\footnote{The
    only additional facts used is that every line class meets $[L]$ and
    therefore contains an empty point class.}  

  By placing one
  additional line segment with direction $[L]$ in each class $[x_i]$
  and arranging the line segments to form an oval of the projective
  plane defined on $[L]$ with nucleus in $p_\infty$, we can extend the
  $(120,8)$-arc to a $(140,9)$-arc.

  For the construction of an $(152,10)$-arc we choose $4$ point
  classes $[y_1]$, $[y_2]$, $[y_3]$, $[y_4]$ of the
  affine plane $\overline{\Pi}\setminus[L_0]$ forming an affine
  subplane whose diagonal points are the three points on $[L_0]$ which are
  distinct from $[z]$ and 
  $[x_0]$.\footnote{Equivalently, choose a projective subplane of $\Pi$
    containing $[L_0]$ and the three points on $[L_0]$ which are
    distinct from $[z]$ and $[x_0]$.}
  In each class $[y_i]$ we place $2$ parallel line segments with
  direction $[y_ix_0]$.
  The classes $[z]$ and $[x_0]$ remain empty. In each of the remaining $15$
  classes $[x]$ we place $2$ parallel line segments with direction
  $[xz]$ and such that the $6$ line segments with direction
  $[L_i]$, $0\leq i\leq 4$, form a hyperoval in the projective plane
  defined on $[L_i]$. (This step is exactly the same as for the
  previous two arcs.) As is easily checked, the so-defined set of
  $19\cdot 8=152$ points forms an $(152,10)$-arc.
\qed

\begin{theorem}
  \label{thml:q=4:n=11,12}
  Suppose $R$ is one the chain rings $\G_4$, $\DN_4$, $\TDN_4$ of
  cardinality $16$. There exist arcs with parameters $(166,11)$ and
  $(186,12)$ in $\PHG(2,R)$.
\end{theorem}

\noindent {\it Proof.}  
  Choose five point classes $[x_i]$, $1\leq i\leq 5$, forming an oval
  in $\overline{\Pi}$, and denote the nucleus of this oval by
    $[x_0]$. Let $[L_i]$, $0\leq i\leq 5$, be the $0$-classes
    (passants) of the hyperoval in $\overline{\Pi}$ 
    formed by $[x_i]$, $0\leq i\leq 5$.
    
  For both constructions we place $3$ parallel line segments in tangent
  direction $[x_ix_0]$ in each oval class $[x_i]$, $1\leq i\leq 5$,
  and a $(6,2)$-arc in the nucleus $[x_0]$.

  In the first construction we choose a $0$-class of lines, say $[L_0]$, of the
  hyperoval and place $1$ line segment with direction $[L_0]$ in each point
  class on $[L_0]$. The five line segments are arranged to form a
  $(5,2)$-arc in the projective plane defined on $[L_0]$. In the
  remaining ten point classes we place $2$ parallel line segments in
  one of the five possible
  passant directions in such a way that each such passant occurs $0$, $1$
  or $3$ times as the direction of a line segment pair. This is
  indeed possible: For example, choose direction $[L_i]$ in
  $[L_i]\cap[L_j]$ for all pairs $(i,j)$, $1\leq i<j\leq 5$, except
  $(1,3)$ and direction $[L_3]$ in $[L_1]\cap[L_3]$. Finally, if
  $[L_i]$ is a direction occuring $3$ times, arrange the corresponding
  $6$ line segments to form a hyperoval in the projective plane
  defined on $[L_i]$. It is then easily checked that the so-defined
  set of $5\cdot 12+6+5\cdot 4+10\cdot 8=166$ points forms a
  $(166,11)$-arc in $\PHG(2,R)$.
  
  In the second construction we place $2$ parallel line segments in
  a passant direction in each of the $15$ point classes outside the
  hyperoval. Again we can arrange the directions in a such a way
  that each of the six passants occurs $0$, $1$
  or $3$ times as the direction of a line segment pair: For example, in
  $[L_i]\cap[L_j]$, where $0\leq i<j\leq 5$, choose direction $[L_j]$
  if $i=0$, direction $[L_i]$ or $j-i\equiv 1,2\pmod{5}$, and
  direction $[L_j]$ if $j-i\equiv 3,4\pmod{5}$.\footnote{In this example
    each passant occurs either $0$ or $3$ times as the direction of a
    line segment pair.}
  The so-defined set of
  $5\cdot 12+6+15\cdot 8=186$ points forms a
  $(186,12)$-arc in $\PHG(2,R)$. 
\qed

\begin{theorem}
  \label{thml:q=4:n=13}
  Suppose $R$ is one of $\G_4$, $\DN_4$, $\TDN_4$.
  There exists a $(201,13)$-arc in $\PHG(2,R)$.
\end{theorem}

\noindent {\it Proof.}
  We choose three point classes $[x_1]$, $[x_2]$, $[x_3]$ of
  $\overline{\Pi}$ forming a triangle,
  and denote the sides of the triangle by $[L_1]$, $[L_2]$,
  $[L_3]$. 

  Now we define a point set $\mset{K}$ in $\PHG(2,R)$ as follows.
  The intersection of $\mset{K}$ with any
  point class on a triangle side $[L_i]$ but different
  from a triangle vertex consists of $2$ parallel line segments with
  direction $[L_i]$. The $6$ line segments with direction
  $[L_i]$, $i=1,2,3$, are arranged to form a hyperoval in the projective plane
  defined on $[L_i]$ and, in addition,
  such that $L_i$ determines a passant of this hyperoval.
  
  The lines $L_i$ determine sets of $7$ points in the classes $[x_j]$
  which look like a ``cross'' (unions of two line segments
  $L_i\cap[x_j]$) and form the intersections $\mset{K}\cap[x_j]$.
  
  Finally, the intersection of $\mset{K}$ with any of the $9$
  remaining point classes $[y_i]$ consists of $3$ parallel line segments,
  whose direction is one of the $9$ line classes $[M_j]$ not passing through a
  triangle vertex. Since each $[y_i]$ is on two line classes $[M_j]$
  and in turn each $[M_j]$ contains two point classes $[y_i]$ (i.e.\ the
  bipartite incidence graph formed by $[y_i]$ and $[M_j]$ is regular),
  we can arrange the directions in such a way that each class $[M_j]$
  occurs exactly once as a direction of a line segment triple (a
  perfect matching of this graph). It is then easy to verify that the
  so-defined set $\mset{K}$ of cardinality $9\times 8+3\times
  7+9\times 12=201$ is a $(201,13)$-arc.
\qed

Finally we turn to constructions for the case $q=5$. The first series
of constructions uses the same basic idea as in Theorems~\ref{thml:s=q+2..2q}
and~\ref{thml:n=2q..3q-2}, but also some special properties
of $(k,3)$-arcs in $\PG(2,5)$ \cite{bramwell-wilson73,soz96}.

\begin{theorem}
  \label{thml:q=5:n=15..17}
  There exist arcs with parameters $(355,15)$, $(375,16)$, and
  $(395,17)$ in both $\PHG(2,\Z_{25})$ and $\PHG(2,\DN_5)$.
\end{theorem}

\noindent {\it Proof.}
  We use the same notation as in the proof of Theorems~\ref{thml:s=q+2..2q}
  and~\ref{thml:n=2q..3q-2}. Accordingly we write $n=15+t$ with $0\leq
  t\leq 2$.

  In each point class $[x]$ contained in $[L_1]\cup\dots[L_{t+2}]$ we
  choose $3$ parallel line segments with direction $[xx_i]$. 
  In each point class $[x]$ contained in $[L_{t+3}]\cup\dots[L_5]$ we
  choose $2$ parallel line segments with direction $[xz]$ and arrange
  the $10$ line segments with direction $[L_i]$, $t+3\leq i\leq q$, to
  form a $(10,3)$-arc in the projective plane defined on $[L_i]$.  
  
  The remaining points (i.e.\ those in $[L_0]$) are chosen
  individually for each $n$. In all cases we specify (unions of)
  line segments, all with direction $[L_0]$.

  For $n=15$ $(t=0$) we choose $1$ line segment in each of $[x_1]$, $[x_2]$
  and $3$ line segments in each of $[x_3]$, $[x_4]$, $[x_5]$ in such a way
  that the $11$ segments chosen
  form an $(11,3)$-arc in the projective plane defined on
  $[L_0]$. This is possible by the structure of the $(11,3)$-arcs. 
  In total there are $71$ line segments, giving the required
  $(355,15)$-arc.

  For $n=16$ ($t=1$)  we choose $1$ line segment in each of $[z]$, $[x_1]$,
  $[x_2]$, $[x_3]$ and $3$ line segments in each of $[x_4]$, $[x_5]$,
  again in such a way that the $10$ segments chosen 
  form a $(10,3)$-arc in the projective plane defined on
  $[L_0]$.
  There are $75$ line segments in total, giving a $(375,16)$-arc.

  Finally, for $n=17$ ($t=2$) we choose a $(9,3)$-arc consisting of
  $2$ line segments in $[z]$, one line segment in each of $[x_1]$,
  $[x_2]$, $[x_3]$, $[x_4]$ and $3$ line segments in $[x_5]$.
  There are $79$ line segments in total, and
  we obtain a $(395,17)$-arc.

  % Finally, for $n=18$ we choose $3$ line segments in $[z]$ and $1$
  % line segment in each of $[x_1]$, \ldots, $[x_5]$.
  % Now there are $83$ line segments in total, giving a $(415,18)$-arc. 
\qed

\begin{theorem}
  \label{thml:q=5:n=18,19}
  There exist arcs with parameters $(425,18)$ and $(455,19)$
  in both $\PHG(2,\Z_{25})$ and $\PHG(2,\DN_5)$.
\end{theorem}

\noindent {\it Proof.}
  We use a ``triangle construction'' similar to that of
  Th.~\ref{thml:q=4:n=13} and the same notation as in its proof.  
  
  A triangle of $\overline{\Pi}\cong\PG(2,5)$ covers $15$ point classes
  and $15$ line classes, and the incidence relation of
  $\overline{\Pi}$ induces a regular bipartite graph of degree $3$ on
  the remaining $16$ point classes $[y]$ and $16$ line classes
  $[M]$. In each class $[y]$ we choose $4$ parallel line segments,
  assigning directions by way of a perfect matching $[y]\to[M]$ of
  this graph (so that each $[M]$ occurs exactly once as a direction of
  some class $[y]$). 

  For the second arc we choose in
  each point class on a triangle side $[L_i]$, $i=0,1,2$, distinct
  from a triangle vertex 
  $2$ parallel line segments with direction $[L_i]$. In each vertex
  $[x_i]$ we choose $1$ line segment with direction $[L_{i+1}]$, where
  indices are taken mod $3$ and $x_i$ is assumed to be the triangle vertex
  opposite to $[L_i]$. We arrange the
  $9$ line segments with direction $[L_i]$ to form an $(9,3)$-arc in the
  projective plane defined on $[L_i]$. The so-defined set of $16\times
  20+12\times 10+3\times 5=455$ points in $\PHG(2,R)$ forms a
  $(455,19)$-arc, as is easily verified.

  For the first arc we replace the double line segments on one
  triangle side, $[L_0]$ say, by single line segments and also empty
  the adjacent vertices $[x_1]$ and $[x_2]$. The resulting set of
  $16\times 20+8\times 10+5\times 5=425$ points forms the required $(425,18)$-arc. 
\qed

\subsection{Dual constructions}
\label{ssec:dualconst}

In general the dual incidence structure
$\Pi^*=(\mathcal{L},\mathcal{P},\supseteq)$ of
$\Pi=\PHG(2,R)=(\points,\lines,\subseteq)$ is
isomorphic to the left projective Hjelmslev plane over $R$. (This
follows from well-known duality properties of $R$; see
\cite{it:chain_survey} for example.) This implies $\Pi^*\cong\Pi$ for
all chain rings of length $2$ and cardinality at most $25$; cf.\ our
introductory remarks.

As in the classical case of
projective planes over finite fields, dualizing certain arcs gives rise
to good constructions. Here we consider only one important application of
this method. The topic is further explored in
\cite{it:dublin_duality}.

\begin{theorem}
\label{thm:dualarcs}
If $q=2^r$ is an even prime power, there exist
$\bigl((q^4-q)/2,q^2/2\bigr)$-arcs in $\PHG(2,\G_q)$. These arcs meet every
line of $\PHG(2,\G_q)$ in either $0$ or $q^2/2$ points and hence are
maximal; i.e., we have $\maxsize_{q^2/2}(\G_q)=(q^4-q)/2$.
\end{theorem} 

\noindent {\it Proof.}
  Define a point set $\mset{K}$ in the dual plane
  $\PHG(2,\G^\opp_q)\cong\PHG(2,\G_q)$ as follows: Take as points the
  passants ($0$-lines) of a $(q^2+q+1,2)$-arc (i.e., a hyperoval)
  $\mset{H}$ in $\PHG(2,\G_q)$. If $x$ is a point of $\PHG(2,\G_q)$
  outside $\mset{H}$, there exists a unique point $y\in\mset{H}$
  satisfying $x\eye y$. All $q$ lines through $x$ in the direction of
  $y$, say $[L]$, are $2$-lines of $\mset{H}$. The $q^2$ points of
  $\mset{H}$ outside $[L]$ split into $q^2/2$ pairs, each determining
  a unique $2$-line through $x$. The remaining $q^2/2$ lines through
  $x$ must then be passants, showing that each point of $\PHG(2,\G_q)$
  is incident with either $0$ or $q^2/2$ passants of $\mset{H}$. In
  other words, every line in the dual plane
  $\PHG(2,\G^\opp_q)\cong\PHG(2,\G_q)$ is incident with $0$ or $q^2/2$
  points of $\mset{K}$.

  The total number of passants of $\mset{H}$ is
  $q^4+q^3+q^2-\binom{q^2+q+1}{2}=(q^2+q+1)(q^2-q)/2=(q^4-q)/2$, with
  $q(q-1)/2$ of them in every line class $[L]$. Optimality then
  follows from the bound in Cor.~\ref{cor:bound}, which is sharp in
  this case.
\qed

Note that in the smallest case $q=2$ the dual of a hyperoval is itself
a hyperoval. For $q=4$ we get $(126,8)$-arcs in $\PHG(2,\G_4)$, which
beat the $(120,8)$-arcs of Th.~\ref{thml:q=4:n=8..10}. 

Now we show that
arcs with corresponding parameters do not exist in the planes over
other chain rings of cardinality $4^r$.

\begin{theorem}
\label{thm:dualimp}
If $R$ is one of the chain rings $\F_q[X;\sigma]/(X^2)$, $q=2^r$,
$\sigma\in\Aut(\F_q)$ (i.e., both length and characteristic of $R$ are
equal to $2$), there do not exist
$(q^4-q)/2,q^2/2)$-arcs. Equivalently, $\maxsize_{q^2/2}(R)<(q^4-q)/2$
for these rings.
\end{theorem}

\noindent {\it Proof.}
  Suppose on the contrary that there exists a $(q^4-q)/2,q^2/2)$-arc
  $\mset{K}$ in $\PHG(2,R)$.
%By Corollary 3.3 from \cite{it:kent}, 
By Cor.~\ref{cor:bound}, all point classes $[x]$
have the same multiplicity $q(q-1)/2$. Moreover, each parallel
class in the affine plane defined on $[x]$ contains one empty line and
$q-1$ lines of multiplicity $q/2$. In other words, the
restriction of $\mset{K}$ to each neighbour class of points has to be
a maximal $\bigl(q(q-1)/2,q/2\bigr)$-arc.  Now all non-empty lines
are $q^2/2$-lines and the empty lines form a $(q^2+q+1,2)$-arc
in the dual plane $\PHG(2,R^\opp)$, which is again a coordinate plane over a
chain ring of length $2$ and characteristic 2.  
But the latter does not exist; cf. \cite{it:ovals}. Contradiction.
\qed

\subsection{Computer constructions}
\label{ssec:bayconst}

As for all combinatorial objects, it is natural to use the computer
for the search for large arcs.  However, a complete backtracking
search as done in \cite{kiermaier06} is only feasable for very small
parameters $q$, $k$ and $n$.  Therefore, the complete search is
dropped in favor of other techniques.

\paragraph{Heuristic search.}
In \cite{Zwanzger-2008-IEEETIT54[5]:2388-2392.pdf} a heuristic approach
for the construction of good linear codes over finite fields was
discussed, which led to the discovery of several new codes.  In the
following, we adapt this approach to the search for arcs:

In each node of the search tree we sort the children nodes according
to some heuristic, such that the heuristic predicts the highest chance
for success when chosing the first child node.  Assuming a $(k,n)$-arc
$\mset{K}$ exists and that the heuristic is not too bad, $\mset{K}$
will be found quite early when only a small fraction of the full
search tree was explored.

At each search node we have a $(k^\prime,n)$-arc $\mset{K}'$ of size
$k^\prime < k$ which should be extended to a $(k,n)$-arc.  Now for
each point $x$ and each line $L$ in $\PHG(2,R)$ we compute the
relative frequency $p_{L}(x)$ that a randomly chosen multiset
$\mset{K}$ of points of size $k$ with
$\mset{K}'\cup\{x\}\subseteq\mset{K}$ does not violate the $n$-arc
condition on the line $L$.  The heuristic value for $x$ is then taken
as $\prod_L p_L(x)$, which resembles the probability of
$\mset{K}'\cup\{x\}$ being contained in a solution $\mset{K}$ if we
assume that the conditions for each individual line $L$ are
stochastically independent.  Of course this assumption is not true at
all: For example it may occur that for all $L$ it holds $p_L(x) > 0$
and therefore $p(x) > 0$, while $\mset{K}'\cup\{x\}$ is not contained
in any $(n,k)$-arc.  But our results indicate that for many instances
this heuristic gives a good hint which point to include next.  As an
example, the $(49,6)$-arc in $\PHG(2,\Z_9)$ as well as the
$(38,5)$-arc and the $(50,6)$-arc in $\PHG(2,\DN_3)$ were found by
this method.

By the fundamental theorem of projective Hjelmslev
geometry~\cite{alex91b} every collineation of a coordinate projective
Hjelmslev plane $\Pi=\PHG(2,R)$ over a finite chain ring $R$ is induced by
a semilinear automorphism of the underlying module $R^3_R$, and the
collineation group of $\PHG(2,R)$ is isomorphic to the projective
semilinear group 
\[\mathrm{P\Gamma L}(3,R)=\mathrm{\Gamma L}(3,R)/\mathrm{Z}(R),\] 
where $\mathrm{Z}(R)$ denotes the center of
the ring $R$.  The size of this group is
\[
q^{11}(q-1)^2(q+1)(q^2 + q + 1)\cdot\abs{\Aut(R)}.
\]

The pure backtrack search as well as the heuristic search are
complicated by the fact that the size of this collineation group is
growing quite rapidly.  Completely ignoring collineations is not an
option, since this would result in one and the same point
constellation showing up in thousands or even millions of isomorphic
copies during the search.  Therefore in the backtrack search in
\cite{kiermaier06} a complete isomorph rejection was done, for the
price of increasing the complexity of the algorithm.  In the heuristic
search a set of invariants is used to get rid of isomorphic copies, in
the same time risking to throw away point constellations which have
the same invariants, but actually are not isomorphic.

\paragraph{Using prescribed automorphisms.}
On the other hand, the large group of collineations of $\PHG(2,R)$ can
be exploited in the following way: We prescribe a group of
collineations $G<\mathrm{P\Gamma L}(3,R)$ to the search, such that we
impose the additional condition that each resulting arc has $G$
contained in its stablilizer.  This approach suggests itself for
several reasons:
\begin{itemize}
\item Experience shows that large or otherwise interesting arcs often
  have a big automorphism group.
\item A large group of automorphisms of an arc benefits a theroretic
  analysis of the results found by the computer.  In fact, many of the
  presented constructions in this section were found by a careful
  investigation of the computer results.
\item The prescription of $G$ means that the atomic building blocks in
  the search are no longer single points, but whole orbits of $G$
  acting on the point set.  The advantage over taking arbitrary
  partitions of the point set is that by using the action of $G$, also
  the line conditions are partitioned into orbits.  This leads a
  further reduction of the computational complexity, since for the
  formulation of the $n$-arc condition, only a single representative
  of each line orbit needs to be considered.
\end{itemize}

The technique of prescribing automorphism was successfully used for
linear codes and is described in detail in \cite{betten-etal06}.  A
part of the new results for arcs in projective Hjelmslev planes over
Galois rings was already published in \cite{kiermaier-kohnert07}.

It is convenient to state the existence problem of a $(k,n)$-arc with
the prescribed group of collineations $G$ as a system of linear
integer inequalities.  Then we try to find a solution for this system
using a solver based on the LLL-algorithm
\cite{wassermann-designs-enumeration-techniques-1998}.  Another
possibility which turned out to be quite successful is again the
application of an heuristic
\cite{Kohnert-Zwanzger-2009-AiMoC3[2]:157-166}.  We start with the
zero vector and try to modify it to a solution vector step-by-step by
repeatedly increasing some component of $x$ by $1$.  Similarly as
above, in each search node we compute the probability of each
increased variable being part in a randomly completed solution
assuming stochastic independence of the single equations.

\paragraph{Modification of existing arcs.}
An $n$-arc $\mset{K}$ found by the above methods not necessarily needs
to be complete, so we check if it can be enlarged to an even bigger
$n$-arc.  Furthermore, we try to find large $n^\prime$-arcs with
$n^\prime > n$ by adding further points to $\mset{k}$, and we try to
find large $n^\prime$-arcs with $n^\prime < n$ by removing points from
$\mset{k}$.  These extension and reduction problems are again
formulated as a system of linear integer inequalities and solved
either by the LLL- or by the heuristic solver algorithm.

The results of the search are available in the online tables \cite{arctables}.

\paragraph{Example}
A lifted Singer cycle of the factor plane $\PG(2,5)$ is a group of
collineations in $\mathrm{P\Gamma L}(3,R)$ of size $31$ partitioning
the point set into $25$ orbits of size $31$, respectively.  Our search
revealed that it is possible to join $10$ of these orbits to get an
$(310,13)$ arc $\mset{K}$.  Then we checked for the maximum number of
points that can be added to $\mset{K}$ to get a $14$-arc.  In this way
we found the $(319,14)$-arc over $\Z_{25}$.

The prescription of a lifted Singer cycle gives good arcs in quite a
few situations and can be treated more theoretically:

\paragraph{Lifted Singer cycles.}
Suppose now that $\Pi$ is a projective Hjelmslev plane over a Galois
ring $\G_q$, represented as $\PHG(\G_{q^3}/\G_q)$.
A generator $\eta$ of the Teichm\"uller
subgroup $\teich_{q^3}$ of $\G_{q^3}$ induces a collineation
$\sigma\in\Aut(\Pi)$ of order $q^2+q+1$, which acts as a Singer cycle
on the factor plane $\PG(2,q)$. There is obviously a one-to-one
correspondence between $\sigma$-invariant multisets in $\Pi$ and
multisets in a fixed point neighbour class of $\Pi$, for example
$[\G_q1]$. For a $\sigma$-invariant multiset $\mset{K}$ in $\Pi$ it is
possible to compute the $\mset{K}$-types of all lines in $\Pi$ from
certain combinatorial data of the corresponding multiset $\mset{k}$ in
$[\G_q1]\cong\AG(2,q)$. As shown in \cite{mt:phgtriangle}, suitable
choices of $\mset{F}$ yield $\sigma$-invariant arcs with good
parameters. As an example of this construction we mention a family of
arcs in the planes over $\G_p$, where $p$ is an odd prime, which
includes an optimal $(39,5)$-arc in the plane over $\Z_9$. A multiset
$\mset{F}$ in $\AG(2,p)$ is called a \emph{triangle set} if it is
affinely equivalent to the set $\bigl\{(x,y)\in\F_p^2;x+y<p-1\bigr\}$.
Here $\F_p=\{0,1,\dots,p-1\}$ is considered as a subset of $\Z$.
\begin{theorem}[\cite{mt:phgtriangle}]
  \label{thm:triangleset}
  For every odd prime $p$ there exists a $\sigma$-invariant
  $\bigl((p^4-p)/2,(p^2+p)/2-1\bigr)$-arc in the projective Hjelmslev
  plane over the Galois ring $\G_p$. The arc is induced from an
  appropriately chosen triangle set in $[\G_p1]\cong\AG(2,p)$.
\end{theorem}

\section{Upper Bounds for the Size of an Arc}
\label{sec:upper_bounds}

In this section we strenghten the known upper bounds for
$\maxsize_n(R)$ (cf.\
\cite{it:deadfin,it:kent}). We start with
the case of arbitrary $q$ and large $n$ ($q^2-q\leq n\leq q^2-1$) in
Section~\ref{ssec:largebound}. This case treated is separately for two
reasons: It serves as a particular example for the improvements in
Section~\ref{ssec:newbound}, and 
it is the only case where we also have a satisfactory lower
bound and hence are able two obtain a rather precise result
(Cor.~\ref{coru:q^2-s}). Then in Section~\ref{ssec:newbound} we
present a refined analysis of the general bound
(Fact~\ref{thm:bound}) and use this to improve the known upper bounds for
small chain rings ($\abs{R}\leq 25$). Finally, in
Section~\ref{ssec:spebounds} we reduce these bounds still further
for particular values of small $q$ and $n$. As a by-product, we are
able to determine the exact value of $\maxsize_n(R)$ in a few
cases.

\subsection{A New Upper Bound for Large $n$}
\label{ssec:largebound}

As mentioned in Section~\ref{ssec:genconst} the exact value of
$\maxsize_n(R)$ is known for $q^2\leq n\leq q^2+q$.  Here it will be
shown that for values of $n$ that are smaller, but close enough to
$q^2$, the value of $\maxsize_n(R)$ is contained within a relatively
small interval. The following theorem determines the upper end of this
interval.

\begin{theorem}
\label{thmu:q^2-s}
Let $\mset{K}$ be a $(k,n)$-arc in
$\PHG(R_R^3)$ with $q^2-q+1\leq n\leq q^2-1$. Then
\begin{equation}
  \label{eq:q^2-s}
  k\le q(q+1)n-q^3 + q^2-q.  
\end{equation}
\end{theorem}

\noindent {\it Proof.}
Let
$u=\max\bigl\{\mset{K}([y]);y\in\points\bigr\}$. If $u\geq q^2-q+1$,
then by the general bound (cf.\ Cor.~\ref{cor:bound})
% \begin{eqnarray*}
% k &\le & \mset{K}([x]) +q(q+1)(q^2-q-s) \\
%   &\le & q^4-q^2s-qs.
% \end{eqnarray*}
\begin{equation*}
   k\leq q(q+1)n-q(q+1)\lceil u/q\rceil+u\leq
   q(q+1)n-q^3.
\end{equation*}
If $q^2-2q+2\leq u\leq q^2-q$, then by the fact
that a blocking set in $\AG(2,q)$ has size at least $2q-1$, 
we have in every neighbour class $[x]$
with $\mset{K}([x])=u$ at least one direction with a line of multiplicity $q$.
Hence we get
\begin{eqnarray*}
k &\le & \mset{K}([x]) +q(n-q)+ q^2(n-q+1) \\
  &\le & q^2-q+q(q+1)n-q^3.
\end{eqnarray*}
In the remaining case $u\leq q^2-2q+1$ we have
\begin{equation*}
  k\leq (q^2-2q+1)(q^2+q+1)=q^4-q^3-q+1,
\end{equation*}
which for $n\geq q^2-q+1$ is less than the bound \eqref{eq:q^2-s}.
\qed

The steps from the above proof can be carried out also for the case
$n=q^2-q$. This implies the following corollary.

\begin{corollary}
\label{coru:q^2-s}
Let $\mset{K}$ be a $(k,q^2-q)$-arc in $\PHG(R_R^3)$.
Then $k\le q^4-q^3-q+1$. In case of equality, all neighbour classes
contain $(q-1)^2$ points each.
\end{corollary}

\begin{remark}
  Probably, it is more natural to state Th.~\ref{thmu:q^2-s} for
  blocking sets.  In this form, it says, that the size of a
  $(k,s)$-blocking multiset, $q+1\leq s\leq 2q-1$, in $\PHG(R_R^3)$ is at
  least $s(q^2+q)-q^2+q$. The result from the corollary can be stated
  as follows: For a $(k,2q)$-blocking multiset $\mset{K}$ 
  we have $k\geq 2q^3+q^2+q-1$;
  in case of equality $\mset{K}$ is projective (i.e.\ a set),
  and every point neighbour class $[x]$ contains exactly $2q-1$
  points. Moreover, the points in $[x]$ form a blocking set in $\AG(2,q)$.
\end{remark}

Fianlly, putting the bounds from Theorems \ref{thml:s=q+1},
\ref{thml:s=q+2..2q}, \ref{thmu:q^2-s} together we obtain the
following rather narrow intervals for the possible values of
$\maxsize_n(R)$.

\begin{corollary}
  \label{cor:s=q+1..2q}
  For $n\in\{q^2-q,q^2-q+1,\dots,q^2-1\}$ we have
  \begin{equation*}
    \maxsize_n(R)\in
    \begin{cases}
      \left[q^4-q^2-q,q^4-2q\right]&\text{if $n=q^2-1$},\\
      \left[q^4-q^3-2q,q^4-q^3-q+1\right]&\text{if $n=q^2-q$},\\
      \left[q^4-q^3+(t-1)q^2-2q,q^4-q^3+(t-1)q^2+(t-1)q\right] &\text{otherwise},\\
    \end{cases}
  \end{equation*}
  where in the last case $n=q^2-q+t$ with $1\leq t\leq q-2$.
\end{corollary}

\subsection{Refined Analysis of the General Bound}\label{ssec:newbound}

In this subsection, we present 
improvements on some of the upper bounds in the tables from 
\cite{it:kent,it:deadfin,ivan-silvia04,it:chain_survey,it:codes_and_hjelmslev_spaces}.

Our first improvement is based on a refinement of the bound in
Fact~\ref{cor:bound}. For a multiset $\mset{F}$ in $\AG(2,q)$ and a
transversal $\mathcal{T}=\{L_1,\dots,L_{q+1}\}$ of $\AG(2,q)$ (i.e., a
set of $q+1$ pairwise non-parallel lines) we set
$\mset{F}(\mathcal{T})=\sum_{i=1}^{q+1}\mset{F}(L_i)$. For integers
$u,n$ with $0\leq n\leq q^2+q$, $0\leq u\leq\maxsize'_n(q)$ we define
$\ell_q(u,n)$ as the minimum of the numbers
\[\max\bigl\{\mset{F}(\mathcal{T});\mathcal{T}\text{ transversal of }
\AG(2,q)\bigr\},\] 
taken over all projective $(u,n)$-arcs $\mset{F}$
in $\AG(2,q)$, and we set $\ell_q(u)=\ell_q(u,q)$ (same minimum, taken
over all point sets $\mset{F}$ in $\AG(2,q)$ with
$\abs{\mset{F}}=u$). Clearly, $\mset{F}(\mathcal{T})$ is
  maximum iff $\mathcal{T}$ contains a line of maximum multiplicity
  in each parallel class of $\AG(2,q)$.
% \begin{equation}
%   \label{eq:l_q}
%   \begin{aligned}
%   \ell_q(u,n)&=\min_{\substack{\mset{k}\\\text{$\mset{k}$ $(k,n)$-arc}}}
%     \max_{\mathcal{T}}\mset{k}(\mathcal{T}),\\
%   \ell_q(u)&=\min_{\substack{\mset{k}\\\abs{\mset{k}}=u}}
%   \max_{\mathcal{T}}\mset{k}(\mathcal{T}),
% \end{aligned}
% \end{equation}
With this definition we obviously have $\ell_q(u,n)=\ell_q(u)$ for all
$n\geq q$, and defining $\ell_q(u,n)$ for $n>q$ is purely a matter of
convenience (see (\ref{eq:Mqn}).

The numbers $\ell_q(u)$, let alone $\ell_q(u,n)$, are unknown in
general. The function $\ell_q(u,n)$ is weakly increasing (decreasing)
in $u$ (resp., in $n$), and we have the bounds
\begin{equation*}
  \label{eq:ell_qbounds}
  \ell_q(u,n)\geq\ell_q(u)\geq u+q\;(u\geq 1)
  \quad\text{and}\quad\ell_q(u)\leq u+2q-1.
\end{equation*}
The lower bound follows from considering line pencils $\mathcal{T}$
through any point of $\mset{F}$, which always return
$\mset{F}(\mathcal{T})=q+\abs{\mset{F}}$. The upper bound follows from
Theorem~\ref{thm:ell_q}(\ref{thm:ell_q:2}) below and
  $\ell_q(u)\leq\ell_q\bigl(q\lceil u/q\rceil\bigr)=q\lceil
  u/q\rceil+q\leq u+q-1+q$.
 
% For our present purpose it is sufficient to know $\ell_q(u)$ (see
% Theorem~\ref{thm:newbound} below). 
We will now collect several explicit
general results on $\ell_q(u,n)$.
\begin{theorem}
  \label{thm:ell_q}
  \begin{enumerate}[(i)]
  \item\label{thm:ell_q:1} $\ell_q(u,n)=u+q$ for $1\leq
    u\leq\min\{n,q\}$;
  \item\label{thm:ell_q:2} $\ell_q(sq)=sq+q$ for $1\leq s\leq
    q$;\footnote{Question: The computations for $q\leq 5$ suggest that
      even $\ell_q(sq-t)=sq+q$ for $1\leq s\leq q$ and $0\leq t\leq
      s-1$ may be true. For $s\in\{q-1,q\}$ this is in fact true and
      stated below in (\ref{thm:ell_q:6}) and
      (\ref{thm:ell_q:7}). Contrary to what these results suggest,
      the general case does not follow from Bruen's
      bound $k\geq(r+1)q-r$ for $(k,r)$-blocking sets
      in $\AG(2,q)$ (cf.\ \cite[Th.~2.1]{bruen92} and the
      related \cite{blokhuis94a}).}
  \item\label{thm:ell_q:3} $\ell_q(q+1,n)=2(q+1)$ for $n\geq
    2$;\footnote{More generally, 
      $\ell_q(sq+t)=(s+1)(q+1)$ holds for integers $s$, $t$ with
      $0\leq s\leq q-1$, $0\leq t\leq s+1$ if and only if there exists
      an $(sq+t,s+1)$-arc in $\AG(2,q)$. This link with the
      maximal arc problem indicates why it is difficult to determine
      the numbers $\ell_q(u)$ in general.}
  \item\label{thm:ell_q:4} $\ell_q(q+2,n)=2(q+1)$ for
    $q\equiv0\pmod{2}$ and $n\geq 2$;\footnote{Here the generalization
      to maximal arcs is valid: For $q=2^r$ and $0\leq\rho\leq r$ we
      have
      $\ell_q\bigl((2^\rho-1)2^r+2^\rho\bigr)=2^\rho(2^r+1)$. This is
      immediate from the proof for $\rho=1$.}
  \item\label{thm:ell_q:5} $\ell_q(q^2-2q+1)=q^2-1$.
  \item\label{thm:ell_q:6} $\ell_q(u)=q^2$ for $q^2-2q+2\leq u\leq q^2-q$;
  \item\label{thm:ell_q:7} $\ell_q(u)=q^2+q$ for $q^2-q+1\leq u\leq q^2$.
  \end{enumerate}
\end{theorem}

\noindent {\it Proof.}
  (\ref{thm:ell_q:1}) is realized by $u$ collinear points in
  $\AG(2,q)$, while (\ref{thm:ell_q:2}) is realized by 
  % $s-1$ parallel lines
  % and $q-t$ further points on a line of the same parallel
  % class.\footnote{In the extreme case $s=1$ this set consists of one
  %   line.} Now let $\mset{k}$ be a set of $sq-t$ points in
  % $\AG(2,q)$. Then, since $sq-t\geq(s-1)q+1$, every parallel class
  % contains a line $T$ of multiplicity $\mset{k}(T)\geq s$. Hence
  % $\mset{k}(\mathcal{T})<sq+q$ is only possible if there exists no
  % line of multiplicity $q$. However, by a result of
  % Bruen~\cite{bruen92}) the complementary point set of $\mset{k}$, a
  % set of size $(q-s)q+t\leq(q-s)q+s-1=(q-s+1)q-(q-s)-1$, cannot be a
  % $q-s$-fold blocking set.
  $s$ parallel lines in $\AG(2,q)$. Parts (\ref{thm:ell_q:3}) and
  (\ref{thm:ell_q:4}) follow from standard results on $(q+1,2)$- and
  $(q+2,2)$-arcs together with the lower bound
  $\ell_q(u)\geq(q+1)\lceil u/q\rceil$. This lower bound also gives
  the inequality $\geq$ in (\ref{thm:ell_q:5}) and
  (\ref{thm:ell_q:7}).  (\ref{thm:ell_q:5}) is realized by the
  complement of two non-parallel lines (or any other
  blocking set of size $2q-1$). Since
  $\AG(2,q)$ contains no blocking set of size $<2q-1$ (cf.\
  \cite{jamison77,brouwer-schrijver78}), we have
  $\ell_q(q^2-2q+2)=q^2$. The rest follows from monotonicity of
  $\ell_q(u)$ and $\ell_q(q^2-q)=q^2$, $\ell_q(q^2)=q^2+q$.
\qed

The bound in Fact~\ref{cor:bound} takes the stronger
form $\maxsize_n(R)\leq M_{q,n}$ with
\begin{equation}
  \label{eq:Mqn}
  M_{q,n}=\max_{1\leq u\leq\maxsize'_n(q)}
  \min\bigl\{u(q^2+q+1),q(q+1)n-q\cdot\ell_q(u,n)+u\bigr\}.
\end{equation}
Theorem~\ref{thm:newbound} below gives a rather explicit, albeit technical,
formula for $M_{q,n}$. In the proof (and from the formula)
we will see that $M_{q,n}$ depends only 
on $\ell_q(u)$ and not on the more general numbers $\ell_q(u,n)$.
Subsequently (in Theorem~\ref{thm:ell_2345}) we will determine $\ell_q(u)$ for 
$q\leq 5$ completely and use Theorem~\ref{thm:newbound} to compute the bounds
$M_{q,n}$ for the tables in Section~\ref{sec:tables}.
% which also restricts the maximum to a much smaller range. 
% We remind the reader that the numbers $\maxsize_n(R)$ are known
% for $n=0,1$ (trivially) and $n\geq q^2$ (cf.\
% Section~\ref{ssec:genconst}).

For the proof of Theorem~\ref{thm:newbound}
will need a lemma on the monotonicity of the second term in
Eq.~\eqref{eq:Mqn}.
%\footnote{In fact we need it only with $\ell_q(u)$
%  in place of $\ell_q(u,n)$, but the general case is useful 
%for deriving sharper bounds than $M_{q,n}$ in particular cases, which
%is the goal of the next subsection.}

\begin{lemma}
  \label{lma:newbound}
  The function $f(u)=q(q+1)n-q\cdot\ell_q(u,n)+u$, defined for $0\leq
  u\leq\maxsize'_n(q)$,
  has the following property: If $f(u+1)<f(u)$ or, equivalently,
  $\ell_q(u+1,n)>\ell_q(u,n)$ then $f(v)<f(u)$ for
  all $v>u$.
  %\footnote{Here we are tacitly assuming that $u$, $u+1$, $v$ belong
  %  to the domain of $f$.}
\end{lemma}

\noindent {\it Proof.}
  For $0\leq u<\maxsize'_n(q)$ we have either $f(u+1)=f(u)+1$ (when
  $\ell_q(u+1,n)=\ell_q(u,n)$) or $f(u+1)\leq f(u)-(q-1)$ (when
  $\ell_q(u+1,n)>\ell_q(u,n)$). Suppose that the lemma is false and
  $u$ is the largest counterexample. Then necessarily
  $\ell_q(u+1,n)=\ell_q(u+2,n)=\dots=\ell_q(u+q,n)$,
  since any $u'\in\{u+1,u+2,\dots,u+q-1\}$ satisfying
  $\ell_q(u'+1,n)>\ell_q(u',n)$ would be a larger
  counterexample\footnote{Note that $u+1\leq u'\leq u+q-1$
    implies $f(u')\leq f(u)-(q-1)+u'-u-1<f(u)$, and hence by
    assumption there exists $v\geq u+q$ such that $f(v)\geq f(u)>f(u')$.}
  
  Now consider a $(u+q,n)$-arc $\mset{F}$ in $\AG(2,q)$ and a
  transversal $\mathcal{T}$ such that
  $\ell_q(u+q,n)=\max_{\mathcal{T}}\mset{F}(\mathcal{T})$. The
  multiplicity $\mset{F}(L)$ must be constant on parallel classes of
  lines. Otherwise we would be able to strictly decrease
  $\max_{\mathcal{T}}\mset{F}(\mathcal{T})$ by removing $q-1$ points
  of $\mset{F}$, contradicting $\ell_q(u+1,n)=\ell_q(u+q,n)$. Hence
  $u+q=qs$, say, is a multiple of $q$ and $\mset{F}$ meets every line
  of $\AG(2,q)$ in the same number $s\geq 1$ of points. It is
  well-known that this implies $s=q$, $u+q=q^2$, and hence also $n\geq
  q$.
%  \footnote{For example, this follows from the identity 
%  	$\sum_{p\in L}\mset{F}(L)=q\*\mset{F}(p)+\abs{\mset{F}}$, which holds for
%    every point $p$ of $\AG(2,q)$ and shows that $\mset{F}(p)$ is
%    independent of $p$.}  
But in this particular case we have
  $\ell_q(q^2-q,n)=\ell_q(q^2-q)=q^2$, 
  $\ell_q(q^2-q+1,n)=\ell_q(q^2-q+1)=q^2+q\geq\ell_q(q^2-q)+2$, so
  that certainly $f(q^2)<f(q^2-q)$ and $u=q^2-q$ is not a
  counterexample. This contradiction completes the proof of the lemma.
\qed

\begin{theorem}
  \label{thm:newbound}
  For $2\leq n\leq q^2+q$ let
  $u_0=n-\left\lceil\frac{n+2q-1}{q+2}\right\rceil$,
  $u_1=n-\left\lfloor\frac{n+q}{q+2}\right\rfloor$ and $t\geq 0$ the
  largest integer satisfying $u_1+t\leq q^2$ and
  $\ell_q(u_1)=\ell_q(u_1+1)=\dots=\ell_q(u_1+t)$. Then the following holds:
  \begin{enumerate}[(i)]
  \item If $n\equiv 2,3,4,5\pmod{q+2}$ then $u_1=u_0+1$ and
    \begin{equation*}
      M_{q,n}=\max\bigl\{u_0(q^2+q+1),q(q+1)n-q\cdot\ell_q(u_1)+u_1+t\bigr\}.
    \end{equation*}
  \item If $n\not\equiv 2,3,4,5\pmod{q+2}$, then $u_1=u_0+2$ and
    \begin{equation*}
      M_{q,n}=\max\bigl\{u_0(q^2+q+1),M,q(q+1)n-q\cdot\ell_q(u_1)+u_1+t\bigr\},
    \end{equation*}
    where
    $M=\min\bigl\{(u_0+1)(q^2+q+1),q(q+1)n-q\cdot\ell_q(u_0+1)+u_0+1\bigr\}$.
  \end{enumerate}  
  % \begin{equation*}
  %   M_{q,n}=\max_{u_0\leq u\leq u_1}\min
  %   \bigl\{u(q^2+q+1),q(q+1)n-q\cdot\ell_q(u)+u\bigr\},
  % \end{equation*}
  % where $u_0=n-\left\lceil\frac{n+2q-1}{q+2}\right\rceil$ and
  % $u_1$ is the smallest multiple of $q$ greater or
  %   equal to $n-\left\lfloor\frac{n+q}{q+2}\right\rfloor$.
\end{theorem}
Note that $0\leq u_0<u_1\leq n\leq q^2$, so that all expressions
occuring in the theorem are well-defined.

\noindent {\it Proof.}
  % The function $(u,n)\mapsto\ell_q(u,n)$ is non-decreasing in $u$ and
  % non-increasing in $n$.
  We show first that replacing $\ell_q(u,n)$ by $\ell_q(u)$ in
  \eqref{eq:Mqn} das not change $M_{q,n}$. This is trivial for $n\geq
  q$, so let $n<q$. Then for $u\leq n$ we also have
  $\ell_q(u,n)=u+q=\ell_q(u)$, and for $u=n$ the minimum in
  \eqref{eq:Mqn} equals
  $q(q+1)n-q(n+q)+n=n(q^2+q+1)-q(n+q)$. For $u>n$ we
  have $q(q+1)n-q\*\ell_q(u,n)+u\leq
  q(q+1)n-q(u+q)+u<q(q+1)n-q(n+q)+n$, and similarly with $\ell_q(u)$
  in place of $\ell_q(u,n)$. Since $u\mapsto u(q^2+q+1)$ is 
  increasing, this clearly implies $M_{q,n}=\max_{1\leq u\leq n}\min
    \bigl\{u(q^2+q+1),q(q+1)n-q\*\ell_q(u)+u\bigr\}=\max_{1\leq u\leq q^2}\min
    \bigl\{u(q^2+q+1),q(q+1)n-q\*\ell_q(u)+u\bigr\}$ as asserted.
  
  % Writing $f_1(u)=u(q^2+q+1)$,
  % $f_2(u)=q(q+1)n-q\cdot\ell_q(u)+u$, we have that $f_1$ is
  % strictly increasing while $f_2$ is ``mostly'' strictly
  % decreasing.\footnote{We shall make this precise below.}
  Now we can prove the theorem. Let us write $f_1(u)=u(q^2+q+1)$,
  $f_2(u)=q(q+1)n-q\cdot\ell_q(u)+u$.
  From $u+q\leq\ell_q(u)\leq u+2q-1$ we
  have $q(q+1)n-2q^2+q-(q-1)u\leq f_2(u)\leq q(q+1)n-q^2-(q-1)u$. The
  function $f_1$ (considered as a real-valued function) meets
  $g_0(x)=q(q+1)n-2q^2+q-(q-1)x$ in
  $x_0=\frac{(q+1)n-2q+1}{q+2}=n-\frac{n+2q-1}{q+2}$ and
  $g_1(x)=q(q+1)n-q^2-(q-1)x$ in
  $x_1=\frac{(q+1)n-q}{q+2}=n-\frac{n+q}{q+2}$. Since $f_1$ is
  increasing and $g_0$ is (strictly) decreasing, we have $f_1(u)\leq
  g_0(u)\leq f_2(u)$ for $u\leq\lfloor
  x_0\rfloor=n-\left\lceil\frac{n+2q-1}{q+2}\right\rceil=u_0$ and
  hence $\max_{0\leq u\leq
    u_0}\min\bigl\{f_1(u),f_2(u)\bigr\}=f_1(u_0)=u_0(q^2+q+1)$.
%\footnote{Changing
%    the range of the max from $\{1,\dots,q^2\}$ to $\{0,\dots,q^2\}$
%    does not change the value of $M_{q,n}$. It is, however, necessary
%    for covering the case $n=2$.}
  Similarly, $f_2(u)\leq g_1(u)\leq f_1(u)$ for $u\geq\lceil x_1\rceil
  =n-\left\lfloor\frac{n+q}{q+2}\right\rfloor=u_1$.
  By definition of $t$ we have $f_2(u_1+t')=f_2(u_1)+t'$ for $0\leq
  t'\leq t$ and by Lemma~\ref{lma:newbound} we have
  $f_2(v)<f_2(u_1+t)=f_2(u_1)+t$ for all
  $v>u_1+t$. %(Here we use the fact derived in the previous paragraph.)
  Hence $\max_{u\geq u_1}\min\bigl\{f_1(u),f_2(u)\bigr\}=\max_{u\geq u_1}f_2(u)
  =f_2(u_1)+t=q(q+1)n-q\*\ell_q(u_1)+u_1+t$.

  Up to this point we have shown that
  $M_{q,n}=\max\bigl\{f_1(u_0),M',f_2(u_1)+t\bigr\}$, where
  $M'=\max_{u_0<u<u_1}\min\bigl\{f_1(u),f_2(u)\bigr\}$. Since the
  numbers $\frac{n+q}{q+2}$ and $\frac{n+2q-1}{q+2}$ differ by
  $0<\frac{q-1}{q+2}<1$, it is clear that
  $u_1-u_0\in\{1,2\}$.\footnote{In other words, there is at most one
    integer strictly between $u_0$ and $u_1$.} We have
  $x_1\in\Z$ iff $n\equiv 2\pmod{q+2}$, and $x_0\in\Z$ iff $n\equiv
  5\pmod{q+2}$. This and a check of signs 
  in the expressions defining $u_0$, $u_1$ gives that the open
  interval $]u_0,u_1[$ contains no integer iff $n\equiv
  2,3,4,5\pmod{q+2}$ and contains exactly one integer (the number
  $u_0+1$) if $n\not\equiv 2,3,4,5\pmod{q+2}$. Hence $M_{q,n}$ has the
  form stated in (i), (ii), and the proof of the theorem is complete.
  %In the latter case $M'=M$ is the number claimed in (ii).
\qed

\begin{remark}
  \label{rmk:newbound-omit}
  Using a more delicate analysis, one can show that the first term
  $u_0(q^2+q+1)$ in the formulas for $M_{q,n}$ in
  Th.~\ref{thm:newbound}(i) and (ii)
  may in fact be omitted. For this it suffices to show
  $f_2(u_0+1)>u_0(q^2+q+1)$. Using the notation of the previous proof,
  if $x_0$ is not an integer then
  $g_0(x_0)=f_1(x_0)\geq\left(u_0+\frac{1}{q+2}\right)(q^2+q+1)$
  and $f_2(u_0+1)\geq g_0(u_0+1)=g_0(x_0)-(q-1)(u_0+1-x_0)
  \geq
  u_0(q^2+q+1)+\frac{q^2+q+1}{q+2}-(q-1)\frac{q+1}{q+2}
  =u_0(q^2+q+1)+1$.\footnote{This inequality is sharp in the case
    $n=q^2$, where $u_0=q^2-q$, $u_1=u_0+1=q^2-q+1$,
    $f_2(u_0+1)=q(q+1)q^2-q(q^2+q)+q^2-q+1=(q^2-q)(q^2+q+1)+1$.}

  Now assume $x_0$ is an integer, so that $n\equiv 5\pmod{q+2}$, 
  $u_0=x_0$, and
  $g_0(u_0+1)=u_0(q^2+q+1)-(q-1)$. The function $\ell_q(u)$ satisfies
  the stronger upper bound $\ell_q(u)\leq u+2q-2$ for all $u\neq
  1\pmod{q}$.\footnote{This is clear from the proof of $\ell_q(u)\leq
    u+2q-1$.} Hence $f_2(u_0+1)\geq u_0(q^2+q+1)+1$ except possibly
  for $u_0\equiv0\pmod{q}$. Substituting $n=(q+2)k+5$ in
  $u_0=n-\frac{n+2q-1}{q+2}$ gives $(q+2)k+5-k-2\equiv 0\pmod{q}$,
  i.e.\ $k\equiv -3\pmod{q}$. We conclude that $k=q-3$,
  $n=q^2-q-1$, $u_0=q^2-2q$. However, in this case
  $\ell_q(u_0+1)=q^2-1=u_0+1+2q-2$ (cf. Th.~\ref{thm:ell_q:5})
  and $f_2(u_0+1)=u_0(q^2+q+1)+1$. Hence $f_2(u_0+1)>u_0(q^2+q+1)$
  holds in all cases.
\end{remark}  

\begin{remark}
  \label{rmk:newbound-explicit}
    For small values of $n$ the known explicit formulas for
    $\ell_q(u)$, $u\leq q+1$ ($u\leq q+2$ for even $q$) can be used to
    determine $M_{q,n}$ explicitly:
    \begin{equation}
      \label{eq:Mqnexpl}
    \begin{aligned}
      M_{q,n}&=(n-1)q^2+q+(n-1)\quad\text{for $2\leq n\leq q+1$},\\
      M_{q,q+2}&=
      \begin{cases}
        q^3+q^2+q+1&\text{if $q\equiv 1\pmod{2}$},\\
        q^3+q^2+q+2&\text{if $q\equiv 0\pmod{2}$},
      \end{cases}\\
      M_{q,q+3}&=\begin{cases}
        q^3+2q^2+2q+1&\text{if $q\equiv 1\pmod{2}$},\\
        q^3+2q^2+2q+2&\text{if $q\equiv 0\pmod{2}$},
        \end{cases}\\
      M_{q,q+4}&=q^3+3q^2+3q+2\quad\text{for $q\equiv 0\pmod{2}$}.
    \end{aligned}
  \end{equation}
\end{remark}

\begin{corollary}
  \label{cor:newbound}
  For $(n,k)$-arcs in $\PHG(2,R)$ we have the general upper bound
  \begin{equation*}
    \maxsize_n(R)
    \leq\left\lfloor\frac{(q+1)n-q}{q+2}\cdot(q^2+q+1)\right\rfloor.
  \end{equation*}
\end{corollary}

\noindent {\it Proof.}
  The minimum in Eq.~\eqref{eq:Mqn} is upper-bounded by
  \[\min\bigl\{u(q^2+q+1),q(q+1)n-q(q+u)+u\bigr\},\] 
  which takes its maximum
  for $x_1=\frac{(q+1)n-q}{q+2}$.
\qed

Now we turn to the determination of the numbers $\ell_q(u)$ for $q\leq
5$.
\begin{theorem}
  \label{thm:ell_2345}
  For small $q$ the values of $\ell_q(u)$ are those listed in the
  following table.
  {\small
  \begin{equation*}
    \setlength{\arraycolsep}{3pt}
    \begin{array}{c|ccccccccccccccccccccccccc}
      \multicolumn{26}{c}{\ell_q(u)}\\[.5ex]
      q\backslash u&1&2&3&4&5&6&7&8&9&10&11&12&13&14&15
      &16&17&18&19&20&21&22&23&24&25\\\hline
      2&3&4&6&6\\
      3&4&5&6&8&9&9&12&12&12\\
      4&5&6&7&8&10&10&12&12&15&16&16&16&20&20&20&20\\
      5&6&7&8&9&10&12&14&14&15&15&18&19&20&20&20
      &24&25&25&25&25&30&30&30&30&30
    \end{array}
  \end{equation*}
}
\end{theorem}

\noindent {\it Proof.}
  Most of the values are covered by Th.~\ref{thm:ell_q}.
  Only the following values need extra justification. 
  % For this we recall that
  % the spectrum $(a_i)_{0\leq i\leq q+1}$
  % of a set $\mset{k}$ 
  % of $u$ points in $\AG(2,q)$ satisfies the equations $\sum
  % a_i=q^2+q$, $\sum ia_i=(q+1)u$, $\sum i(i-1)a_i=u(u-1)$.
  \smallskip
  
  $\bullet$ $\ell_4(7)=12$: Suppose $\abs{\mset{F}}=7$. If $\mset{F}$ has a
  $4$-line then
  \[\ell=\max\bigr\{\mset{F}(\mathcal{T});\mathcal{T}\text{ transversal
    of }\AG(2,q)\bigr\}\geq 4+4\cdot 2=12.\]  
Otherwise every point of
  $\mset{F}$ is on at least one $3$-line. Since $\mset{F}$ is not the
  union of parallel $3$-lines, there exist $3$-lines in at least two
  directions and hence $\ell\geq 12$.
%  \footnote{It goes without saying
%    that $\ell_4(7)\leq 12$. This follows from $\ell_4(8)=12$ and the
%    monotonicity of $\ell_4$; a set of $7$ points contained in two
%    parallel lines has $\ell=12$.}
  
  $\bullet$ $\ell_5(7)=14$: $\ell=14$ is realized by $7$ points contained in a
  $2\times 4$-grid.\footnote{By an \emph{$r\times s$-grid} in
    $\AG(2,q)$, $1\leq r,s\leq q$, we mean any set of $rs$ points
    obtained as the intersection points of $r$ parallel lines with $s$
    lines from another parallel class. In our case it actually
    suffices to take $4$ points on one line together with $3$ points
    on another parallel line. Such a set need not form a $2\times
    4$-grid. The equality $\ell=14$ is also realized by a $(7,3)$-arc
    (the union of an oval and one of its external points), showing that
    $\ell_5(7,3)=\ell_5(7,4)=14$.}  Now let $\abs{\mset{F}}=7$. If
  there is a $4$-line or a $5$-line we have $\ell\geq 14$. Otherwise,
  % we must have $3$-lines in at least two directions (and so again
  % $\ell\geq 14$), since every parallel class contains a $1$-line or a
  % $3$-line and points on $1$-lines must also be on at least one
  % $3$-line.
  if there is only one $3$-line $L$ then removing a point in
  $\mset{F}\cap L$ produces an oval and hence a
  contradiction.
  %\footnote{Points outside an oval are on at least two
  % $2$-lines.}
  If there exist two parallel $3$-lines $L_1$ and $L_2$ then the point
  of $\mset{F}$ outside $L_1\cup L_2$ is on a $1$-line parallel to
  $L_1$, $L_2$ and hence $\ell\geq q+u+2=14$.
%  \footnote{One can also
%    use the fact that a point $p\in\mset{k}$
%    on a $1$-line must be on at least one
%    $3$-line, since $\sum_{p\in L}\mset{k}(L)=12$.}

  $\bullet$ $\ell_5(8)=14$: $\ell=14$ is realized by a $2\times
  4$-grid. Using the standard representation of 
  $\AG(2,5)=\AG(\F_5^2)$, $\ell=14$ is also realized by
  $\{0,1,2\}\times\{0,1,2\}\setminus\{(1,1)\}$ (a special punctured
  $3\times 3$-grid), showing that $\ell_5(8,3)=\ell_5(8,4)=14$.

  $\bullet$ $\ell_5(9)=15$: $\ell=15$ is realized by $9$ points contained in a
  $2\times 4$-grid. Now suppose $\abs{\mset{F}}=9$. The inequality
  $\ell\geq q+u+1=15$ clearly holds unless all lines of every parallel
  class contain the same number of points of $\mset{F}$. This means
  that every line is a $3$-line, which is impossible.

  $\bullet$ $\ell_5(12)=19$: Since the maximum number of points on a $(k,3)$-arc
  in $\AG(2,5)$ is $11$, we have $\ell_5(12)\geq 19$. The equality
  $\ell=19$ is realized by a $3\times 4$-grid.
  
  $\bullet$ $\ell_5(13)=20$: Let $\abs{\mset{F}}=13$. If $\mset{F}$ has a
  $5$-line then $\ell\geq 20$. Otherwise, since there is no
  $(12,3)$-arc in $\AG(2,5)$, there exist at least two $4$-lines. If
  these are non-parallel, we have $\ell\geq 2\cdot 4+4\cdot 3=20$, and
  if they are parallel, we have a $2$-line or a $1$-line in the same
  parallel class and hence again $\ell\geq q+u+2=20$.

  $\bullet$ $\ell_5(14)=20$: This follows from $\ell_5(13)=20$ (and $\ell_5(15)=20$).
\qed

Combining Theorems~\ref{thm:newbound} and~\ref{thm:ell_2345} yields
improved upper bounds for $\maxsize_n(R)$ for several pairs $(q,n)$ in
the range considered (cf. the tables in Section~\ref{sec:tables}).

\subsection{Special Non-Existence Proofs for $q=3,4,5$}
\label{ssec:spebounds}

In the remainder
we will further improve the upper bounds for $q=3,4,5$ and particular
small values of $n$. Several of the improvements rely on the following Lemma.

\begin{lemma}
\label{lma:upper_bound}
Let $\mset{K}$ be a $(k,n)$-arc in $\PHG(2,R)$, and let $[L]$ be a neighbour 
class of lines in $\PHG(2,R)$ such that $\mset{K}\bigl([x]\bigr)\leq
n$ for all point neighbour classes $[x]$ incident with $[L]$. Then
\begin{equation*}
  \mset{K}\left([L]\right)\leq\maxsize_n(q)-(n-m),\quad\text{where}
    \quad m=\max\bigl\{\mset{K}\bigl([x]\bigr);x\in L\bigr\}
\end{equation*}
is the maximum multiplicity of a point class on $[L]$.
% \begin{equation}
% \label{eq:ub}
% \mset{K}\left([L]\setminus[x]\right)\le 
% \left\{\begin{array}{ll}
% m_n(q)-\mset{K}([x]) & \text{ if } \mset{K}([x])\le n, \\
% m_n(q)-n             & \text{ if } \mset{K}([x])\le n.
% \end{array}\right.
% \end{equation}
% Here $m_n(q)$ denotes the maximal number of points contained in a 
% $(\kappa,n)$-arcs in $\PG(2,q)$.
\end{lemma}

\noindent {\it Proof.}
  The restriction of $\mset{K}$ to $[L]$ defines a
  $\bigl(\mset{K}(L),n\bigr)$-arc in the punctured
  projective plane induced on
  $[L]$. Assigning multiplicity $n-m$ to the missing point $p_\infty$
  we obtain a $\bigl(\mset{K}([L])+n-m,n\bigr)$-arc in $\PG(2,q)$. The
  result follows immediately.
\qed

Note that the inequality of the lemma may also be stated in the form
$\mset{K}\bigl([L]\setminus[x]\bigr)\leq\maxsize_n(q)-n$ (which is
independent of $m$), provided that $[x]$ has maximum multiplicity
among the point neighbour classes incident with $[L]$.

In (\cite{it:kent}) it was proved that $m_3(R_R^3)\le2q^2$
for all $q\ge4$. Below we improve this bound.

\begin{theorem}
  \label{thm:n=3}
  \begin{equation}
    \label{eq:n=3}
    \maxsize_3(R)\leq
    \begin{cases}
      2q^2-2q+4 & \text{ for } q \text{ odd,} \\
      2q^2-q+4 & \text{ for } q \text{ even.}
    \end{cases}
  \end{equation}
  for every prime power $q\ge5$.
\end{theorem}
% \nocite{barlotti56b,ball-blokhuis-mazzocca97,ball-blokhuis98,ball-hirschfeld05}

 \noindent {\it Proof.}
  Let $\mset{K}$ be a $(k,3)$-arc of maximum cardinality in
  $\PHG(2,R)$, $\abs{R}=q^2\geq 25$, and $[x]$
  a neighbour class of points with
  $u:=\mset{K}([x])=\max\bigl\{\mset{K}([y]);y\in\points\bigr\}$. By
  Cor.~\ref{cor:bound} we have
 \[ k\leq 2q^2+3q-(q-1)u.\]
 If we assume that $u\geq 6$, then $k\leq 2q^2-3q+6\leq 2q^2-2q+1$.

Assume now that $u=5$. By considering the possible types of $(5,3)$-arcs in 
$[x]\cong\AG(2,q)$, we find %(using the notation of Fact~\ref{thm:bound})
that $\sum u_i\geq q+6$ in each case and
hence by Fact~\ref{thm:bound}
\begin{equation*}
  k\leq 3q^2+3q-q\sum_{i=1}^{q+1}u_i+5\leq 2q^2-3q+5\leq 2q^2-2q.
 \end{equation*}

Next consider the case $u=4$. If three of the points in
$[x]\cap\supp\mset{K}$ are collinear, we have again $\sum u_i\ge q+6$
and hence $k\leq 2q^2-3q+4\leq 2q^2-2q-1$.
Assume there exists a class $[x]$ with $\mset{K}([x])=4$. 
There are two possibilities.
 \begin{enumerate}[i.]
 \item three of the points in $[x]$ are collinear, or
 \item no three of the points in $[x]$ are collinear.
 \end{enumerate}

 i.\ Let $x,y,z,t$ be four points in $[x]\supp\mset{K}$ and
 let $y,z,t$ be collinear.
 In this case the neighbour class of lines $[\overline{yz}]$
 contains no points of $\mset{K}$ (apart from those in $[x]$).
 The  classes $[\overline{xy}]$, $[\overline{xz}]$, $[\overline{xt}]$
 are different and the multiplicity of these classes is at most
 $q+4$. Finally, Lemma~\ref{lma:upper_bound} implies that
 for all the remaining classes of lines $[L]$ thorugh $[x]$
 we have $\mset{K}([L]\setminus[x])\le 2q-2$. Hence we obtain
 \[k\le 4 +3\cdot q+ (q-3)\cdot(2q-2)=2q^2-5q+10<2q^2-2q+3.\]

 Otherwise the four points in $[x]\cap\supp\mset{K}$ form a
 quadrangle $Q=\{x,y,z,t\}$, and there are two subcases to distinguish:

 (i) $q$ is odd. The points in $Q$ determine at least four
 directions. If $[L]=[\overline{xy}]$, say, is determined by $Q$, then
 $\mset{K}\bigl([L]\setminus[x]\bigr)\leq q$. If $[L]$ is incident with
 $[x]$ but not determined by $Q$, then
 $\mset{K}\bigl([L]\setminus[x]\bigr)\leq 2q$. 

 By the above argument,
 \[k\le 4 +4\cdot q+(q-3)\cdot2q=2q^2-2q+4.\]

 (ii) $q$ is even. Here the same argument applies with the difference 
 that the points in $Q$ can determine exactly three directions. Now
 \[k\le 4 +3\cdot q+(q-2)\cdot2q=2q^2-q+4.\]

 It remains to consider the case $u=3$. Here we will use the refined bound of
 Lemma~\ref{lma:upper_bound}. It is well-known that $\maxsize_3(q)\leq
 2q+1$ for all $q$.
 %\footnote{This follows from the results in
 %  \cite{barlotti56b,ball-blokhuis-mazzocca97} (see also
 %  \cite{ball-blokhuis98,ball-hirschfeld05,hirschfeld98}).}

 If the three points in $[x]\cap\supp\mset{K}$ are collinear then
 \[k\le 3+1\cdot 0+q\cdot(2q-2)=2q^2-2q+3.\]
Otherwise they determine a triangle and 
\[k\le 3 +3\cdot q + (q-2)\cdot(2q-2)=2q^2-3q+7\leq 2q^2-2q+2.\]
\qed

For chain rings of order $q^2=4$ the numbers $\maxsize_n(R)$ were
determined completely in \cite{it:deadfin}. 

Next we briefly review the case $q^2=9$; cf. Table~\ref{tbl:q=3}.  The
numbers $\maxsize_2(\Z_9)=\maxsize(\DN_3)$ were determined in
\cite{it:deadfin}. For $n=3$ the bounds $18\leq\maxsize_3(R)\leq 19$
were derived in \cite{it:deadfin}. Arcs with Parameters $(18,3)$ are
easy to construct, for example we can take six point neighbour classes
forming a minimal blocking set in the quotient plane $\PG(2,3)$ and
place a line segment in tangent direction into each class. The case
$n=3$ was finally completed using an exhaustive computer search
\cite{ivan-silvia04}, which revealed $\maxsize_3(\Z_9)=19$,
$\maxsize_3(\DN_3)=18$.\footnote{Is there a direct proof that the
  existence of a $(19,3)$-arc in $\PHG(2,R)$, $\abs{R}=9$, is
  incompatible with the existence of a proper subplane?}  The
situation for $n=5$ is similar: Here we have $\maxsize_5(\Z_9)=39$,
$\maxsize_5(\DN_3)=38$; cf.~\cite{imt:two_exact_values}. The structure
of a $(39,5)$-arc in $\PHG(2,R)$, $\abs{R}=9$, can be determined
rather explicitly, and the existence of such an arc turns out to be
incompatible with the existence of a proper subplane. The numbers
$\maxsize_8(\Z_9)=\maxsize_8(\DN_3)=60$ are also due to
\cite{imt:two_exact_values}. Only two cases, $n=6$ and $n=7$, remain
open. The standard bound (Fact~\ref{cor:bound} or
Th.~\ref{thm:newbound}) gives $\maxsize_6(R)\leq 52$,
$\maxsize_7(R)\leq 63$. Both bounds can be improved by $1$.  In the
former case this was done in~\cite[Th.~5.4]{it:kent}.  Here we provide
a proof for the latter, which was sated without proof in
\cite{it:deadfin}.  

In what follows we denote the number of point classes of multiplicity
$i$ (relative to a given multiset $\mset{K}$ of points in $\PHG(2,R)$)
by $f_i$. The integers $f_i$ satisfy the identities $\sum
f_i=\abs{\pclasses}=q^2+q+1$ and $\sum
if_i=\sum_{[x]\in\pclasses}\mset{K}\bigl([x]\bigr)=\abs{\mset{K}}$.

\begin{theorem}
  \label{thmu:q=3:n=7}
  $\maxsize_7(R)\leq 62$ for $R=\Z_9$ or $\DN_3$.
\end{theorem}

\noindent {\it Proof.}
  Assume that $\mset{K}$ is a $(63,7)$-arc in
  $\PHG(2,R)$. Fact~\ref{thm:bound} gives that
  $u=\max_{[x]}\mset{K}\bigl([x]\bigr)=6$ and the points of $\mset{K}$ in
  every neighbour class $[x]$ with $\mset{K}\bigl([x]\bigr)=6$ form
  two parallel line segments. Fix such a class $[x]$ and denote by
  $[L_0]$, $[L_1]$, $[L_2]$, $[L_3]$ the line classes incident with
  $[x]$, the class $[L_0]$ being the direction of the line segments in
  $[x]$. Then necessarily $\mset{K}\bigl([L_0]\bigr)=18$,
  $\mset{K}\bigl([L_1]\bigr)=\mset{K}\bigl([L_2]\bigr)
  =\mset{K}\bigl([L_3]\bigr)=21$, and by assigning multiplicity $1$ to
  the missing point $p_\infty$ we can extend the multisets induced by
  $\mset{K}$ in $[L_i]$ to $(19,7)$- resp.\ $(22,7)$-arcs in the
corresponding projective planes, which have order $3$.
%\footnote{The
%  $(22,7)$-arcs are maximal arcs, and the $(19,7)$-arc is maximal
%  subject to having a point of multiplicity $3$.}

As $(22,7)$-arcs in $\PG(2,3)$ are just complements of lines; more
precisely, the multiset $2-\mset{K}\vert_{[L_i]}$ is a
line, and hence every line class $[L_i]$, $i=1,2,3$, has type
$(6,6,6,3)$. Note that the line formed by the point classes of
multiplicity $1$ must pass through $p_\infty$ (see the figure below). 

\begin{figure}[htp!]
	\begin{center}
		\begin{tikzpicture}[line width=1pt, scale=0.4]
		\draw[black] (3,1.5)--(15,1.5);
		\draw[gray] (0,0)--(3,0)--(3,3)--(0,3)--(0,0) [fill=yellow!30];
		\draw[gray] (5,0)--(8,0)--(8,3)--(5,3)--(5,0) [fill=yellow!30];
		\draw[gray] (10,0)--(13,0)--(13,3)--(10,3)--(10,0) [fill=yellow!30];
		\draw[gray] (15,0)--(18,0)--(18,3)--(15,3)--(15,0) [fill=yellow!30];
		
		\draw[gray] (1.5,3) .. controls (1.5,4.5) and (8,5) .. (9,5);
		\draw[gray] (16.5,3) .. controls (16.5,4.5) and (10,5) .. (9,5);
		\draw[gray] (6.5,3) .. controls (6.5,4.5) and (8,4.5) .. (9,5);
		\draw[gray] (11.5,3) .. controls (11.5,4.5) and (10,4.5) .. (9,5);

		\draw[gray] (0.6,2.25) circle (0.2cm) [fill=black];
		\draw[gray] (0.6,1.5) circle (0.2cm) [fill=black];
		\draw[gray] (0.6,0.75) circle (0.2cm) [fill=black];
		\draw[gray] (1.2,2.25) circle (0.2cm) [fill=black];
		\draw[gray] (1.2,1.5) circle (0.2cm) [fill=black];
		\draw[gray] (1.2,0.75) circle (0.2cm) [fill=black];
		\draw[gray] (5.6,2.25) circle (0.2cm) [fill=black];
		\draw[gray] (5.6,1.5) circle (0.2cm) [fill=black];
		\draw[gray] (5.6,0.75) circle (0.2cm) [fill=black];
		\draw[gray] (6.2,2.25) circle (0.2cm) [fill=black];
		\draw[gray] (6.2,1.5) circle (0.2cm) [fill=black];
		\draw[gray] (6.2,0.75) circle (0.2cm) [fill=black];
		\draw[gray] (10.6,2.25) circle (0.2cm) [fill=black];
		\draw[gray] (10.6,1.5) circle (0.2cm) [fill=black];
		\draw[gray] (10.6,0.75) circle (0.2cm) [fill=black];
		\draw[gray] (11.2,2.25) circle (0.2cm) [fill=black];
		\draw[gray] (11.2,1.5) circle (0.2cm) [fill=black];
		\draw[gray] (11.2,0.75) circle (0.2cm) [fill=black];
		\draw[gray] (15.6,2.25) circle (0.2cm) [fill=black];
		\draw[gray] (15.6,1.5) circle (0.2cm) [fill=black];
		\draw[gray] (15.6,0.75) circle (0.2cm) [fill=black];

		\draw[gray] (9,5) circle (4mm) [fill=black];

\draw(1.5,-.8) node{\small {\color{blue} $[x]$}};	
\draw(9.5,6) node{\small {\color{blue} $p_{\infty}$}};	
\draw(20,1.5) node{\small {\color{blue} $[L_i]$}};		
\end{tikzpicture}
\end{center}
\caption{One of the neighbor classes $[L_i]$, $i=1,2,3$}
\end{figure}
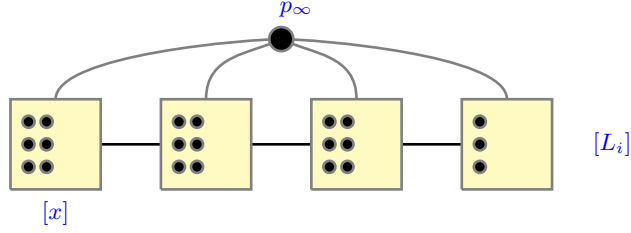

Hence we have $f_6\geq 7$.
Next we show that in fact $f_6=7$.
%there cannot be more than $7$ such classes. 
The line segments in two distinct point classes $[x]$ with
$\mset{K}\bigl([x]\bigr)=6$ cannot have the same direction $[L]$,
since this would imply $\mset{K}\bigl([L]\bigr)\leq 16$. Hence there
are at least $f_6$ line classes of multiplicity $18$ (the directions
of the line segments in those point classes). The remaining line
classes have multiplicity at most $21$ (at most $20$ if they are not
incident with a point class of multiplicity $6$). Hence we obtain
\begin{equation}
  \label{eq:f_6}
  252=|\mset{K}|\cdot 4=\sum_{[L]}\mset{K}\bigl([L]\bigr)\leq 18f_6+21(13-f_6),
  \quad\text{i.e.}\quad f_6\leq 7.
\end{equation}
With $f_6=7$ it is now easy to derive the desired contradiction. The
set $S=\bigl\{[x]\in\pclasses;\mset{K}([x])=7\bigr\}$, viewed as a set
of $7$ points in $\PG(2,3)$, has the following property: Every point
of $S$ is incident with exactly three $3$-lines and exactly one 1-line. 
This means that $S$
is isomorphic to $\PG(2,2)$, which is impossible.
\qed

The next case $q^2=16$ is much more involved and a resolution of the
remaining open cases is not to be expected soon. Since the publication
of \cite{it:kent}, however, considerable progress has been made. The
values of $\maxsize_n(R)$ have been determined completely
for $n\in\{2,4,5,9\}$ (and in the case of $\G_4$ also for $n\in\{6,8\}$),
and much sharper bounds are now known for the remaining open cases;
cf.\ Table~\ref{tbl:q=4} in Section~\ref{sec:tables}.

Here we complete the cases $n=5$ (for all three rings) and $n=6$ (for
$\G_4$) by sharpening
the standard bounds $\maxsize_5(R)\leq 72$ to $\maxsize_5(R)\leq 68$,
respectively $\maxsize_8(R)\leq 86$ to $\maxsize_8(R)\leq
84$.\footnote{The latter was already stated (without proof) in
  \cite{it:kent}.}

\begin{theorem}
  \label{thmu:q=4:n=5}
  Let $R$ be one of the rings $\G_4$, $\DN_4$, $\TDN_4$. Then
  $\maxsize_5(R)\leq 68$.
\end{theorem}

\noindent {\it Proof.}
  Let $\mset{K}$ be a $(k,5)$-arc in $\PHG(2,R)$.
  In this case the standard bound (Facts~\ref{thm:bound}, \ref{cor:bound})
  gives $k\leq 72$, and equality can hold
  only if $\mset{K}\bigl([x]\bigr)\leq 4$ for each point neighbour class
  $[x]$ and the points of $\mset{K}$ in each class $[x]$ with
  $\mset{K}\bigl([x]\bigr)=4$ form a line segment or an affine
  subplane (i.e.\ a quadrangle with three pairs of parallel
  sides). In all other cases $k\leq 68$. 
  
  Now we assume $k=69$. Then there exists a point class of multiplicity
  $4$ (in fact at least six such classes), and we
  are left with the following two cases:
  % The two cases are exhaustive, since there exist (at least
  % six) classes $[x]$ with $\mset{K}\bigl([x]\bigr)=4$.
  
  \textsl{Case~1:} $\mset{K}$ contains a line segment. The class $[x]$
  containing a line segment is uniquely determined, since two such
  classes $[x_1]$, $[x_2]$ would give rise to two distinct line classes
  $[L_1]$, $[L_2]$ of multiplicity $\leq 8$ and
  $k\leq 2\cdot 8+(21-9)\cdot 4=64$, a contradiction.
  
  We denote the direction of the line segment in $[x]$ by
  $[L_0]$. Every point class $[y]\neq[x]$ on $[L_0]$ satisfies
  $\mset{K}\bigl([y]\bigr)\leq 1$, since otherwise an application of
  the standard bound to $[y]$ yields a contradiction. For
    example, $\mset{K}\bigl([y]\bigr)=2$ implies $k\leq 2+1\cdot
    6+1\cdot 12+3\cdot 16=66$, since the direction determined by the
    two points in $[y]$ is certainly not $[L_0]$.  
    
Let $[y]$ be  a neighbor class of points that is of multiplicity 1.
Let $[L]$ be a neighbor class  of lines through $[y]$ that is of multiplicity 17,
i.e. it contains apart from $[y]$ four point classes of maximal multiplicity 4. 
The restriction $\mathcal{K}\vert_{\Pi_{[L]}}$ is a $(17,5)$-arc in the projective plane of order four which is extendable to a $(18,5)$-arc and hence does not have points of multiplicity greater than 1. The numeber of such lines is at least 4.
For inctance if $\mathcal{K}([L_0])=5$ these are all lines through $[y]$. If 
$\mathcal{K}([L_0])=8$ there is on such line through every point class on $[L_0]$ different from $[x]$.
  
  The proof can now be finished by counting the number of $2$-points
  in the geometries $\Pi_{[L]}$ for the $16$ line classes $[L]$ not
  incident with $[x]$. Each class $[z]$ outside
  $[L_0]$ determines exactly six such $2$-points
  and there are at least 13 such classes.
  On the other hand there are at most 12 planes  $\Pi_{[L]}$ with 2-points.
  Hence there must be neighbor class of lines that contains at least
  $\displaystyle  {13\cdot 6}{12}>6$ points. This implies the existence of a 
  plane $\Pi_{[L]}$ with  three colliniear 2 points which in turn implies the
  existence of a line in $\PHG(2,R)$ of multiplicity at least 6, a contradiction. 

  \textsl{Case~2}: $\mset{K}$ does not contain a line segment.  In
  this case we first also dispose of (shortened) line segments of
  length $3$ in $\mset{K}$, so that the points of $\mset{K}$ in line
  classes of multiplicity $3$ are necessarily triangles: Suppose
  $\mset{K}\cap[x]$ is a line segment of length $3$ with direction
  $[L_0]$.  Using essentially the same reasoning as in Case~1, all
  point classes $[y]\neq[x]$ on $[L_0]$ have multiplicity $\leq 2$. 
  This leads to a contradiction if a neighbor class of points is of 
  multiplicity 3 with three collinear points. 

  Now we argue as follows: We call a line class $[L]$
  \emph{projective} if $\mset{K}$ induces a projective arc in
  $\Pi_{[L]}$. Non-projective line classes $[L]$ have
  $\mset{K}\bigl([L]\bigr)\leq 16$, since $\mset{K}$ induces, by
  assigning multiplicity $1$ to the missing point $p_\infty$, a
  $\bigl(\mset{K}([L])+1,5\bigr)$-arc in $\Pi_{[L]}$.\footnote{For
    $(l,5)$-arcs in $\PG(2,4)$ containing a multiple point we have the
    obvious bound $l\leq 2+5(5-2)=17$.}  By the standard bound and in
  view of the number of directions they determine, point classes of
  multiplicity $3$ or $4$ are incident with exactly two projective
  line classes.\footnote{For a point class $[x]$ of multiplicity $3$
    (``triangle'') the standard bound is $k\leq 3+3\cdot 12+2\cdot
    16=71$, and $k=69$ implies that $[x]$ is incident with two line
    classes of multiplicity $\geq 17$.}
  We will now show that there are exactly six projective
  line classes forming a dual hyperoval in $\overline{\Pi}=\PG(2,4)$,
  leaving only a small number of possibilities for the multiset
  $\overline{\mset{K}}$. Subsequently we will rule these out, thus completing
  the proof of Theorem~\ref{thmu:q=4:n=5}. 
  
  Let us denote the number of projective
  line classes by $t$.
  Since $13\cdot 4+8\cdot 2<69$, we must have at least $14$ neighbor classes 
  on points of
  multiplicity $3$ or $4$; in other words, the set of projective
  line classes must have at least $14$ double points. This forces
  $t\geq 6$, and $t=6$ is possible only if the projective
  line classes form a hyperoval in the dual plane.
%  \footnote{This is clear from the
%    3rd standard equation $a_2+3a_3+6a_4+\dots=\binom{t}{2}$ for the
%    corresponding spectrum.}
  On the other hand, consider the total number $D$ of double
  points determined by $\mset{K}$ in the various planes
  $\Pi_{[L]}$.
%  \footnote{Equivalently, the number of line segments of
%    multiplicity $2$. By the standard bound there are no line segments
%    of multiplicity $\geq 3$.} 
We have $D=f_2+3f_3+6f_4\geq 81$ with
  equality iff $f_0=f_1=f_2=0$, $f_3=13$,
  $f_4=6$.\footnote{$f_0+f_1+f_2+f_3+f_4=21$, $f_1+2f_2+3f_3+4f_4=69$
    implies $f_2+2f_3+3f_4\geq 48$, $f_3+2f_4\geq 27$, $f_4\geq
    6$. Hence $D=(f_2+2f_3+3f_4)+(f_3+2f_4)+f_4\geq 81$.}
  Since for a non-projective line class the corresponding
  plane $\Pi_{[L]}$ contains at most $6$ double points,
  this forces $21-t\geq 14$, i.e.\ $t\leq 7$. But seven lines in
  $\PG(2,4)$ have no more than $12<14$ double points, leaving $t=6$ as the
  only possibility.
%  \footnote{This follows again from the standard
%    equations for a set of $7$ points in $\PG(2,4)$. Equality holds
%    for a dual $(7,3)$-arc containing a dual hyperoval and for sets
%    $7$ lines determining a $3\times 4$-grid.}
  
  Since a dual hyperoval in $\PG(2,4)$ has $15$ double points and six $0$-points
  forming an ordinary hyperoval, we know now that
  there are either $14$ or $15$ point classes of multiplicities $3$ or
  $4$ and that the point classes of multiplicity $\leq 2$ 
  form, with one possible exception, a hyperoval in
  $\overline{\Pi}$. In particular this leaves the following
  possible multiplicities for $\overline{\mset{K}}$ and $D$.
  \begin{equation}
    \label{eq:thmu:q=4:n=5-p1}
    \begin{array}{c|cccccc}
      &f_0&f_1&f_2&f_3&f_4&D\\\hline
      \text{(a)}&0&0&6&3&12&93\\
      \text{(b)}&0&1&5&2&13&89\\
      \text{(c)}&0&0&7&1&13&88\\
      \text{(d)}&1&0&5&1&14&92\\
      \text{(e)}&0&2&4&1&14&91\\
      \text{(f)}&0&1&6&0&14&90\\
      \text{(g)}&1&1&4&0&15&94\\
      \text{(h)}&0&3&3&0&15&93
    \end{array}
  \end{equation}
  Cases (a), (d), (e), (g), (h) are ruled out immediately using the
  bound $D\leq
  6(21-t)=90$ (which expresses the fact that each of the planes $\Pi{[L]}$
  has at most six double points). In each of the remaining cases there must exist a line
  class $[L]$ of type $(4,4,4,2,2)$ (a secant class to the hyperoval formed
  by six point classes of multiplicity $\leq 2$). This is
  impossible---the $(17,5)$-arc induced by $\mset{K}$ in $\Pi_{[L]}$ would
  have six $2$-points (which must form a hyperoval) and five
  non-collinear $1$-points (which necessarily determine a secant
  of the hyperoval, resulting in a $6$-line). The proof is now complete.
\qed

For our next theorem we will need the classification of
$(22,6)$-arcs in $\PG(2,4)$ from \cite{ball-hill-ivan-ward01}. There
are precisely $6$ isomorphism classes of such arcs, which may be
described as follows: 

The so-called 
\emph{line construction} with $D\in\{1,2,3\}$ double points and $D-1$
zero points on a fixed line of $\PG(2,4)$. (Points not
mentioned are understood to have multiplicity $1$.) The corresponding
isomorphism classes will be referred to as \underline{1L},
\underline{2L},  \underline{3L}. The \emph{quadrangle construction}
\underline{Q} with $D=4$ double points forming a 
quadrangle and the diagonal points as zero points). The
\emph{hyperoval construction} \underline{H} with $D=6$ double points
forming a hyperoval and one passant containing $5$ zero points.
Finally, the \emph{dual arc construction} \underline{DA} with $D=7$
double points forming the complement of the union of four lines in
general position and $6$ zero points forming the intersection points of
these lines. An alternative description of the dual arc construction
is the following: There are
$7$ double points on three concurrent lines, one $2$-point being 
the point of concurrency of the three lines and the six remaining
$2$-points forming a hyperoval. The  six $(22,6)$-arcs are presented on the figure below
(Red, blue and white points are points of multiplicitiy 2, 1, and 0, respectively.
Blue (resp. yellow) background means that the default multiplicity for 
the points is 1 (resp. 0).)

\begin{figure}[htp!]
	\begin{center}
		\begin{tikzpicture}[line width=1pt, scale=0.4]
	
		\draw[gray] (0,0)--(8,0)--(8,5)--(0,5)--(0,0) [fill=blue!20];
		\draw[gray] (11,0)--(19,0)--(19,5)--(11,5)--(11,0) [fill=blue!20];
		\draw[gray] (22,0)--(30,0)--(30,5)--(22,5)--(22,0) [fill=yellow!20];
		\draw[gray] (0,8)--(8,8)--(8,13)--(0,13)--(0,8) [fill=blue!20];
		\draw[gray] (11,8)--(19,8)--(19,13)--(11,13)--(11,8) [fill=blue!20];
		\draw[gray] (22,8)--(30,8)--(30,13)--(22,13)--(22,8) [fill=blue!20];		
		\draw[gray] (12,10.5)--(18,10.5);
		\draw[gray] (23,10.5)--(29,10.5);
		\draw[gray] (1,1)--(3.25,2.5)--(2.5,4)--(1,1);
		\draw[gray] (4,1)--(1.75,2.5)--(2.5,4)--(4,1);
		\draw[gray] (0.6,1)--(7.4,1);
		
		\draw[gray] (4,10.5) circle (0.3cm) [fill=red];
%		\draw[black] (4,10.5) circle (0.3cm);
		\draw[gray] (15,10.5) circle (0.3cm) [fill=red];
%	    \draw[black] (15,10.5) circle (0.3cm);
		\draw[gray] (14,10.5) circle (0.3cm) [fill=red];
%		\draw[black] (14,10.5) circle (0.3cm);
		\draw[gray] (16,10.5) circle (0.3cm) [fill=white];

		\draw[gray] (24,10.5) circle (0.3cm) [fill=red];
%		\draw[black] (24,10.5) circle (0.3cm);
		\draw[gray] (25,10.5) circle (0.3cm) [fill=red];
%		\draw[black] (25,10.5) circle (0.3cm);
		\draw[gray] (27,10.5) circle (0.3cm) [fill=white];
		\draw[gray] (26,10.5) circle (0.3cm) [fill=red];
%		\draw[black] (26,10.5) circle (0.3cm);
		\draw[gray] (28,10.5) circle (0.3cm) [fill=white];
		
		\draw[gray] (1,1) circle (0.3cm) [fill=white];	
		\draw[gray] (2.5,1) circle (0.3cm) [fill=white];	
		\draw[gray] (4,1) circle (0.3cm) [fill=white];	
		\draw[gray] (5.5,1) circle (0.3cm) [fill=blue];	
		\draw[gray] (7,1) circle (0.3cm) [fill=blue];	

		\draw[gray] (12,1) circle (0.3cm) [fill=white];	
		\draw[gray] (13.5,1) circle (0.3cm) [fill=white];	
		\draw[gray] (15,1) circle (0.3cm) [fill=white];	
		\draw[gray] (16.5,1) circle (0.3cm) [fill=white];	
		\draw[gray] (18,1) circle (0.3cm) [fill=white];
			
		\draw[gray] (23,1) circle (0.3cm) [fill=blue];	
		\draw[gray] (24.5,1) circle (0.3cm) [fill=blue];	
		\draw[gray] (26,1) circle (0.3cm) [fill=blue];	
		\draw[gray] (27.5,1) circle (0.3cm) [fill=blue];	
		\draw[gray] (29,1) circle (0.3cm) [fill=red];		
\draw[gray] (29,1.8) circle (0.3cm) [fill=blue];	
\draw[gray] (29,2.6) circle (0.3cm) [fill=blue];	
\draw[gray] (29,3.4) circle (0.3cm) [fill=blue];	
\draw[gray] (29,4.2) circle (0.3cm) [fill=blue];	

	\draw[line width=1pt, color=black] (15,3) ellipse (2.25cm and 1.2cm);
\draw[line width=1pt, color=black] (25.25,3) ellipse (2.25cm and 1.2cm);

		\draw[gray] (12.8,2.9) circle (0.3cm) [fill=red];		
		\draw[gray] (17.2,2.9) circle (0.3cm) [fill=red];		
		\draw[gray] (14,2) circle (0.3cm) [fill=red];		
		\draw[gray] (16,2) circle (0.3cm) [fill=red];		
		\draw[gray] (14,4) circle (0.3cm) [fill=red];		
		\draw[gray] (16,4) circle (0.3cm) [fill=red];		

		\draw[gray] (23,2.9) circle (0.3cm) [fill=red];		
		\draw[gray] (27.4,2.9) circle (0.3cm) [fill=red];		
		\draw[gray] (24.2,2) circle (0.3cm) [fill=red];		
		\draw[gray] (26.2,2) circle (0.3cm) [fill=red];		
        \draw[gray] (24.2,4) circle (0.3cm) [fill=red];		
        \draw[gray] (26.2,4) circle (0.3cm) [fill=red];		

		\draw[gray] (2.5,2) circle (0.3cm) [fill=red];		
		\draw[gray] (1.75,2.5) circle (0.3cm) [fill=red];		
		\draw[gray] (3.25,2.5) circle (0.3cm) [fill=red];		
		\draw[gray] (2.5,4) circle (0.3cm) [fill=red];

		\draw(4,7.2) node{\small {\color{blue} $\underline{1L}$}};	
		\draw(15.5,7.2) node{\small {\color{blue} $\underline{2L}$}};	
		\draw(26,7.2) node{\small {\color{blue} $\underline{3L}$}};	
		\draw(4,-0.8) node{\small {\color{blue} $\underline{Q}$}};	
		\draw(15.5,-0.8) node{\small {\color{blue} $\underline{H}$}};	
		\draw(26,-0.8) node{\small {\color{blue} $\underline{DA}$}};		
			
		\end{tikzpicture}
	\end{center}
	\caption{The $(22,6)$-arcs in $\PG(2,4)$}
\end{figure}
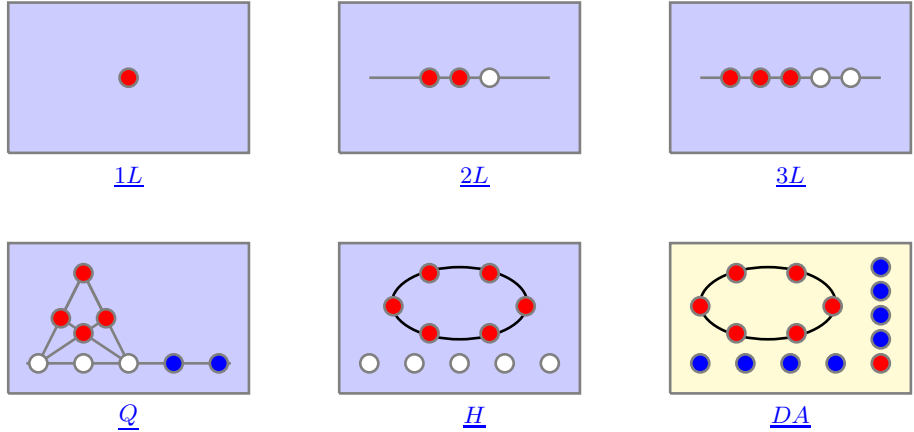

Table~1 lists the number of lines in $\PG(2,4)$ of each of 
the possible line types in the
six constructions described above.

\begin{table}[ht]
%  \label{tbl:(22,6)}
  \centering
  \begin{tabular}{c|ccccccc}
    &22200&22110&21111&11111&11110&11000&00000\\\hline
    \underline{1L}&&&5&16&&&\\
    \underline{2L}&&1&8&8&4&&\\
    \underline{3L}&1&&12&&8&&\\
    \underline{Q}&&6&8&&6&1&\\
    \underline{H}&&15&&&5&&1\\
    \underline{DA}&3&12&2&&&4&\\
  \end{tabular}
  \caption{Line types of the $(22,6)$-arcs in $\PG(2,4)$}
\end{table}

We will also use the following lemma about $(k,6)$-arcs in $\PG(2,4)$.
\begin{lemma}
  \label{lma:(k,6)_4}
  A $(k,6)$-arc in $\PG(2,4)$ with $D=8$ double points has $k\leq 19$.
\end{lemma}

\noindent {\it Proof.}
  Let $S$ be the set of $2$-points of the $(k,6)$-arc, which is
  obviously an $(8,3)$-arc. The possible spectra for such an arc are 
  $(a_0,a_1,a_2,a_3)=(3,2,10,6)$, $(2,5,7,7)$, and
  $(1,8,4,8)$, which follows easily from the standard equations.
  so that in particular the number of $3$-lines of $S$ satisfies
  $6\leq a_3\leq 8$. Let $b_i$ (resp.\ $b_i'$) be the number of points
  in $S$ (resp.\ outside $S$) which are on exactly $i$ three-lines of
  $S$. Every point in $S$ is on either two or three $3$-lines of $S$,
  so that $b_i=0$ for $i\neq 2,3$. Similarly $b_i'=0$ for $i\neq 0,1,2$. We
  have $b_2+b_3=8$, $2b_2+3b_3=3a_3$, $b_0'+b_1'+b_2'=13$,
  $b_1'+2b_2'=2a_3$, and $b_2+b_2'+3b_3=\binom{8}{2}=28$. Solving for
  $b_i$, $b_i'$ we obtain
  \begin{equation*}
    \begin{array}{c|ccccc}
      a_3&b_2&b_3&b_0'&b_1'&b_2'\\\hline
      6&6&2&4&6&3\\
      7&3&5&2&8&3\\
      8&0&8&1&8&4
    \end{array}
  \end{equation*}
  Since $b_0'\leq 4$, the arc has at most four $1$-points and so $k\leq 20$.
  Now we show that $k=20$ is impossible. A $(20,6)$-arc with $D=8$
  must have $a_3=6$, $b_3=2$. Let $p_1$, $p_2$ be the points in $S$ which are
  on three $3$-lines of $S$. We claim that $L=p_1p_2$ is a
  $2$-line of $S$. Indeed, if $L$ were a $3$-line and $p_3\in L\cap S$
  then every $3$-line would intersect $\{p_1,p_2,p_3\}$, so that the
  remaining $5$ points of $S$ would form an oval which has only
  one secant through $p_3$; a contradiction. Hence the two sets of
  $3$-lines through $p_1$, $p_2$ are disjoint and hence cover all points of
  $\PG(2,4)$ except for three points on $p_1p_2$ and one further point
  outside $p_1p_2$. Since at least one point on $p_1p_2$ must be a
  $0$-point, we have $k\leq 19$ also in this case.
  
  % \textsl{Case~1}: $S$ contains a hyperoval. Then the remaining two
  % points in $S$ determine a passant $L$ of the hyperoval, and the set of
  % points not covered by the $3$-lines of $S$ consists of the remaining
  % three points on $L$ and one additional point outside $S$. Adding
  % all four points to $S$ turn $L$ into a $7$-line, so $k\leq 19$ in
  % this case.

  % \textsl{Case~2}: $S$ contains an oval but no hyperoval. Let $O$ be
  % such an oval and $p$ its nucleus. Then $p\notin S$ and the three remaining
  % points of $S$ must be collinear. Otherwise there exist three
  % lines $L_1,L_2,L_3$ (passants of the hyperoval $O\cup\{p\}$)
  % disjoint from $S$ and forming a triangle, and hence
  % $S\cup\{p\}$ decomposes into three $3$-sets, any two of which form a
  % hyperoval; a contradiction.

  % \textsl{Case~3}: $S$ contains a quadrangle but no oval. Let $p\in S$
  % be on three $3$-lines of $S$, say $L_1$, $L_2$, $L_3$. Considering
  % the quadrangle formed by the four points of $S$ in $(L_1\cup
  % L_2)\setminus\{p\}$, we have from the hypothesis that the $6$-set $(L_1\cup
  % L_2\cup L_3)\setminus\{p\}$ is contained in two lines $L_4$, $L_5$
  % meeting outside $S$.
\qed

\begin{remark}
  \label{rmk:(k,6)_4}
Together with the well-known fact that the maximum size of a
$(9,3)$-arc in $\PG(2,4)$ is $9$, the lemma and its proof yield more
generally that the maximum size of a $(k,6)$-arc in $\PG(2,4)$ having
$D\geq 8$ double points is $k=19$, and that a $(19,6)$-arc of this
type necessarily
has six $2$-points forming a hyperoval, a passant $L$ to this hyperoval of
type $22110$ and one further $1$-point determing a triangle of three
passants with the two $1$-points on $L$.
\end{remark}

Now we are ready to state and prove the theorem.

\begin{theorem}
  \label{thmu:q=4:n=6}  
  Let $R$ be one of the rings $\G_4$, $\DN_4$, $\TDN_4$.
  \begin{enumerate}[(i)]
  \item\label{thmu:q=4:n=6:R} $\maxsize_6(R)\leq 84$;
  \item\label{thmu:q=4:n=6:ST} $\maxsize_6(R)<84$ if $R\neq\G_4$.
  \end{enumerate}
\end{theorem}

Theorem~\ref{thmu:q=4:n=5}\,(\ref{thmu:q=4:n=6:R}) was stated without proof in
\cite{it:kent}.

\noindent {\it Proof.}
  Let $\mset{K}$ be a $(k,6)$-arc in $\PHG(2,R)$.
  The standard bound (Facts~\ref{thm:bound}, \ref{cor:bound})
  gives $k\leq 82$ except in those cases where $u\in\{4,5,6\}$ and the
  points of $\mset{K}$ in every point 
  class of multiplicity $5$ and $6$ form a hyperoval resp.\ oval.
  Our first goal is to show that $u\in\{5,6\}$ implies $k<84$. Since
  for $u=4$ obviously $k\leq 21\cdot 4=84$, this already proves
  (\ref{thmu:q=4:n=6:R}), and moreover that equality in
  (\ref{thmu:q=4:n=6:R}) is possible only if
  $\mset{K}\bigl([x]\bigr)=4$ for every point class $[x]$. For this we
  will employ the classification of $(22,6)$-arcs in $\PG(2,4)$,
  Lemma~\ref{lma:(k,6)_4}, and
  we will need the fact that for hyperovals and ovals in $\AG(2,4)$
  all parallel classes have the same type, which is $2220$
  respectively $2210$.

  % Using a refined argument we will
  % first show that every $(84,6)$-arc has precisely $4$ points (of
  % multiplicity $1$) in each
  % neighbour class $[x]$ and meets every line in either $2$ or $6$
  % points. Case (\ref{thmu:q=4:n=6}) will then be completed by ruling
  % out such an arc for the rings $\DN_4$ and $\TDN_4$, and by
  % showing that it is not extendable to a $(85,6)$-arc in the case of
  % $\G_4$.
  
  \textsl{Case~1}: $u=6$. Here the standard bound gives only $k\leq 6+5\cdot
  16=86$. Assuming $k\geq 84$,
  there exist at least three line classes of multiplicity $22$
  incident with $[x]$. For such a line class $[L]$ the arc $\mset{K}$
  determines a $(22,6)$-arc in $\Pi_{[L]}$ satisfying
  $\mset{K}(p_\infty)=0$. A check of Table~1 (together
  with the fact that $0$-points are actually incident with lines of
  all types permitted by the table) reveals that this arc must be of
  type \underline{3L} and $[L]$ therefore of type
  \begin{equation*}
    \label{eq:q=4:n=6-p1}
    \begin{array}{ccccc}
      2&1&1&1&1\\
      2&1&1&1&1\\
      2&1&1&1&1\\
      0&1&1&1&1
    \end{array}.
  \end{equation*}
  (All other types of $(22,6)$-arcs force a point class $[y]$ of type
  $2211$ on $[L]$. Then the points of $\mset{K}$ in $[y]$ cannot form a
  hyperoval, we can apply the standard bound to $[y]$ to obtain $k\leq
  82$, a contradiction.) It follows that at least three line classes
  $[L]$ incident with $[x]$ contain, apart from $[x]$, only point
  classes of multiplicity $4$ which do not determine $[L]$ as a
  direction.

  Assume that there exists a second point class $[y]$ of multiplicity
  $6$. Then by applying the preceding argument to $[y]$ as well, we
  obtain at least $15$ point classes of multiplicity $4$ outside
  $[xy]$. The points of $\mset{K}$ in each such class $[z]$, as they
  determine neither $[zx]$ nor $[zy]$ as a direction, must form either
  a line segment or an affine subplane. In the first the remaining
  point class outside $[xy]$ must be empty (otherwise the requirement
  $\mset{K}\bigl([L]\bigr)\leq 12$ for the direction $[L]$ of the line
  segment could not be satisfied) and we obtain $k\leq 22+15\cdot
  4=82$, a contradiction. Hence the points of $\mset{K}$ in $[z]$ form
  an affine subplane and determine each of the remaining three
  directions through $[z]$ twice. Now let $[w]$ be a point class on
  $[xy]$ distinct from $[x]$ and $[y]$. There exists a line class
  $[L]$ through $[w]$, all of whose remaining point classes have
  multiplicity $4$ (in fact at least $3$ such line classes), and hence
  of type
  \begin{equation*}
    \label{eq:q=4:n=6-p2}
    \begin{array}{ccccc}
      *&2&2&2&2\\
      *&2&2&2&2\\
      *&0&0&0&0\\
      *&0&0&0&0
    \end{array},
  \end{equation*}
  where the '$*$'s indicate the (yet to be determined)
  multiplicities of the line segments in $[w]$. It is easily seen that
  at most one of these multiplicities can be nonzero. Hence
  $\mset{K}\bigl([w]\bigr)\leq 2$ (a point of multipliciity
    $3$ would imply $\mset{K}\bigl([L]\setminus[w]\bigr)\leq 12$,
    which is not the case).
  and $k<84$ unless there exists a third point class $[z]$ of multiplicity
  $6$. The latter, however, leads to further point classes of
  multiplicity $\leq 2$, a contradiction.

  Next assume that there exists a point class $[y]$ of multiplicity
  $5$. Applying the standard bound to $[y]$ we obtain $k\leq 5+5\cdot
  16=85$. Under the present assumption $k\geq 84$ there exist at least
  four line classes of multiplicity $21$ incident with $[y]$, of which
  at least three do not contain the point class $[x]$ of multiplicity
  $6$. For these three line classes $[L]$ the arc $\mset{K}$
  determines a $(22,6)$-arc in $\Pi_{[L]}$ satisfying
  $\mset{K}(p_\infty)=1$. Another check of Table~1
  reveals that this arc must be of type \underline{2L} with $p_\infty$
  on the line of type $22110$ or of type \underline{H}, and $[L]$
  therefore of type
  \begin{equation}
    \label{eq:q=4:n=6-p3}
    \begin{array}{ccccc}
      2&1&1&1&1\\
      2&1&1&1&1\\
      1&1&1&1&1\\
      0&1&1&1&1
    \end{array}\quad\text{or}\quad 
    \begin{array}{ccccc}
      2&2&2&1&1\\
      2&2&2&1&1\\
      1&1&1&1&1\\
      0&0&0&0&0
    \end{array}.
  \end{equation}
  All other types of $(22,6)$-arcs force a point class $[y]$ of type
  $2111$ on $[L]$. The second configuration is not possible, since
  $[L]$ intersects the line classes of multiplicity $22$ through $[x]$
  and therefore must be incident with point classes of multiplicity $4$.
  Hence we have at least $15$ point classes of multiplicity $4$
  outside $[xy]$, and proceeding as before leads to a contradiction.
  
  It remains to consider the case where all point classes outside
  $[x]$ have multiplicity $\leq 4$ and, in view of $k\geq 84$, at
  least $18$ of them have multiplicity exactly $4$. Clearly the points of
  $\mset{K}$ in a class of multiplicity $4$ must form a quadrangle
  (determing either $3$ or all $5$ directions). Point classes
    of multiplicity $4$ containing $3$ but not $4$ collinear points
    are excluded by the standard bound and
    Lemma~\ref{lma:upper_bound}, which give $k\leq 4+1\cdot 12+3\cdot
    16+1\cdot 18=82$.
  Hence the point classes of multiplicity $4$ contain a total
  of at least $18\cdot 6=108$ line segments of multiplicity $2$ 

  Consider now a line class $[L]$ not incident with $[x]$. The arc
  $\mset{K}$ induces in $\Pi_{[L]}$ a $(22,6)$-arc with
  $\mset{K}(p_\infty)=2$. Since such an arc has at most $7$ double
  points (see Table~1), at most $6$ line segments have
  $[L]$ as their direction.

  This leaves at least $108-6\cdot 16=12$ line segments of
  multiplicity $2$ in point classes of multiplicity $4$ whose
  direction passes through $[x]$. However, as we have seen, at most two line
  classes through $[x]$ can be directions of line segments of
  multiplicity $2$. Hence there must be exactly two, 
  each being the direction of $9$ line segments ($3$ in $[x]$ and $6$
  outside $[x]$). If we consider now for these line classes $[L]$ the
  arc induced by $\mset{K}$ in $\Pi_{[L]}$, we obtain the final
  contradiction: Since $(9,3)$-arcs in $\PG(2,4)$ are complete, we
  have $\mset{K}\bigl([L]\bigr)=18$ and hence $k\leq 6+19\cdot
  4=80$.
%  \footnote{In other words, apart from $[x]$ and $18$ point
%    classes of multiplicity $4$ there must be two empty classes.}

  \textsl{Case~2}: $u=5$. Assuming $k\geq 84$, the class $[x]$ 
  must be incident with at least four line
  classes $[L]$ of multiplicity $21$ and hence of one of the two types
  in \eqref{eq:q=4:n=6-p3}.

  Assume that there exists a second class $[y]$ of multiplicity
  $5$. Then among the line classes incident with $[x]$ or $[y]$ but
  different from $[xy]$ only
  one of the types in \eqref{eq:q=4:n=6-p3} can occur. So either we
  have at least three line classes of the first type incident $[x]$
  (resp. with $[y]$) and can then proceed as in Case~1; or both $[x]$
  and $[y]$ are incident with at least four line classes 
  of the second type, one of which may be the class $[xy]$. Continuing
  with the latter, we claim that there are at most $9$ line classes
  of multiplicity $5$. Indeed, suppose there exist more than $9$ such
  line classes. Then, since $\maxsize_3(4)=9$, there exists a line
  class of type $5555\ast$. But this means that $\mset{K}$
  determines in $\Pi_{[L]}$ a $(k,6)$-arc with $\mset{K}(p_\infty)=1$,
  $D\geq 8$, and $k\geq 21$, contradicting Lemma~\ref{lma:(k,6)_4}.
  
  This implies that both $[x]$ and $[y]$ are incident with exactly four
  line classes of type $55533$, one of which is $[xy]$. These
  seven line classes cover all but one point class, so that
  $k\leq 9\cdot 5+11\cdot 3+4=82$, a contradiction.

  The remaining subcase where all point classes except $[x]$ have
  multiplicity $\leq 4$ is reduced ad absurdum in essentially the same
  way as in Case~1.

  This completes the proof of Part~(\ref{thmu:q=4:n=6:R}) of the
  theorem.

  From now on we assume $k=84$, so that by the first part of the proof
  every point class contains exactly $4$ points of $\mset{K}$.
  
  Our next goal is to show that $\mset{K}$ meets every line of
  $\PHG(2,R)$ in either $2$ or $6$ points. For every line class $[L]$
  the arc $\mset{K}$ now determines a $(22,6)$-arc in $\Pi_{[L]}$ with
  $\mset{K}(p_\infty)=2$ and hence $D\leq 6$ double-points distinct
  from $p_\infty$; cf.\
  Table~{tbl:(22,6)}. On the other hand every point class $[x]$
  accounts for six $2$-points (line segments of multiplicity $2$),
  since the points of $\mset{K}$ in
  $[x]$ now necessarily form a quadrangle (of one of two possible
  types, determing either $3$ or $5$ directions). A double count
  of all line segments of multiplicity $2$ reveals that we must have $D=7$
  for each induced $(22,6)$-arc and hence that all these arcs are of
  type \underline{DA}. This in turn gives the following two
  possibilities (depending on whether $p_\infty$ is the point of
  concurrency of the three lines of type $22200$ or not)
  for the type of a line class $[L]$:
  \begin{equation*}
    \label{eq:q=4:n=6-p4}
    \begin{array}{ccccc}
      2&2&2&1&1\\
      2&2&2&1&1\\
      0&0&0&1&1\\
      0&0&0&1&1
    \end{array}\quad\text{or}\quad 
    \begin{array}{ccccc}
      2&2&2&2&2\\
      2&1&1&1&1\\
      0&1&1&1&1\\
      0&0&0&0&0
    \end{array}.
  \end{equation*}
  Regardless of the types which actually occur each line of
  $\PHG(2,R)$ of multiplicity $<6$ must have multiplicity $2$ (since
  the $(22,6)$-arcs of type \underline{DA} have this property, cf.\
  \cite{ball-hill-ivan-ward01}).\footnote{Different types of point and
    line classes are conceivable, for example $3$ point classes of the
    rectangle type (i.e.\ determining only $3$
    directions) forming a triangle in $\PG(2,4)$ and $18$ point
    classes of the skew quadrangle type (i.e.\ determining $5$
    directions). The triangle sides must then have the
    first type in \eqref{eq:q=4:n=6-p4} and coincide with the $2200$
    directions of its three point classes of the skew quadrangle
    type. All other line classes must have the second type in
    \eqref{eq:q=4:n=6-p4}, and the
    $2200$ directions of the remaining $9$ point classes must match
    the $9$ passant directions of the triangle.}

  Now we are ready to prove Part~(\ref{thmu:q=4:n=6:ST}). The argument
  will be essentialy the same as in the proof of
  \cite[Th.~4(ii)]{imt:two_exact_values} and will consist of showing
  that the existence of a $(84,6)$-arc in $\PHG(2,R)$ is incompatible
  with the existence of a subplane isomorphic to $\PG(2,4)$ (a
  set of $21$ points forming a transversal for the point classes and
  meeting every line in either $1$ or $5$ points). Such subplanes
  exist in $\PHG(2,\DN_4)$ and $\PHG(2,\TDN_4)$, but not in
  $\PHG(2,\G_4)$. An example is the set of all points of $\PHG(2,R)$
  with coordinates in a fixed subfield $F\cong\F_4$ of $R$.

  Let $B$ be a Baer subplane of $\PHG(2,R)$,
  $\abs{\mset{K}\cap B}=t$ and $b_2$ (resp.\ $b_6$) the number of
  $2$-lines ($6$-lines) of $\mset{K}$ meeting $B$ in $5$ points. We
  have the equations
  \begin{equation*}
    \label{eq:q=4:n=6-p5}
    \begin{array}{rcrcl}
      b_2&+&b_6&=&21,\\
      2b_2&+&6b_6&=&5t+1\cdot(84-t)=84+4t.
    \end{array}
  \end{equation*}
But this system has no solution modulo 4, and hence no solution in
$\Z$, a contradiction.
\qed

Similar to the case $n=6$, the standard bound $\maxsize_7(R)\leq 106$
can be sharpened by emplyoing the classification of $(27,7)$-arcs in
$\PG(2,4)$ (also due to \cite{ball-hill-ivan-ward01}). There exist
precisely two isomorphism classes of such arcs, having $D=6$ resp.\
$D=9$ double points: A \emph{hyperoval construction} \underline{H}
having six $2$-points forming a hyperoval and fifteen $1$-points, and
the \emph{dual arc} (or \emph{triangle}) \emph{construction}
\underline{DA} having three $0$ points forming a triangle, nine
$1$-points on the sides of the triangle and nine
$2$-points.\footnote{Using the language of minihypers, \underline{H}
  corresponds to the unique indecomposable $(15,3)$-minihyper (the
  complement of a hyperoval), and \underline{DA} corresponds to the
  $(15,3)$-minihyper which is a sum of three distinct, non-concurrent
  lines.}  Table~\ref{tbl:(27,7)} gives the corresponding line types
and their frequencies. Below we will use the fact that
points are actually incident with
lines of all types permitted by the table.
%\footnote{This can be seen
%  as a consequence of the fact that in both cases the automorphism
%  group of the arc acts transitively on the $0$-, $1$-, and $2$-points
%  of the arc.}

\begin{table}[ht]
  \label{tbl:(27,7)}
  \centering
  \begin{tabular}{c|cccc}
    &22210&22111&11111&11100\\\hline
    \underline{H}&&15&6&\\
    \underline{DA}&9&9&&3\\
  \end{tabular}
  \caption{Line types of the $(27,7)$-arcs in $\PG(2,4)$}
\end{table}
\begin{theorem}
  \label{thmu:q=4:n=7}
  Let $R$ be one of the rings $\G_4$, $\DN_4$, $\TDN_4$. Then
  $\maxsize_7(R)\leq 101$.
\end{theorem}
\noindent {\it Proof.}
  Let $\mset{K}$ be a $(k,7)$-arc in $\PHG(2,R)$.
  The standard bound (Theorem~\ref{thm:bound}) yields $k\leq 101$
  except in the following three cases, which are then ruled out with the help of
  the classification of the $(27,7)$-arcs in $\PG(2,4)$.

  \textsl{Case~1}: $u=6$, and the points of $\mset{K}$ in $[x]$ form a
  hyperoval. Here the standard bound gives only $k\leq 6+5\cdot
  20=106$. Assuming $k\geq 102$, at least one line class $[L]$
  incident with $[x]$ has multiplicity $26$. The arc $\mset{K}$
  determines a $(27,7)$-arc in $\Pi_{[L]}$ with
  $\mset{K}(p_\infty)=1$. An inspection of Table~2
  % Table~\ref{tbl:(22,6)}
  shows that this arc must be of type \underline{DA} and that $[L]$
  must contain a point class $[y]$ of type $2211$. The points of
  $\mset{K}$ in $[y]$ can neither form a hyperoval nor two parallel
  line segments, so that applying the standard bound to $[y]$ gives
  $k\leq 6+2\cdot 16+3\cdot 20=98$, a contradiction.

  \textsl{Case~2}: $u=6$, and the points of $\mset{K}$ in $[x]$ form two
  parallel line segments. Here the standard bound gives $k\leq
  6+1\cdot 16+4\cdot 20=102$. If this inequality is sharp, there exist
  line classes of multiplicity $26$ incident with $[x]$ and we can
  argue as in Case~1.

  \textsl{Case~3}: $u=5$, and the points of $\mset{K}$ in $[x]$ form
  an oval. Here the standard bound is $k\leq 105$. Assuming $k\geq
  102$ (in fact $k\geq 101$ suffices here), there exists a line class
  $[L]$ of multiplicity $25$ incident with $[x]$, and $\mset{K}$
  determines a $(27,7)$-arc in $\Pi_{[L]}$ with $\mset{K}(p_\infty)=2$. 
  Table~1 shows that this arc is of type
  \underline{H} and that $[L]$ contains a point class $[y]$ of type
  $2111$. As the points of $\mset{K}$ in such a class cannot form an
  oval, we obtain $k\leq 101$ by applying the standard bound to $[y]$.
\qed

%\begin{longver}

\begin{theorem}
  \label{thmu:q=4}
  Let $R$ be one of the rings $\DN_4$, $\TDN_4$. Then
  $\maxsize_8(R)\leq 125$.
\end{theorem}
For the Galois ring $\G_4$ we have
$\maxsize_8(\G_4)=126$; cf.\ Th.~\ref{thm:dualarcs}.

\noindent {\it Proof.}
  Suppose $\mset{K}$ is a $(126,8)$-arc in
  $\PHG(2,R)$, where $R=\DN_4$ or $\TDN_4$. We will show that the
  passants ($0$-lines) of $\mset{K}$ form a hyperoval (i.e., an
  $(21,2)$-arc) in the dual 
  plane $\PHG(2,R)^\ast$, which is isomorphic to $\PHG(2,R)$. The
  result then follows from the known non-existence of hyperovals in
  the planes $\PHG(2,\DN_4)$ and $\PHG(2,\TDN_4)$; cf.
  \cite[Th.~13]{it:ovals}.

  Counting points on all lines through a fixed point $x$ gives
  \begin{equation}
    \label{eq:20mal8}
    16\mset{K}(x)+3\mset{K}\bigl([x]\bigr)+126=\sum_{\substack{L\in\lines\\
        x\in L}}\mset{K}(L)\leq 20\cdot 8=160.
  \end{equation}
  This implies $\mset{K}(x)\in\{0,1\}$ and
  $\mset{K}\bigl([x]\bigr)\leq 6$. On the other hand
  $\mset{K}(\points)=126=21\cdot 6$, so necessarily
  $\mset{K}\bigl([x]\bigr)=6$ for every class $[x]$. Further, since
  \eqref{eq:20mal8} is sharp for $1$-points of $\mset{K}$, we must have
  $\mset{K}(L)\in\{0,8\}$ for every line $L$. 

  Consequently the number of $8$-lines ($0$-lines) through a $0$-point of
  $\mset{K}$ is ($3\cdot 6+126)/8=18$ (resp.\ $20-18=2$), and the total
  number of $8$-lines ($0$-lines) is $18\cdot 10\cdot 21/12=315$
  (resp.\ $336-315=21$). In particular the $0$-lines form a dual
  $(21,2)$-arc as claimed.

  % The standard bound gives $\mset{K}\bigl([x]\bigr)\leq 6$ for every
  % neighbour class $[x]$, and that the points in each class $[x]$ with
  % $\mset{K}\bigl([x]\bigr)=6$ form a hyperoval.
\qed

\begin{remark}
In the cases $q=2,3$ the exact value of $m_3(R)$ is known: For
$q=2$ the values are
$\maxsize_3(\Z_4)=\maxsize_3(\DN_2)=10$. For $q=3$ we have
$\maxsize_3(\Z_9)=19$, $\maxsize_3(\DN_3)=18$. 
\end{remark}

%\end{longver}

%\begin{theorem}
%\label{thm:n=4}
%\[m_4(R_R^3)\le\left\{
%\begin{array}{ll}
%5+8q + (q-3)(m_4(q)-4) & \text{ for } q \text{ even},\\
%5+10q + (q-4)(m_4(q)-4) & \text{ for } q \text{ odd}.
%\end{array}\right.\]
%\end{theorem}

%For $q=5$, this theorem implies $m_4(R_R^3)\le69$.
%The next upper bounds follow easily from Lemma~\ref{lma:upper_bound}.

Finally we treat the case $q^2=25$. Here our first goal is to derive a
tight upper bound for $\maxsize_3(R)$. In the case where
$u=\max_{[x]}\mset{K}\bigl([x]\bigr)\leq 3$ we will use
Lemma~\ref{lma:upper_bound} and an adaption of the techniques used in
the proofs for $\abs{R}<25$ (cf.\ \cite[Th.~5.1,5.2]{it:kent}).
For the case  $u=4$, which is not covered by Lemma~\ref{lma:upper_bound},
we need an additional Lemma.

\begin{lemma}
  \label{lmau:q=5:n=3}
  If a $(k,4)$-arc in $\PG(2,5)$ has the property that all its
  $4$-lines are concurrent in a point, then $k\leq
  12$.
\end{lemma}

\noindent {\it Proof.}
  Let $\mset{F}$ be such an arc. If there is at most one $4$-line or
  the point of concurrency of the $4$-lines (henceforth denoted by
  $p$) is not a $0$-point, the
  assertion follows form the fact that the maximum size of a $3$-arc
  in $\PG(2,5)$ is $11$. Likewise, if $\mset{F}$ has a multiple point
  the assertion is clear.
%  \footnote{A multiple point different from $p$
%    is on at most one $4$-line and so $k\leq 2+1\cdot 2+5\cdot 1=9$.}
  So we may assume that $\mset{F}$ has two
  $4$-lines $L_1$, $L_2$ concurrent in a $0$-point $p$, and $\mset{k}$
  is projective. We denote the second $0$-point on $L_i$ by $p_i$
  ($i=1,2$).

  Suppose $\mset{k}$ contains at least $4$ points $q_1$, $q_2$, $q_3$, $q_4$
  outside $L_1\cup L_2$. If three of them, say $q_1$, $q_2$, $q_3$,
  are on a further line $L_3$ through $p$, we must also have $q_4\in L_3$, 
  since otherwise one of the
  three lines $q_iq_4$, $i=1,2,3$, would be a $4$-line. This shows
  $\mset{F}=L_1\cup L_2\cup L_3$ and $k=12$. If $q_1$, $q_2$, $q_3$
  are on a line not through $p$, this line must be $p_1p_2$. Using the
  preceding argument we get a contradiction, since $q_4\in p_1p_2$ is
  impossible. Finally suppose $q_1$, $q_2$, $q_3$, $q_4$ form a
  quadrangle. Arguing in a similar way we find that $p$, $p_1$, $p_2$
  must be the diagonal points of this quadrangle. Adding a further
  point to $\mset{k}$ is not possible,
%  \footnote{There are two
%    classes of points to distinguish: points on diagonals through $p$
%    (these are on a $3$-line through a quadrangle vertex),
%    and points on the remaining diagonals (which are themselves
%    $3$-lines).}
  so that $k=12$ holds in this case as well, and we are done.
\qed

\begin{remark}
	\label{rmk:q=5:n=3}
The proof shows that a $(12,4)$-arc with
	this property either (i) is a (possibly non-projective) extension
	of a maximal $(11,3)$-arc, or (ii) consists of a quadrangle and
	$8$ points, distinct from the diagonal points of the quadrangle,
	on two of the three lines joining the diagonal points, or (iii) of
	three $4$-lines concurrent in a $0$-point.
\end{remark}

\begin{theorem}
  \label{thm:q=5:n=3}
  $\maxsize_3(R)\leq 43$ for $R=\Z_{25}$ and $R=\DN_5$.
\end{theorem}

\noindent {\it Proof.}
  Let $\mset{K}$ be a $(k,3)$-arc in $\PHG(2,R)$, and fix a point
  class $[x]$ with
  $\mset{K}\bigl([x]\bigr)=u=\max\bigl\{\mset{K}([y]);y\in\points\bigr\}$. The
  standard bound (Fact~\ref{thm:bound}) gives $k\leq 36$ for $u\geq
  6$, $k\leq 40$ for $u=5$ (noting that a $(5,3)$-arc in $\AG(2,5)$
  determines at least four directions
%  \footnote{A rectangle together with the
%    intersection point of its diagonals determines at least exactly
%    four directions.} 
and a $(5,2)$-arc at
  least five directions), and $k\leq 4+4\cdot 5+2\cdot 10=44$ for $u=4$
  (since a quadrangle in $\AG(2,5)$ determines at least four
  directions). 

  Continuing with the case $u=4$,
  let $[L]$ be a line class incident with $[x]$ which
  is not determined by two points in $\mset{K}\cap[x]$. By
  Lemma~\ref{lmau:q=5:n=3} we have $\mset{K}\bigl([L]\bigr)\leq 12$,
  i.e.\ $\mset{K}\bigl([L]\setminus [x]\bigr)\leq 8$,
  so that we can strengthen the preceding bound to
  $k\leq 4+4\cdot 5+2\cdot 8=40$.

  Next let $u=3$. Applying Lemma~\ref{lma:upper_bound} and using
  $\maxsize_3(5)=11$ we obtain $k\leq 3+1\cdot 0+5\cdot(11-3)=43$ if
  the points in $\mset{K}\cap[x]$ are collinear, respectively $k\leq
  3+3\cdot 5+3\cdot(11-3)=42$ if they form a triangle.

  It remains to consider the case $u=2$. Obviously every line class $[L]$
  is determined at most once as a direction by the points of
  $\mset{K}$ in the various point classes $[x]$. If this is the case, we have
  $\mset{K}\bigl([L]\bigr)\leq 7$. Otherwise
  $\mset{K}\bigl([L]\bigr)\leq 10$, since the multiset induced by
  $\mset{K}$ in $[L]$ can be extended to a $3$-arc with
  $\mset{K}(p_\infty)=1$ in the projective
  plane corresponding to $[L]$. Proceeding as in the
  proof of Th.~\ref{thmu:q=3:n=7} we obtain
  \begin{equation*}
    6(f_1+2f_2)=6k=\sum_{[L]}\mset{K}\bigl([L]\bigr)\leq 7f_2+10(31-f_2),
  \end{equation*}
  i.e.\ $6f_1+15f_2\leq 310$. 
  % We have the following two cases.

  % \textsl{Case~1}: $f_0=0$. In this case
  % $\mset{t}=\bigl\{[y];\mset{K}([y])=2\bigr\}$ is a tangency set in
  % $\PG(2,5)$, since all point classes different from $[y]$ on the line
  % class (direction) determined by $[y]$ must have multiplicity $1$. As
  % remarked above, this implies $f_2\leq 12$ and hence $k\leq
  % 19+2\cdot 12=43$.

  For $f_2\leq 12$ we have $k\leq
  19+2\cdot 12=43$, and for $f_2\geq 16$ we also have
  $k\leq\left\lfloor\frac{310-3f_2}{6}\right\rfloor
  \leq\left\lfloor\frac{262}{6}\right\rfloor=43$. This leaves the following
  non-negative integer solutions of $f_1+f_2\leq 31$, $6f_1+15f_2\leq
  310$, and $k=f_1+2f_2>43$:
  \begin{equation*}
    \begin{array}{c|cccc}
      &f_0&f_1&f_2&k\\\hline
      (a)&0&18&13&44\\
      (b)&1&16&14&44\\
      (c)&2&14&15&44
    \end{array}
  \end{equation*}
  Now we check that in none of the three cases there exists a
  corresponding multiset $\overline{\mset{K}}$ in $\PG(2,5)$
  satisfying the restrictions above (i.e.\ $2$-points are incident
  with a line of multiplicity at most $7$ and all lines have
  multiplicity at most $10$). The third standard equation for the
  spectrum of $\overline{\mset{K}}$ is $\sum
  i^2a_i=5(f_1+4f_2)+k^2$. We compare this with the ``ideal'' line
  distribution $(b_i)$ given by $b_7=f_2$, $b_{10}=31-f_2$ and $b_i=0$
  otherwise:
  \begin{equation*}
    \begin{array}{c||ccc|ccc}
      &\sum a_i&\sum ia_i&\sum i^2a_i&\sum b_i&\sum ib_i&\sum
      i^2b_i\\\hline
      (a)&31&264&2286&31&271&2437\\
      (b)&31&264&2296&31&268&2386\\
      (c)&31&264&2306&31&265&2335
    \end{array}
  \end{equation*}
  Intuitively, since the first moments $\sum ia_i$ and $\sum ib_i$ are
  close to each other, the same should be true of the second moments
  $\sum i^2a_i$ and $\sum i^2b_i$, which is not the case. We thus have
  a contradiction. 

  To make the argument precise, note that
  $(a_i)$ must satisfy the additional restrictions
  $a_0+a_1+\dots+a_7\geq f_2$ and $a_i=0$ for $i>10$, from which
  $a_0+a_1+\dots+a_i\geq b_0+b_1+\dots+b_i$ for $0\leq i\leq
  10$. Hence, putting $A_i=a_i+a_{i+1}+\dots+a_{10}$,
  $B_i=b_i+b_{i+1}+\dots+b_{10}$ ($0\leq i\leq 10$), we have $A_i\leq
  B_i$ for all $i$. Moreover, $\sum_{i=1}^{10}A_i=\sum ia_i$,
  $\sum_{i=1}^{10}(2i-1)A_i=\sum i^2a_i$ and similarly for
  $B_i$. We obtain
  \begin{align*}
    \sum i^2(b_i-a_i)&=\sum_{i=1}^{10}(2i-1)(B_i-A_i)\leq
    19\sum_{i=1}^{10}(B_i-A_i)\\
    &=19\sum i(b_i-a_i).
  \end{align*}
  This inequality is violated in all three cases (cf.\ the table
  above), and hence $\overline{\mset{K}}$ does not exist.
  
  This contradiction completes the proof of the theorem.
\qed

\begin{remark}
  \label{rmk:tset}
  In Case~(a) above the set
  $\mset{T}=\bigl\{[y];\mset{K}([y])=2\bigr\}$ is a \emph{tangency
    set} (a point set all of whose points are on at least one tangent)
  in $\PG(2,5)$, since all point classes different from $[y]$ on the
  line class (direction) determined by $[y]$ must have multiplicity
  $1$. This case can also be ruled out using the bound
  $f_2=\abs{\mset{T}}\leq\left\lfloor5\sqrt{5}+1\right\rfloor=12$
  cf.~\cite{illes-szonyi-wettl91,bruen-drudge99a}.
\end{remark}

\begin{theorem}
\label{thmu:q=5}
Let $R$ be a chain ring of nilpotency index 2 with
$|R|=25$. Then

\begin{enumerate}[(i)]
\item\label{thmu:q=5:n=4} $\maxsize_4(R)\le 70$;
\item\label{thmu:q=5:n=5} $\maxsize_5(R)\le 102$;
\item\label{thmu:q=5:n=6} $\maxsize_6(R)\le 130$.
\end{enumerate}
\end{theorem}

\noindent {\it Proof.}
  (\ref{thmu:q=5:n=4}) Let $\mset{K}$ be a $(k,4)$-arc in $\PHG(2,R)$,
  and fix a point class $[x]$ with
  $\mset{K}\bigl([x]\bigr)=u=\max\bigl\{\mset{K}([y]);y\in\points\bigr\}$.
  For $u\geq 5$ we use the standard bound (Fact~\ref{thm:bound}) to
  obtain $k\leq 70$. Equality can hold only if $u=5$ and the points of
  $\mset{K}$ in $[x]$ (and hence in every neighbour class $[y]$ with
  $\mset{K}\bigl([y]\bigr)=5$) form a $(5,2)$-arc.\footnote{For this
    we need the fact that $5$ points in $\AG(2,5)$ either possess two
    $3$-lines or determine at least $5$ directions. $(70,4)$-arcs
    realizing this configuration are conceivable, for example four
    $5$-classes and two $0$-classes on a line class $[L]$, and
    twenty-five $2$-classes outside $[L]$ directed towards a
    $0$-class.}  For $u\leq 4$ we use Lemma~\ref{lma:upper_bound}
  together with $\maxsize_4(5)=16$ to obtain $k<70$. For example, if
  $u=3$ and the three points of $\mset{K}$ in $[x]$ form a triangle,
  we get $k\leq 3+3\cdot
  10+3\cdot\bigl(\maxsize_4(5)-4)=69$.
%  \footnote{Apart from the above,
%    this case is in fact the only obstruction to the stronger bound
%    $k\leq 68$.}

  (\ref{thmu:q=5:n=5}) Assuming $\mset{K}$ to be a $(k,5)$-arc and
  using the same notation as before, the standard bound gives $k\leq
  100$ except when $u\in\{4,5\}$ and the points in $\mset{K}\cap[x]$
  are either collinear or form a quadrangle determining only four
  directions.\footnote{Equality can hold only if $u=5$ and the points
    of $\mset{K}$ in $[x]$ form a $(5,2)$-arc with taking
    the unique point incident with five tangents to lie on the line at infinity.}  

  First assume that $u=5$ and the points of $\mset{K}$ in every
  neighbour class $[x]$ with $\mset{K}\bigl([x]\bigr)=5$ are
  collinear. In this case the standard bound is $k\leq 5+1\cdot
  0+5\cdot 20=105$.  Assuming $k\geq 103$, there certainly exists a
  line class $[L]$ through $[x]$ of multiplicity $25$ ($20$ without
  the points of $[x]$).  Since a $(25,5)$-arc in $\PG(2,5)$ is
  equivalent to $\AG(2,5)$ (i.e. has the points of the plane with a
  line deleted), the line class $[L]$ contains five point classes of
  multiplicity $5$ and an empty class. However, already three point
  classes of multiplicity $5$ are enough to get the bound $k\leq
  3\cdot 5+16\cdot 5=95$ (obtained from the three corresponding line
  classes of multiplicity $5$, which together cover at least $15$
  point classes), a contradiction.
%  \footnote{Collinearity of the
%    $5$-classes on $[L]$ and the known size of the second smallest
%    blocking set in $\PG(2,5)$ yield in fact the much stronger bound
%    $k\geq 90$ for the special case under consideration.}
  Thus we have at least $21$ point classes of
  multiplicity $5$.  This is impossible, since the point classes of
  multiplicity $5$ clearly form a tangency set; cf.\
  Remark~\ref{rmk:tset}.\footnote{Remark~\ref{rmk:tset} and the known
  size of the second smallest blocking set in $\PG(2,5)$ yield in
  fact a much stronger bound ($k\leq 98$) for the case under
  consideration.}

  % Now let $u=4$ with $4$ collinear points in every point neighbour
  % class of multiplicity 
  % $4$. Here the standard bound is $k\leq 4+1\cdot 5+5\cdot 20=109$,
  % and line classes of multiplicity $24$ incident with $[x]$ contain
  % six point classes of multiplicity $4$. (For the latter note that
  % $\mset{K}$ now induces a $(25,5)$-arc containing $p_\infty$ in the
  % projective plane corresponding to $[L]$.) There exists at least one
  % such line class unless $k\leq 104$. Hence, assuming $k\geq 105$ we
  % have at least six point classes of multiplicity $4$. But already two
  % such point classes, as they are incident with a line class of
  % multiplicity $9$, force $k\leq 20\cdot 4+2\cdot 9=98$, a contradiction.

  Now let $u=4$. The standard bounds are $k\leq 4+1\cdot 5+5\cdot
  20=109$ if the $4$ points of $\mset{K}$ in $[x]$
  are collinear, respectively $k\leq 4+4\cdot 15+2\cdot
  20=104$ if they form a quadrangle. 
  
  Let the four points of $\mset{K}$ in $[x]$ be collinear and assume
  $k\geq 103$. If there exists a second point class containing $4$
  collinear points of $\mset{K}$, we have two line classes of
  multiplicity at most $9$ and hence the bound $k\leq 2\cdot 9+20\cdot
  4=98$, a contradiction. If the points of $\mset{K}$ in some point
  class $[y]$ form a quadrangle, we also obtain a contradiction by the
  following argument: Consider the line class $[L]$ of multiplicity
  $\leq 9$ through $[x]$. All point classes $[z]\neq[x]$ on $[L]$ have
  multiplicity $\leq 2$, since otherwise the standard bound applied to
  $[z]$ yields a contradiction. (The critical case is that of $3$
  collinear points in $\mset{K}\cap[z]$, in which $k\leq 3+1\cdot
  10+4\cdot 20+1\cdot 9=102$.) This in turn gives that all line
  classes through $[y]$ except $[xy]$ have multiplicity $\leq 22$.
  Applying the standard bound to $[y]$ we obtain $k\leq 4+4\cdot 15+
  1\cdot 20+1\cdot 18=102$, the said contradiction. In all we have
  shown that $[x]$ is the only point class of multiplicity $4$. Hence
  $k\leq 4+30\cdot 3=94$, a contradiction ruling out the case of $4$
  collinear points in $\mset{K}\cap[x]$.

  This leaves the case $u=4$ with the points of $\mset{K}$ in all
  point classes of multiplicity $4$ forming a quadrangle. Assuming
  $k\geq 103$, the standard bound yields that every such class is
  incident with at least one line class of multiplicity $24$, which
  necessarily contains six point classes of multiplicity $4$. (For the
  latter note that $\mset{K}$ now induces a $(25,5)$-arc containing
  $p_\infty$ in the projective plane corresponding to $[L]$.) Hence
  the set of point classes of multiplicity $4$ is equal to the union
  of some line classes of multiplicity $24$. These line classes must
  be in general position (since each point class of multiplicity $4$
  is incident with at most two of them), bounding their number by $6$.
  Also by the standard bound, the second largest line class through a
  point class of multiplicity $4$ must have multiplicity $\geq 23$ and
  hence contain at least five point classes of multiplicity
  $4$. However, any union $S$ of $t\leq 6$ lines in $\PG(2,5)$ in general
  position contains points which are on exactly one line of multiplicity
  $\geq 5$ (by a direct check).
%  \footnote{For $t\neq 5$ any $1$-point of
%    $S$ has this property; for $t=5$ five of the $1$-points are
%    collinear and the remaining five $1$-points have this property.}). 
Thus we have the final contradiction to the assumption $k\geq 103$. 

  (\ref{thmu:q=5:n=6}) For a $(k,6)$-arc $\mset{K}$ the standard bound
  yields $k<130$ except when $u=5$ and the points of $\mset{K}$ in
  $[x]$ form a line segment or a $(5,2)$-arc. In the latter case we
  get $k\leq 130$. In the former case we only get $k\leq 135$, but the
  proof can be finished using essentially the same argument as in
  (\ref{thmu:q=5:n=5}): Assuming that $k\geq 131$ and the points of
  $\mset{K}$ in all classes of multiplicity $5$ form line segments, we
  have that $[x]$ is incident with at least one line class $[L]$ of
  multiplicity $30$, on which $\mset{K}$ then necessarily induces the
  trivial $(30,6)$-arc (assigning multiplicity $1$ to every line
  segment in $[L]$). Thus we have at least six point classes of
  multiplicity $5$, and a similar count as before, $k\leq 20\cdot
  5+2\cdot 10=120$, gives the desired contradiction.

\qed

\section{Tables with bounds on the maximal sizes of arcs in 
Hjelmslev planes over small chain rings}
\label{sec:tables}

In the tables below, we summarize our knowledge about the values of
$\maxsize_n(R)$ for the chain rings $R$ with $|R|=q^2\le25$, $R/\rad
R\cong\mathbb{F}_q$. We give information about all values of $n$ with
$2\le n\le q^2-1$. As remarked at the beginning of
Section~\ref{ssec:genconst}, the cases $n=q^2,q^2+1,\ldots,q^2+q$ are
covered by \cite[Cor.~2]{it:deadfin}.
%The new values reported in this paper are given in boldface.
We want to point out the fact that there are 
lots of examples with
$\maxsize_n(R)\ne\maxsize_n(S)$ for nonisomorphic chain rings $R$, $S$ with
$|R|=|S|$, $R/\rad R\cong S/\rad S$.
% (cf. Theorems~\ref{thm:teichcap}, \ref{thm:ovalimp}
% and the results in Section~\ref{ssec:dualconst}). 
However, in all these examples $\Char R\ne\Char S$, and we do not have
a single example of chain rings $R$ and $S$ of the same order, length
and characteristic, in which the values of $\maxsize_n(R)$ and
$\maxsize_n(S)$ are different.

\begin{table}[htbp]
  \centering
  $\begin{array}{|c||r|r|}
    \hline
    n/R&\Z_9&\DN_3\\\hline\hline
    0&0&0\\\hline
    1&1&1\\\hline
    2&7^a&6^a\\\hline
    3&10^a&10^a\\\hline
    4&16^a&16^a\\\hline
    5&22^a&22^a\\\hline
    6&28&28\\\hline
  \end{array}$
  \caption{The numbers $\maxsize_n(R)$ for chain rings of order $4$}
  \label{tbl:q=2}
\end{table}

\begin{table}[htbp]
  \centering
  $\begin{array}{|c||c@{\;-\;}c|c@{\;-\;}c|}
    \hline
    n/R& \multicolumn{2}{|c|}{\Z_9}&\multicolumn{2}{|c|}{\DN_3}\\\hline\hline
    0 & \multicolumn{2}{|c|}{0} & \multicolumn{2}{|c|}{0} \\\hline
    1 & \multicolumn{2}{|c|}{1} & \multicolumn{2}{|c|}{1} \\\hline
    2 & \multicolumn{2}{|c|}{9^a} & \multicolumn{2}{|c|}{9^a} \\\hline
    3 & \multicolumn{2}{|c|}{{}^c19^a} & \multicolumn{2}{|c|}{18^c} \\\hline
    4 & \multicolumn{2}{|c|}{{}^f30^b} & \multicolumn{2}{|c|}{{}^D30^b} \\\hline 
    5 & \multicolumn{2}{|c|}{{}^f39^h} & \multicolumn{2}{|c|}{38^h} \\\hline 
    6 & \multicolumn{2}{|c|}{{}^D49^i} & \multicolumn{2}{|c|}{{}^D50^i} \\\hline
    7 & \multicolumn{2}{|c|}{^{fB}60^i} & \multicolumn{2}{|c|}{{}^B60^i} \\\hline
    8 & \multicolumn{2}{|c|}{{}^{fh}69^h} & \multicolumn{2}{|c|}{69^h} \\\hline
    9 & \multicolumn{2}{|c|}{81^a} & \multicolumn{2}{|c|}{81^a} \\\hline
    10 & \multicolumn{2}{|c|}{93^a} & \multicolumn{2}{|c|}{93^a} \\\hline
    11 & \multicolumn{2}{|c|}{105^a} & \multicolumn{2}{|c|}{105^a} \\\hline
    12 & \multicolumn{2}{|c|}{117} & \multicolumn{2}{|c|}{117} \\\hline
  \end{array}$
  \caption{Bounds for arcs in projective Hjelmslev planes over chain rings 
    of order  $q^2=9$}
  \label{tbl:q=3}
\end{table}

\begin{table}[htbp]
  \centering
  $\begin{array}{|c||c@{\;-\;}c|c@{\;-\;}c|c@{\;-\;}c|}
    \hline
    n/R
    &\multicolumn{2}{|c|}{\G_4}
    &\multicolumn{2}{|c|}{\DN_4}
    &\multicolumn{2}{|c|}{\TDN_4}
    \\\hline\hline
    0 & \multicolumn{2}{|c|}{0} & \multicolumn{2}{|c|}{0}
      & \multicolumn{2}{|c|}{0}\\\hline 
    1 & \multicolumn{2}{|c|}{1} & \multicolumn{2}{|c|}{1}
      & \multicolumn{2}{|c|}{1}\\\hline 
    2 & \multicolumn{2}{|c|}{{}^d21^a} & \multicolumn{2}{|c|}{{}^b18^e}
      & \multicolumn{2}{|c|}{18^e}\\\hline 
    3 & {}^D29&30^b & {}^D29&30^b & {}^D29&30^b \\\hline
    4 & \multicolumn{2}{|c|}{{}^f52^b} &  \multicolumn{2}{|c|}{{}^D52^b} 
      & \multicolumn{2}{|c|}{{}^D52^b} \\\hline 
%    5 & {}^f68&70^F & {}^D68&70^F & {}^D68&70^F \\\hline
      5 & \multicolumn{2}{|c|}{{}^f68^F} &  \multicolumn{2}{|c|}{{}^f68^F} 
      & \multicolumn{2}{|c|}{{}^f68^F} \\\hline 
    6 & \multicolumn{2}{|c|}{{}^f84^F} & {}^D81&83^F & {}^D81&83^F \\\hline
    7 & {}^D94&101^F &  {}^D99&101^F &  {}^D96&101^F \\\hline
    8 & \multicolumn{2}{|c|}{{}^{fC}126^b} & {}^B120&125^C 
      & {}^B120&125^C \\\hline
    9 & \multicolumn{2}{|c|}{{}^{fB}140^E} & \multicolumn{2}{|c|}{{}^B140^E}
      & \multicolumn{2}{|c|}{{}^B140^E} \\\hline 
    10 & {}^{fB}152&160^b & {}^B152&160^b & {}^B152&160^b \\\hline 
    11 & {}^B166&169^E & {}^B166&169^E & {}^B166&169^E \\\hline 
    12 & {}^B186&189^E & {}^B186&189^E & {}^B186&189^E \\\hline 
    13 & {}^B201&208^E & {}^D202&208^E & {}^D202&208^E \\\hline 
    14 & {}^D224&228^E & {}^D216&228^E & {}^D219&228^E \\\hline 
    15 & {}^{fA}236&248^b & {}^{fA}236&248^b & {}^{fA}236&248^b \\\hline 
    16 & \multicolumn{2}{|c|}{256^a} & \multicolumn{2}{|c|}{256^a}
       & \multicolumn{2}{|c|}{256^a} \\\hline
    17 & \multicolumn{2}{|c|}{276^a} & \multicolumn{2}{|c|}{276^a}
       & \multicolumn{2}{|c|}{276^a} \\\hline
    18 & \multicolumn{2}{|c|}{296^a} & \multicolumn{2}{|c|}{296^a}
       & \multicolumn{2}{|c|}{296^a} \\\hline
    19 & \multicolumn{2}{|c|}{316^a} & \multicolumn{2}{|c|}{316^a}
       & \multicolumn{2}{|c|}{316^a} \\\hline
    20 & \multicolumn{2}{|c|}{336} & \multicolumn{2}{|c|}{336}
       & \multicolumn{2}{|c|}{336} \\\hline
  \end{array}$
  \caption{Bounds for arcs in projective Hjelmslev planes over chain rings 
    of order $q^2=16$}
  \label{tbl:q=4}
\end{table}

\begin{table}[htbp]
  \centering
  $\begin{array}{|c||c@{\;-\;}c|c@{\;-\;}c|}
    \hline
    n/R
    &\multicolumn{2}{|c|}{\Z_{25}}
    &\multicolumn{2}{|c|}{\DN_5}
    \\\hline\hline
    0 & \multicolumn{2}{|c|}{0} & \multicolumn{2}{|c|}{0} \\\hline 
    1 & \multicolumn{2}{|c|}{1} & \multicolumn{2}{|c|}{1} \\\hline 
    2 & {}^g21&25^a & \multicolumn{2}{|c|}{{}^D25^a} \\\hline
    3 & {}^D40&43^F & {}^D42&43^F \\\hline
    4 & {}^D64&70^F & {}^D64&70^F \\\hline
    5 & {}^D85&102^F & {}^D90&102^F \\\hline
    6 & {}^D114&130^F & \multicolumn{2}{|c|}{{}^D130^F} \\\hline
    7 & {}^D135&156^E & {}^D152&156^E \\\hline
    8 & {}^D162&186^E & {}^D162&186^E \\\hline
    9 & {}^D186&208^E & {}^D190&208^E \\\hline
    10 & {}^D210&238^E & {}^D225&238^E \\\hline
    11 & {}^D234&265^E & {}^D250&265^E \\\hline
    12 & {}^D259&295^b & {}^D280&295^b \\\hline
    13 & {}^f310&311^E & {}^D297&311^E \\\hline
    14 & {}^D319&341^E & {}^D318&341^E \\\hline
    15 & {}^B355&367^E & {}^B355&367^E \\\hline
    16 & {}^B375&395^E & {}^B375&395^E \\\hline
    17 & {}^D400&425^E & {}^D405&425^E \\\hline
    18 & {}^B425&455^b & {}^D433&455^b \\\hline
    19 & {}^f465&466^E & {}^B455&466^E \\\hline
    20 & {}^A490&496^E & {}^A490&496^E \\\hline
    21 & {}^A515&525^E & {}^A515&525^E \\\hline
    22 & {}^A540&555^E & {}^A540&555^E \\\hline
    23 & {}^{fA}565&585^E & {}^A565&585^E \\\hline
    24 & {}^A595&615^b & {}^A595&615^b \\\hline
    25 & \multicolumn{2}{|c|}{625^a} & \multicolumn{2}{|c|}{625^a} \\\hline
    26 & \multicolumn{2}{|c|}{655^a} & \multicolumn{2}{|c|}{655^a} \\\hline
    27 & \multicolumn{2}{|c|}{685^a} & \multicolumn{2}{|c|}{685^a} \\\hline
    28 & \multicolumn{2}{|c|}{715^a} & \multicolumn{2}{|c|}{715^a} \\\hline
    29 & \multicolumn{2}{|c|}{745^a} & \multicolumn{2}{|c|}{745^a} \\\hline
    30 & \multicolumn{2}{|c|}{775^a} & \multicolumn{2}{|c|}{775^a} \\\hline
  \end{array}$
  \caption{Bounds for arcs in projective Hjelmslev planes over chain rings 
    of order $q^2=25$}
  \label{tbl:q=5}
\end{table}

\noindent
\begin{footnotesize}
Notes to the tables: 

\noindent
${}^a$ \cite{it:deadfin}
${}^b$ \cite{it:kent}
${}^c$ \cite{ivan-silvia04}
${}^d$ \cite{it:ovals}
${}^e$ \cite{kiermaier06}
${}^f$ \cite{kiermaier-kohnert07}
${}^g$ \cite{kiermaier-koch09}
${}^h$ \cite{imt:two_exact_values}
${}^i$ \cite{honold-kiermaier12}

${}^A$ Subsection~\ref{ssec:genconst}
${}^B$ Subsection~\ref{ssec:speconst}
${}^C$ Subsection~\ref{ssec:dualconst}
${}^D$ Subsection~\ref{ssec:bayconst}

${}^E$ Subsection~\ref{ssec:newbound}
${}^F$ Subsection~\ref{ssec:spebounds}

% $^a$ Subsection~\ref{subsec:smallsize}
% $^b$ Theorem~\ref{thm:triangleset}
% $^c$ Corollary~\ref{cor:q2-q}
% $^d$ Theorem~\ref{thml:q=3:n=7}
% $^e$ Remark~\ref{rmkl:q=3:n=8}

% \noindent
% $^A$\ Theorem~\ref{thm:40_5},
% $^B$\ Theorem~\ref{thm:39_5} 
% $^C$\ Theorem~\ref{thm:71_8}
% $^?$\ Theorem~\ref{thm:n=3}
\end{footnotesize}

\subsection{Some arcs of small size}
\label{subsec:smallsize}

\noindent
A $(30,4)$-arc in $\PHG(2,R)$, $R=\mathbb{Z}_9$

{\small \[
\begin{array}{rrrrr rrrrr rrrrr rrrrr rrrrr rrrrr}
 0& 0& 0& 3& 3& 6& 3& 6& 6& 1& 1& 1& 1& 1& 1& 1& 1& 1& 1& 1& 1& 1& 1& 1& 1& 1& 1& 0& 3& 6\\
 0& 1& 1& 0& 1& 1& 1& 1& 6& 0& 1& 3& 7& 8& 8& 4& 5& 6& 6& 7& 8& 0& 1& 3& 4& 5& 5& 1& 1& 1\\
 1& 0& 4& 1& 6& 7& 4& 6& 1& 3& 8& 7& 3& 4& 5& 2& 2& 3& 7& 6& 1& 6& 6& 4& 8& 1& 5& 2& 8& 5    
\end{array}\]}

\noindent
A $(30,4)$-arc in $\PHG(2,R)$, $R=F_3[X]/(X^2)$

\noindent
A $(39,3)$-arc in $\PHG(2,R)$, $R=\mathbb{Z}_9$
{\small
\[\begin{array}{rrrrr rrrrr rrrrr rrrrr}
1& 1& 1& 1& 1& 1& 1& 1& 1& 3& 6& 6& 6& 0& 1& 1& 1& 1& 1& 1\\
0& 1& 3& 5& 5& 6& 7& 7& 8& 1& 0& 1& 1& 1& 1& 2& 6& 6& 6& 7\\
0& 3& 4& 0& 2& 5& 2& 7& 1& 0& 1& 2& 7& 0& 1& 2& 2& 3& 4& 5\\ 
\end{array}\]}
{\small
\[\begin{array}{rrrrr rrrrr rrrrr rrrr}
1& 1& 1& 3& 3& 6& 0& 1& 1& 1& 1& 1& 1& 1& 1& 1& 3& 3& 3 \\
7& 8& 8& 1& 1& 6& 1& 1& 1& 3& 3& 3& 4& 5& 5& 8& 1& 1& 6 \\
6& 3& 7& 1& 5& 1& 5& 5& 7& 0& 1& 5& 6& 4& 8& 6& 3& 7& 1 
\end{array}\]}

\noindent
A $(28,3)$-arc in $\PHG(2,R)$, $R=\GR(4^2,2^2)$

{\small
\[\begin{array}{rrrrrrrrrrrrrr}
0 & 0 & 8 & 8 & 8 & 1 & 1 & 1 & 1 & 1 & 1 & 1 & 1 & 1 \\
0 & 1 & 8 & 1 & 1 & 0 & 0 & 4 & 4 & 8 &12 &12 & 1 & 5 \\
1 & 4 & 1 & 5 &11 & 4 & 8 & 5 & 6 & 0 &10 &11 & 0 & 4  
\end{array}\]
\[\begin{array}{rrrrrrrrrrrrrr}
1 & 1 & 1 & 1 & 1 & 1 & 1 & 1 & 1 & 1 & 1 & 1 & 1 & 2 \\
5 & 9 &13 &13 & 2 & 2 & 6 & 6 &10 & 3 & 7 & 7 &11 & 1 \\
2 &15 & 8 & 1 & 1 &15 & 2 &11 &15 &11 & 4 & 5 & 6 & 0 
\end{array}\]}

\noindent
A $(68,5)$-arc in $\PHG(2,R)$, $R=\GR(4^2,2^2)$

{\footnotesize
\[\begin{array}{rrrrr rrrrr rrrrr rrrrr rrr}
 0& 0& 1& 1& 1& 1& 2&10& 8& 1& 1& 1& 1& 1& 1&10& 0& 0& 1& 1& 1& 1& 2\\ 
 1& 2& 4&12& 1&11& 2& 1& 1& 0& 0& 8& 8&12&14& 1& 1&10& 4&12& 9& 3&10\\
 1& 1& 9&11&11& 3& 1&11& 5& 9& 2& 2& 3&13&13&15&11& 1&11& 9& 3&11& 1\\
\end{array}\]

\[\begin{array}{rrrrr rrrrr rrrrr rrrrr rrr}
10& 1& 1& 1& 1& 1& 1& 2& 2& 0& 1 & 1& 1& 1& 1& 1& 2& 1& 1& 1& 1& 1 & 1\\
 1& 4& 4& 5&13&13& 7& 1& 1& 1& 0 & 0& 4& 8& 8& 6& 1& 1& 9& 9& 2&10 & 7\\
 1& 8&10& 5& 9&13& 9& 4&14&15& 0 &11& 5& 0& 1& 5& 5& 8& 5& 2&13& 7 & 0
\end{array}\]

\[\begin{array}{rrrrr rrrrr rrrrr rrrrr rr}
 1& 1& 1& 1& 1& 1& 1& 1& 1& 1& 8& 8& 1& 1& 1& 1& 1& 1& 1& 1& 1& 1\\
 7&11& 0& 8& 1& 5& 5& 3& 3&11& 1& 1&12&12& 5&13&13& 7& 1& 9& 3&11\\ 
 8& 7& 7&13& 5& 0& 8& 2& 7& 8&12& 6& 0& 2& 7& 3&15& 3& 4& 4& 4& 4 
\end{array}\]}

\subsection{Open Problems}

\begin{enumerate}[(1)]
\item What is the maximal size of a $(k,2)$-arc in $\PHG(R_R^2)$ for
  geometries over chain rings of length 2 and characteristic different
  from 4?  For chain ring of characteristic 2 the upper bound is
  $q^2+q$ and for chain rings of odd characteristic it is
  $q^2$. Numerical data suggests that for large $q$ these bounds are
  not met. For instance for chain rings of size 16 the largest size of
  such an arc is 18 (and not 20). Find a general construction for
  $(k,2)$-arcs in $\PHG(2,R)$, $R$ a chain ring of length 2,
  $\abs{R}=q^2$ and $\Char R\ne4$, which gives arcs of size $\approx
  cq^2$ for some constant $c\ge1/2$ (preferably $c$ should be close to
  1).

\item It is known that there exists a $(2^{2r}+2^r+1,2)$-arc in the
projective Hjelmslev plane over the Galois ring $R=\GR(2^{2r},2^2)$,
for all $r\ge1$. We call these arcs hyperovals since they have jjust two
intersection numbers: 0 and 2. They were constructed via Witt vectors
by taking the Techm\"uller set of points corresponding to the 
(unique)subgroup of units of order $2^{3r}-1$ in the multiplicative  group
of the extension of $R$ of order 3. It is not known whether these arcs 
are unique up to collineation.

\item The upper bound for the size of an arc
with $n=q^2-q$ (in an arbitrary projective Hjelmslev plane) is
$(q-1)^2(q^2+q+1)$. It is easily proved that each point class
has to have multiplicity $(q-1)^2$. The missing $2q-1$ points should 
form a blocking set in $\AG(2,q)$. Construct such a class of arcs
or prove its nonexistence.

\item Find an example of chain rings $R$ and $S$
with $\abs{R}=\abs{S}=q^2$, $R/\rad R\cong\ S/\rad S\cong\mathbb{F}_q$ and the 
same characteristic for which $m_n(R_R^3)\ne m_n(S_S^3)$ for some $n$.
The same question for rings that are not commutative
and of the same cardinality, length and characteristic.
More generally, can the maximal problem distinguish all pairs of chain rings.

\item Prove (or disprove) that for $0\leq n\leq q^2+q$ all
  maximal $(k,n)$-arcs in $\PHG(2,R)$ are projective.
  
  \item (A. A. Nechaev) Are chain rings of size $q^2$ geometrically distinguishable?
  In other words, let $R$ and $R'$ be non-isomorphic chain rings with
  $|R|=|R'|=q^2$, $R/\rad R\cong R'/\rad R'\cong\mathbb{F}_q$, $\Char R=\Char R'$.
  Is it true that $m_n(R)=m_n(R')$ for all $n$? 
\end{enumerate}

\noindent
{\it Acknowledgments.}\
The research of the first author was supported by the National Natural
Science Foundation of China under Grant 61571006. 
%The research of the second author was supported by ..... . 
The research of the third author was supported 
by the Bulgarian National Science Foundation under Contract KP-06-N72/6-2023.

%\noindent
%{\it Acknowledgements.}\
%The research of the first author was supported by the Research Fund of
%Sofia University under contract No. 80-10-51/17.04.2018.
%The research of the second author was supported by
%the project "Finite geometries, coding theory and cryptography" 
%between the Research Foundation - Flanders (FWO) and
%Bthe Bulgarian Academy of Sciences. 

%\section{Conclusions and remarks}
%Here some conclusions and remarks.

\end{document}